\documentclass[12pt]{amsart}
\usepackage{amsmath,amssymb,amsfonts,amscd,graphicx,xypic}
\usepackage{relsize}

\usepackage{changes}

\newcommand{\slantedcup}{\mathbin{\rotatebox[origin=c]{50}{$\cup$}}}

\newcommand{\Picture}[1]{
\begin{minipage}{.8in}
\includegraphics[scale=.15]{#1}
\end{minipage}
}

\newcommand{\HopfLink}{\ensuremath{[\!\Picture{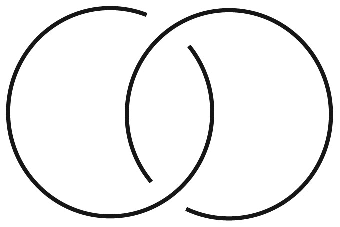}}}

\newtheorem{lemma}{Lemma}[section]

\newtheorem{proposition}[lemma]{Proposition}
\newtheorem{corollary}[lemma]{Corollary}
\theoremstyle{definition}
\newtheorem{definition}[lemma]{Definition}

\theoremstyle{remark}
\newtheorem{remark}[lemma]{Remark}
\numberwithin{equation}{section}
\theoremstyle{plain}
\newtheorem{thm}{Theorem}
\newtheorem{question}{Question}[section]
\newtheorem{cor}{Corollary}[section]
\newtheorem{prop}[cor]{Proposition}
\theoremstyle{definition}
\newtheorem{defi}[cor]{Definition}
\newtheorem{condition}[lemma]{Condition}

\newcommand{\co}{\colon\thinspace}

\numberwithin{figure}{section}

\begin{document}

\title[Engel relations in $4$-manifold topology]{Engel relations in ${\mathbf 4}$-manifold topology}
\author{Michael Freedman and Vyacheslav Krushkal}

\address{Microsoft Station Q, University of California, Santa Barbara, CA 93106-6105}
\email{michaelf\char 64 microsoft.com}
\address{Department of Mathematics, University of Virginia, Charlottesville, VA 22904}
\email{krushkal\char 64 virginia.edu}

\begin{abstract}
We give two applications of the $2$-Engel relation, classically studied in finite and Lie groups, to the $4$-dimensional topological surgery conjecture.
The \mbox{A-B} slice problem, a reformulation of the surgery conjecture for free groups, is 
shown to admit a homotopy solution.
We also exhibit a new collection of universal surgery problems,
defined using ramifications of  homotopically trivial links. 
More generally we show how $n$-Engel relations arise from higher order double points of surfaces in $4$-space.
\end{abstract}

\dedicatory{Dedicated to Andrew Casson}

\maketitle

\section{Introduction} \label{introduction}
Forty years ago Andrew Casson taught us \cite{Casson} that singularities of surfaces and the fundamental group of their complements are intimately related. 
We study a classical group relation, $2$-Engel, and the corresponding surface singularities. The results include two surprises (to us) regarding topological surgery.
What direction they point is presently unknown. They might later be seen as: a step in proving the full surgery conjecture, or contrariwise as pointing toward a surgery obstruction, or possibly as mere curiosities. The purpose of this paper is to explain these surprises and reconsider fundamental conjectures and constructions in this new light.

Topological surgery is known to work in dimension $4$ for a class of ``good'' fundamental groups. Originally this was established in the simply-connected setting by the first author in \cite{F0}. It has since been shown \cite{F1} that elementary amenable groups, and more recently \cite{FT, KQ} the groups of subexponential growth are good in this sense. The validity of surgery for arbitrary fundamental groups remains a central open problem. Surgery may be reduced to a collection of {\it universal} problems \cite{CF, F1} with free fundamental groups, therefore the validity of surgery for (non-abelian) free groups is the key open question. It has been reformulated \cite{F2, F3} in terms of the {\it A-B slice} problem for a family of links, the ``generalized Borromean rings''.

We give applications of the group-theoretic {\it $2$-Engel relation} both to the A-B slice problem and to construction of model surgery problems. The study of the universal relation, stating that all $3$-fold commutators of the form $[[y,x],x]]$ are trivial in a group $G$, dates back to the work of Burnside \cite{Burnside}. It is easily seen to be equivalent to the relation that every element $x$  in $G$ commutes with all of its conjugates $x^y$. A restricted version of this relation is familiar in low-dimensional topology: when applied to a set of preferred normal generators $x$ of a group $G$, it is a defining relation of the {\it Milnor group} $MG$, see \cite{M} and section \ref{Engel section} below. The results of imposing the relation in these two settings turn out to be quite different: the free Milnor group on $n$ generators, $MF_n$, is nilpotent of class $n$. On the other hand, the free group $F_n$ modulo the universal $2$-Engel relation is nilpotent of class $3$, independent of $n$ (see section \ref{Engel section}). This is the property of the Engel relation that we exploit in our applications.

To formulate our first result, we briefly recall the A-B slice problem (a detailed discussion is given in section \ref{ABslice section}). Surgery for free groups predicts the existence of topological $4$-manifolds $M$ which are homotopy equivalent to a wedge of circles and whose boundary is the zero-framed surgery on a Whitehead double of $L$, for each $L$ in the collection of generalized Borromean rings.  These links (GBRs) are obtained from the Borromean Rings by ramification and Bing doubling. Following \cite{F2, F3} consider the resulting free group action on the end-point compactification of the universal cover $\widetilde M$, which is homeomorphic to the $4$-ball. Choosing a fundamental domain for this action, one is led to the notion of a {\it decomposition} $D^4=A \cup B$ of the $4$-ball into two codimension zero smooth submanifolds, extending the standard genus one Heegaard decomposition of  $\partial D^4$. Given an $n$-component GBR $L$, the existence of the free group action is then equivalent to the existence of $n$ decompositions $D^4=A_i\cup B_i$ and a disjoint embedding problem for these $2n$ submanifolds into $D^4$, 
with the boundary condition given by the link $L$ and its parallel copy. If this embedding problem has a solution, the link $L$ is called A-B slice.

Considering handle decompositions of the submanifolds, one gets a pair of links, 
which we call a ``stabilization'', corresponding to the $1$- and $2$-handles.  The embedding question can then be reformulated \cite{FL} as a {\it relative-slice problem} for a certain collection of link pairs corresponding to a GBR $L$.  
A key feature of the GBRs is that they are homotopically essential in the sense of Milnor. 
Therefore it is a natural question whether there is a link-homotopy obstruction in the A-B slice problem, in other words whether  the 
relevant relative-slice problems do not even admit a link-homotopy solution. The evidence thus far has pointed to an affirmative answer: partial obstructions of this type have been found for many families of decompositions; see for example  \cite{FL}, \cite{K1}. 
Surprisingly, here we construct the first examples of decompositions giving rise to a homotopy solution to the A-B slice problem. 

\begin{thm} \label{ABslice theorem} \sl 
The generalized Borromean rings, a collection of links forming universal surgery problems, are homotopy A-B slice. \sl
\end{thm}

We present two possible notions of a ``homotopy solution'', one in the sense of link-homotopy, and a stronger one in terms of disjoint homotopy of $2$-handles; see
definitions \ref{link homotopy ABslice definition}, \ref{homotopy AB slice} in section \ref{ABslice section}. The theorem is true for both notions.

The action of the free group on $D^4$ by covering transformations is encoded in the requirement that the disjoint embeddings
of the $A_i, B_i$ in $D^4$ are {\it standard}, in other words isotopic to the original embeddings corresponding to the given decompositions $D^4=A_i\cup B_i$.
It was observed in \cite{K} that there exist solutions to the embedding problem
if this requirement is omitted. (However the existence of a solution without the equivariant feature does not have a direct implication for surgery.) Our proof of theorem \ref{ABslice theorem} satisfies 
the homotopy analogue of the standard embedding requirement, see definition \ref{homotopy AB slice} and the proof of theorem \ref{ABslice theorem} in section \ref{homotopy solution}.

One way to view theorem \ref{ABslice theorem} is as evidence towards the validity of the surgery conjecture. There is a well established hierarchy of $2$-complexes, defined in terms of gropes and capped gropes (cf. \cite{FQ}), extrapolating between disjoint surfaces and disjoint embedded disks. It seems possible that a homotopy solution to the AB slice problem may be further improved using group-theoretic methods. For example $n$-Engel relations, $n>2$, are candidates for such an approach, however these higher relations are not as well understood algebraically.  We refer the reader to \cite{Traustason} for a recent survey of the subject. It is an open question whether a homotopy solution may be improved to a stage that would imply an actual embedded solution to the AB slice problem. To assist the reader who would like to solve this problem we discuss in the Appendix how $n$-Engel relations relate to higher order self-intersections of a disk.

It has been shown in  \cite{FT2} that Whitehead doubles of (homotopically trivial)$^+$ links (a class of links just slightly smaller than homotopically trivial links) are topologically slice. Therefore (even a strong version of)  the AB slice problem has a solution for (homotopically trivial)$^+$ links. Viewing our present work in the context of the relative-slice problem, discussed in section \ref{ABslice section}, for each GBR we find a stabilization so that the resulting link is homotopically trivial. As remarked above, starting with a (homotopically trivial)$^+$ link a stabilization may be found so that the result is slice. There is gap corresponding to the $+$ assumption, but it is an interesting question whether the two stabilizations may be combined to give a solution.

Overall, the key open problem is to determine whether there still is an obstruction to the AB slice problem in terms of nilpotent invariants of links, specifically Milnor's $\mu$-invariants.
Our results here suggest $\bar\mu$ invariants with repeating coefficients  are not necessarily more fragile than non-repeating and should play a role in a surgery obstruction - if in fact there is one. The upshot is that the repeating/non-repeating dichotomy now seems false.
  Of course such an obstruction would give a counterexample to surgery for free groups. Conversely, as discussed above the ability to ``improve'' a homotopy solution could lead to the resolution of the surgery conjecture in the affirmative.  
An axiomatic framework in terms of {\it topological arbiters} for an obstruction in the AB slice program has been introduced in \cite{FK}.
Since our theorem \ref{ABslice theorem} constructs a solution up to homotopy, there is no topological arbiter satisfying an extended ``Bing doubling axiom'' \cite{FK} defined in terms of  $\bar\mu$-invariants with non-repeating coefficients. Since the method of proof of theorem \ref{ABslice theorem} does not extend to the relevant stabilized link together with parallel copies of its components, $\mu$-invariants with repeating coefficients remain a candidate for a surgery obstruction. But if such an obstruction  is produced it will {\it not} be to homotopy solutions but actual standardly embedded solutions to the strong AB slice problem.

It is interesting to compare the complexity of the homotopy solution to the AB slice problem constructed in theorem \ref{ABslice theorem} with the current state of knowledge about general decompositions $D^4=A\cup B$. 
A recent paper  \cite{K1} gave a thorough analysis of the decompositions of the $4$-ball where $A$ has two $2$-handles and {\it one} $1$-handle. The answer is quite subtle and the analysis relies on delicate\footnote{It is important that a certain system of cubic equations has no integral solutions whereas it manifestly has a solution over ${\mathbb Z}[1/4]$.} calculations in commutator calculus. 
In the relevant decomposition $D^4=A\cup B$ used for the Borromean rings in the proof of theorem \ref{ABslice theorem} (see section \ref{homotopy solution}), the side $A$ has two $2$-handles and 36 $1$-handles. (The $B$-side has a handle decomposition with the number of $1$- and $2$-handles reversed.) 
It seems likely that a novel algebraic structure will be needed to gain further insight into the problem.

In section \ref{s-cobordism section} we describe a slicing problem for a link in a $4$-manifold, the ``Round Handle Problem'' (RHP), where the existence of a solution depends not just on surgery but also on the $5$-dimensional s-cobordism conjecture. At first sight this problem appears similar to the relative-slice formulation of the AB slice problem for GBRs, however
the proof of theorem \ref{ABslice theorem} does not extend to this setting. This suggests a subtle distinction between the two problems, with the possibility that a link-homotopy obstruction is still possible to the combination of surgery and s-cobordism conjectures. That is, non-repeating $\bar\mu$ invariants might still form the basis of an obstruction to the RHP for the Generalized Borromean Rings.

Another application of the $2$-Engel relation yields a new set of universal surgery problems. The ``usual'' model surgery kernels \cite{FQ} are given by $S^2\vee S^2$-like capped gropes. They are universal in the sense that if solvable they imply solutions to all $4$-dimensional surgery problems with the vanishing Wall obstruction, see section \ref{model section}.
There is a corresponding collection of slicing problems for links {\it $\{\!$Wh(Bing(Hopf))$\!\}$} (where the slice 
complement in the $4$-ball is required to have free fundamental group, generated by the meridians). The links in question are Whitehead doubles of the generalized Borromean rings mentioned above and discussed further in section \ref{model section}. We introduce a new collection of universal slicing problems:
\begin{thm} \label{new models theorem} \sl
There is a family of links $\{ K\}$  for which the problem of constructing free slices constitutes a universal problem, where each $K\in \{ K\}$ is of the form: $${\text{D(Ram(h-triv)),}}$$ a generalized (genus one)  double of a ramified homotopically trivial link.
\end{thm}

The ``double'' in this statement is a generalization of the notion of a Whitehead double of a link, introduced in section \ref{general double section}. The key (and surprising) feature of this new collection of links is that they are defined starting from {\it homotopically trivial} links, see remark \ref{universal theorem remark}.

The organization of the paper is as follows.
Section \ref{Engel section} discusses the $2$-Engel relation.
A proof is given that $2$-Engel groups are $3$-nilpotent. We wish to note that Peter Teichner was already aware of this fact in the 1990s.
In section \ref{ABslice section} we recall the formulation of the A-B slice problem with a particular focus on the notions of a homotopy solution and a standard embedding, which are important for theorem \ref{ABslice theorem}.
The proof of theorem \ref{ABslice theorem} is given in section \ref{homotopy solution}.
Section \ref{s-cobordism section} formulates the Round Handle Problem, providing a comparison of our results with the setting of the s-cobordism theorem.
In section \ref{weak homotopy section} we define a very strong equivalence relation on links in $S^3$, called weak link homotopy (WLH$_0$). Using this relation, 
new model surgery problems are constructed in section \ref{model section}.
In the Appendix we show how $n$-Engel relations correspond to higher order intersections of disks.

\section{The $2$-Engel relation} \label{Engel section}

The Milnor group provides a convenient setting for the analysis of the $2$-Engel relation and for the main results of the paper.
We start by briefly reviewing the Milnor group and link homotopy in section \ref{Milnor group section}; the reader is referred to \cite{M} for a detailed introduction. Section \ref{Engel subsection} presents the $2$-Engel relation and shows that $2$-Engel groups are $3$-nilpotent. A geometric realization of this relation, weak homotopy of links, is discussed in section \ref{weak homotopy section}.

\subsection{The Milnor group.} \label{Milnor group section}

\begin{defi}
Let $G$ be a group normally generated by a fixed finite collection of elements $g_1,\ldots, g_n$. The 
{\it Milnor group} of $G$, defined with respect to the given normal generating set $\{ g_i\}$,
is given by
\begin{equation} \label{eq:Milnor group}
MG := G \, /\, \langle\! \langle \, [g_i, g_i^y] \;\; i=1,\ldots, n,\;\, y\in G \rangle\!\rangle.
\end{equation}
\end{defi}

The Milnor group $MG$ is generated by $g_1, \ldots, g_n$. Moreover, it is a finitely presented nilpotent group of class $\leq n$, see \cite{M}. We remind the reader that a group $G$ is said to be nilpotent of class $n$ if the $(n+1)^{\rm st}$ stage $G^{n+1}$ of its lower central series vanishes, see section \ref{Engel subsection}.

Given an $n$-component link $L$ in $S^3$, let $G$ denote ${\pi}_1(S^3\smallsetminus L)$.
Consider meridians
$g_i$ to the components $l_i$ of $L$: $g_i$ is an
element of $G$ obtained by following a path $\alpha_i$ in $S^3\smallsetminus L$
from the basepoint to the boundary of a regular neighborhood of $L$, followed
by a small circle (a fiber of the circle normal bundle) linking $l_i$, then followed
by ${\alpha}_i^{-1}$. $G$ is normally generated by the elements $g_1, \ldots, g_n$. Then 
$MG$, defined with respect to the meridians, is called the Milnor group $ML$ of the link $L$.

Denoting by $F_{g_1,\ldots, g_n}$ the free group generated by the $\{g_i\}$, $i=1,\ldots, n$, consider the Magnus
expansion
\begin{equation}\label{Magnus}
M\co F_{g_1,\ldots, g_n}\longrightarrow {\mathbb Z}[\! [x_1,\ldots, x_n]\! ]
\end{equation}
into the ring
of formal power series in non-commuting variables $\{ x_i\}$, defined by $$M(g_i)=1+x_i, \;\,
M(g_i^{-1})=1-x_i+x_i^2-x_i^3\pm\ldots$$ The Magnus expansion induces a homomorphism
\begin{equation} \label{MagnusMilnor}
MF_{g_1,\ldots, g_n}\longrightarrow R_{x_1,\ldots,x_n},
\end{equation}
into the quotient $R_{x_1,\ldots,x_n}$ of ${\mathbb Z}[\! [x_1,\ldots, x_n]\!  ]$ by the ideal generated by all monomials
$x_{i_1}\cdots x_{i_k}$ with some index occuring at least twice. The homomorphism (\ref{MagnusMilnor}) is well-defined and injective \cite{M}.
Using the Magnus expansion it is not difficult to see that the Milnor group $MF_n$ of the free group $F_n$ on $n$ generators is nilpotent of class precisely equal to $n$.

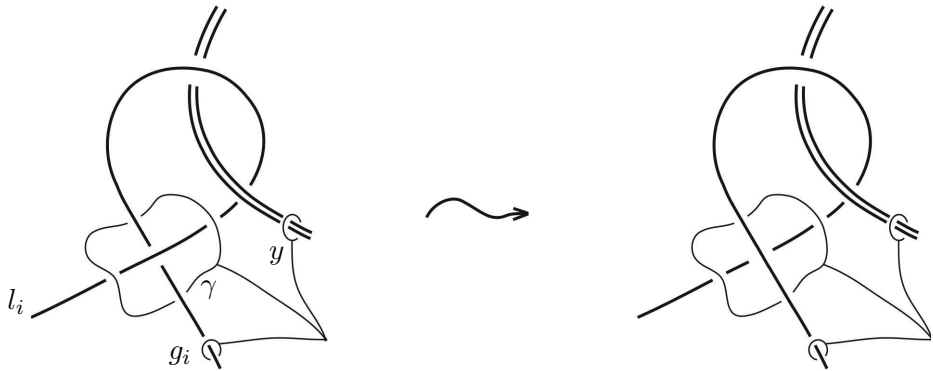
\begin{figure}[ht]
\includegraphics[height=5cm]{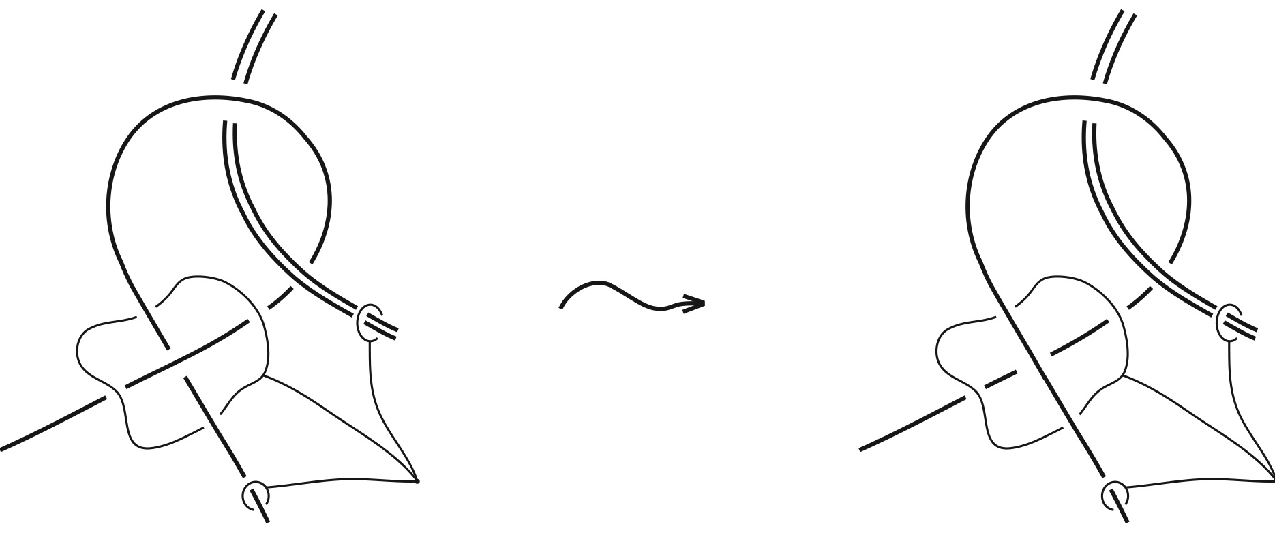}
{\small
\put(-289,6){$g_i$}
\put(-277,31){$\gamma$}
\put(-251,44){$y$}
\put(-350,25){$l_i$}
}
 \caption{An illustration of a non-generic ``crossing-time'' during a link homotopy. The based curve ${\gamma}$ in the link complement, corresponding to the defining relation $[g_i, g_i^y]$ of the Milnor group, becomes trivial after a self-intersection of the component $l_i$.}
\label{link homotopy figure}
\end{figure}

The defining relations of the Milnor group (\ref{eq:Milnor group}) are well suited for studying 
links $L$ in $S^3$ up to link homotopy. Recall that two links are {\it link-homotopic} if they
are connected by a $1$-parameter family of link maps where different components
stay disjoint for all values of the parameter, see figure \ref{link homotopy figure}.
If $L$, $L'$ are link-homotopic then 
their Milnor groups $ML$, $ML'$ are isomorphic, and moreover an $n$-component link $L$ is
homotopically trivial (link-homotopic to the $n$-component unlink) if and only if $ML$ is isomorphic to the free Milnor group $MF_{m_1,\ldots, m_n}$.

The Milnor group is also useful for studying surfaces $\Sigma$ in the $4$-ball where the components are disjoint but may have self-intersections. In this context the Clifford tori (cf. \cite[p. 32]{FQ})   linking
the double points in $D^4$ give rise to the relations (\ref{eq:Milnor group}) in $M{\pi}_1(D^4\smallsetminus {\Sigma})$.
Link homotopy theory may be interpreted as the study of links up to singular concordance
(links $L\subset S^3\times\{0\}$, $L'\subset S^3\times\{1\}$ bounding disjoint maps of annuli into $S^3\times[0,1]$).
In particular, a link $L$ is homotopically trivial if and only if its components bound disjoint immersed disks $\Delta$ in $D^4$, and in this case $M{\pi}_1(S^3\smallsetminus L)\cong M{\pi}_1(D^4\smallsetminus {\Delta})\cong MF_n$.

\subsection{2-Engel groups}  \label{Engel subsection}
We start off by fixing the notation. The lower central series of a group $G$ is defined inductively by $G^1=G, G^n=[G^{n-1}, G]$.
Given $g_1,\ldots, g_n\in G$, the commutator $[[\ldots [g_1, g_2],\ldots, g_{n-1}], g_n]$ will be concisely denoted $[g_1, g_2, \ldots, g_n]$.

The main focus of this section is on $2$-Engel groups, that is groups satisfying the universal relation $[y,x,x]=1$, or equivalently $[x,x^y]=1$. Unlike the setting of the Milnor group (\ref{eq:Milnor group}), this relation is {\it universal} in the sense that it holds for {\it all} elements $x, y$ of a $2$-Engel group. 

In reference \cite{Burnside} (which is at the foundation of the subject of Engel groups)  W.~Burnside showed any elements of a $2$-Engel group $G$ satisfy the identities 
$$
[x,y,z]\; \; =\; \; [y,z,x], \; \; [x,y,z]^3\;\; = \; \; 1.
$$
In a later paper \cite{Hopkins} C. Hopkins showed that $2$-Engel groups $G$ have nilpotency class $\leq 3$, that is $G^4=\{1\}$. (Also see \cite{Levi}.)
We give a proof of this result below in the context of the Milnor group, to establish a reference point for geometric applications in later sections. Corollary 
\ref{Engel corollary} summarizes the relevant facts.

It is interesting to note that the $2$-Engel relation is functorial, that is any group homomorphism $G\longrightarrow H$ induces a homomorphism of groups modulo the relation. This contrasts the Milnor group setting (\ref{eq:Milnor group}): only homomorphisms taking chosen generators of $G$ to the chosen generators of $H$ are guaranteed to induce a homomorphism of the Milnor groups, $MG\longrightarrow MH$, defined with respect to these generators. It follows that the quotient needed for making the theory functorial necessarily kills most of non-abelian information:

\begin{lemma} \label{Engel nilpotent} \sl
Any $2$-Engel group is nilpotent of class $\leq 3$.
\end{lemma}

We will use the following basic result about Milnor groups. Given a group $G$ normally generated by $g_1, \ldots, g_n$, consider ``basic commutators'' $[g_{i_1}, \ldots, g_{i_m}]$.

\begin{proposition} \label{repeated proposition} \sl
Any basic commutator
$[g_{i_1}, \ldots, g_{i_m}]$ where at least two of the indices coincide, $i_j=i_k$ for some $j\neq k$, is trivial in the Milnor group $MG$.
\end{proposition}

One way to prove this fact is to use the Magnus expansion. Every monomial (other than $1$) in the expansion 
of $[g_{i_1}, \ldots, g_{i_m}]$ has a variable $x_{i_j}=x_{i_k}$ occurring at least twice. Since the Magnus expansion  (\ref{MagnusMilnor}) is injective it follows that
$[g_{i_1}, \ldots, g_{i_m}]=1\in MF_{g_1, \ldots, g_n}$.

One can also use the commutator identities (\ref{product id}) below  to show directly that $[g_{i_1}, \ldots, g_{i_m}]$ with repeating labels is a product of
the defining relations (\ref{eq:Milnor group}) of the Milnor group.

{\it Proof of lemma} \ref{Engel nilpotent}. 
We will use the Hall-Witt identity (\ref{Hall-Witt eq}) and basic commutator identities (\ref{product id}), cf. \cite[Theorem 5.1]{MKS}.
\begin{equation} \label{Hall-Witt eq}
[x,y,z^x]\cdot [z,x,y^z]\cdot [y,z,x^y]=1,
\end{equation}
\begin{equation} \label{product id}
[x,yz]\, =\, [x,z]\; [x,y]^z, \; \;\,[xz,y]\, =\, [x,y]^z\; [z,y], \; \;\, [x^{-1},y]=[y,x]^{x^{-1}}.
\end{equation}
It suffices to show that the free group $F_n=F_{g_1, \ldots, g_n}$ modulo the $2$-Engel relation is nilpotent of class $3$.
This quotient factors through the Milnor group $MF_n$. Using the identities (\ref{product id}) and proposition \ref{repeated proposition}, $(MF_n)^4$ is seen to be normally generated by 
commutators $[g_{i_1}, \ldots, g_{i_4}]$ with non-repeating indices.
It suffices to show that $F_4$,  the free group on $4$ generators is nilpotent of class $3$, after dividing out by the $2$-Engel relations. 

Denote the generators of $F_4$ by $x,y,z,w$.
First we focus on $3$-fold commutators. 
Denoting by ``$\equiv$'' the equivalence up to the $2$-Engel relation and reserving ``$=$'' for equality in the free Milnor group $MF_{4}$, one has, by definition:
\begin{equation} \label{first eq}  1\equiv [z, xy, xy].
\end{equation}
Expanding this commutator according to (\ref{product id}) yields a product of four terms (where the conjugations are omitted, for a reason discussed below) in (\ref{second eq}). The second equality follows from proposition \ref{repeated proposition}:
\begin{equation} \label{second eq}
[z, xy, xy]\; =\; [z,x,x]\cdot [z,x,y] \cdot [z,y,x] \cdot [z,y,y]\; = \; [z,x,y] \cdot [z,y,x].
\end{equation}
It is a basic fact that conjugation as in (\ref{Hall-Witt eq}), (\ref{product id}) may be disregarded in Milnor group calculations of this type, as they contribute corrections in the kernel ($F_4\longrightarrow MF_4$). One way to see this is to consider the Magnus expansion (\ref{MagnusMilnor}). The effect of conjugation is an introduction of higher order terms. Each higher order monomial that comes up in applications of the identities 
(\ref{product id}) to (\ref{second eq}) has repeated indices, so is trivial in the target ring $R$ of the Magnus expansion. Since the expansion (\ref{MagnusMilnor}) is injective, it follows that conjugation resulting from the commutator identities (\ref{product id}) does not change the terms appearing in (\ref{second eq}). 

It follows from (\ref{first eq}), (\ref{second eq}) that 
\begin{equation} \label{another equation}
[z,x,y]\; \equiv \; [z,y,x]^{-1} \; =\; [y,z,x].
\end{equation}
Similarly,
\begin{equation} \label{similarly}
1\equiv [x, yz, yz]=[x,y,z]\cdot [x,z,y], \; {\rm so} \;  [x,y,z]\; \equiv \; [x,z,y]^{-1} \; =\; [z,x,y].
\end{equation}

Then the Hall-Witt identity (where conjugation is again irrelevant) implies:
\begin{equation} \label{torsion}
1\; = \; [x,y,z]\cdot [z,x,y] \cdot [y,z,x] \; \equiv\; [x,y,z]^3.
\end{equation}
Using (\ref{product id}) and disregarding conjugation in the Milnor group as above, it follows that $4$-fold commutators are also 
of order $3$:
\begin{equation} \label{another equation2}
[x,y,z,w]^3\; = \; [[x,y,z]^3, w] \; \equiv \; 1.
\end{equation}

The first equality is obtained from two applications of the middle identity of (\ref{product id}).
Next we show that $4$-fold commutators are also of order $4$. The  Hall-Witt identity (\ref{Hall-Witt eq}) 
(applied to $[x,y], z, w$) implies in the Milnor group $MF_4$:
$$[x,y, z, w]\cdot [w,[x,y],z]\cdot [z,w,[x,y]]=1.$$
Now interchange the order of the terms $w,[x,y]$ in the second commutator above (this inverts the term) and expand the last commutator using the Hall-Witt identity (applied to $x,y$ and $[z,w]$):
\begin{equation} \label{commutators}
[x,y,z,w]\cdot [x,y,w,z]^{-1}\cdot [z,w,x,y]\cdot [z,w,y,x]^{-1}=1.
\end{equation}

Again applying (\ref{product id}) twice and dropping terms with repeated letters, the first two terms of (\ref{commutators}) reduce to 
$[x,y,zw,zw]\equiv 1$, a 2-Engel relation, and therefore cancel in the quotient. Similarly, it follows from 
$[z,w,xy,xy]\equiv 1$ that the last two terms are equal.
Denoting $P:=[x,y,z,w]$, $Q:=[z,w,x,y]$, the equation (\ref{commutators}) then asserts: 
$P^2\, Q^2=1$. 

It follows from the 2-Engel relation$[y,xz,xz]$ that $[xz,y,xz]\equiv 1$ and $[xz,y,xz,w]\equiv 1$. Therefore
$[x,y,z,w]\equiv [z,y,x,w]^{-1}$. 
Similarly, $y$ and $w$ can also be interchanged at the expense of inverting the term. 
This implies $P\equiv Q$, so $P^4=[x,y,z,w]^4\equiv 1$.

Therefore, $[x,y,z,w]$ is both of order $3$ and $4$, so is trivial. It follows that all $4$-fold commutators in $F_{x,y,z,w}$ are trivial mod the $2$-Engel relation.
This concludes the proof of lemma \ref{Engel nilpotent}. \qed 

The following corollary of the proof of lemma \ref{Engel nilpotent} will be used in later sections. 

\begin{corollary} \label{Engel corollary} \sl
Suppose $G$ is a group normally generated by $g_1,\ldots, g_n$. Let $g\in G^k$ be an element of the $k$-th term of the lower central series, $4\leq k\leq n$. Then $g$ may be represented in the Milnor group $MG$ as a product of (conjugates of) $k$-fold commutators of the form $[h_1,\ldots, h_k]$ where two of the elements $h_i$ are equal to each other and to a product of two generators, $h_j = h_m=g_{i_1} g_{i_2}$ for some $j\neq m$, and each other element $h_i$ is one of the generators $g_1, \ldots, g_n$.  
\end{corollary}

{\it Proof of corollary \ref{Engel corollary}.}
First consider the case $k=4$. Then $g$, considered as an element of $(MG)^4$, equals a product of conjugates of $4$-fold basic commutators $[g_{i_1},\ldots,g_{i_4}]$ with distinct indices. Therefore it suffices to consider $g=[x,y,z,w]$ where 
$x,y,z,w$ are distinct normal generators.
The triviality modulo the $2$-Engel relation, $g\equiv 1$,  was established in the proof of lemma \ref{Engel nilpotent} as a consequence of $g=g^4\cdot g^{-3}$, where 
$g^4\equiv 1$ and $g^{3}\equiv 1$.

The proof that $g$ may be represented as stated in the corollary follows from a direct inspection of the instances where the $2$-Engel relation is used: 
equations (\ref{first eq}), (\ref{similarly}) and two paragraphs following  (\ref{commutators}). 
A priori there are two types of usage of the $2$-Engel relation. One type is immediately in the form of a $4$-fold commutator, as claimed in the statement of the corollary, cf. $[x,y,zw,zw]\equiv 1$ in the paragraph following (\ref{commutators}). The other type is a $3$-fold commutator, for example (\ref{first eq}), considered within a $4$-fold commutator in (\ref{another equation2}). More precisely, the equations (\ref{first eq})-(\ref{similarly}) express $[x,y,z]^3$ in the Milnor group as a product of $3$-fold commutators of the form $[z,xy,xy]$ and $[x,yz,yz]$. Using (\ref{product id}) as in (\ref{another equation2}), with the conjugation disregarded in the Milnor group, 
$g^3 =  [[x,y,z]^3, w]$ is then expressed as a product of $[z,xy,xy, w]$ and $[x,yz,yz,w]$.

For any $k\geq 4$, $g$ (considered as an element of $(MG)^k$) equals a product of conjugates of $k$-fold basic commutators $[g_{i_1},\ldots,g_{i_k}]$ with distinct indices.
The argument for $k=4$ shows that ``the initial segment'' $[g_{i_1},\ldots, g_{i_4}]$ of each factor is a product $\prod [h_{j_1},\ldots, h_{j_4}]$.
Then as in the previous paragraph, using the identities (\ref{product id}), in the Milnor group one has 
$[g_{i_1},\ldots,g_{i_k}]=[\prod [h_{j_1},\ldots, h_{j_4}], g_{i_5},\ldots, g_{i_k}]=\prod [h_{j_1},\ldots, h_{j_4}, g_{i_5},\ldots, g_{i_k}].$
 \qed

The following observations are useful for estimating the number of commutators
$[h_1,\ldots, h_k]$ needed for a given element $g\in G^k$.
Since the statement takes places in the Milnor group, it may be assumed that all generators $g_i$ that appear in each commutator $h_1,\ldots, h_k$ are distinct. 
Also, it suffices to consider commutators where (in the notation of the corollary) $1\leq j, m\leq 4$.

\begin{figure}[t]
\includegraphics[width=6.5cm]{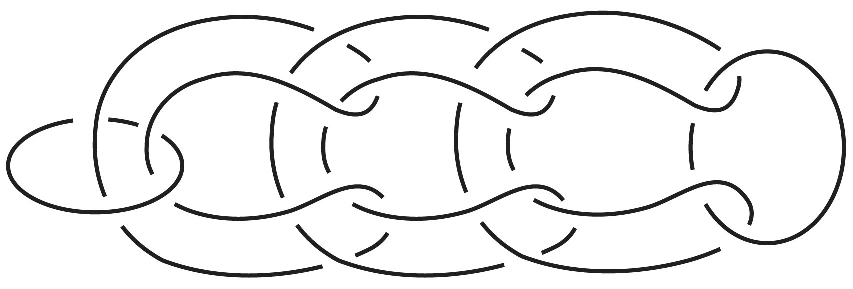}
{\small
\put(-194,15){$\gamma$}
\put(-140,-6){$x$}
\put(-95,-6){$y$}
\put(-55,-6){$z$}
\put(-5,8){$w$}
}
\caption{${\gamma}=[x,y,z,w]$.}
\label{Milnor link}
\end{figure}

It is useful to illustrate the links which are a geometric analogue of the commutators appearing in the statement of Corollary \ref{Engel corollary}.
Figure \ref{Milnor link} shows a $5$-component link (obtained by iterated Bing doubling the Hopf link) where the left-most component is denoted $\gamma$ and the meridians to the other $4$ components are labeled $x, y, z, w$. Then $\gamma$ reads off the commutator $[x,y,z,w]$ in the complement of the other $4$ components. 

\begin{figure}[ht]
\includegraphics[width=14.3cm]{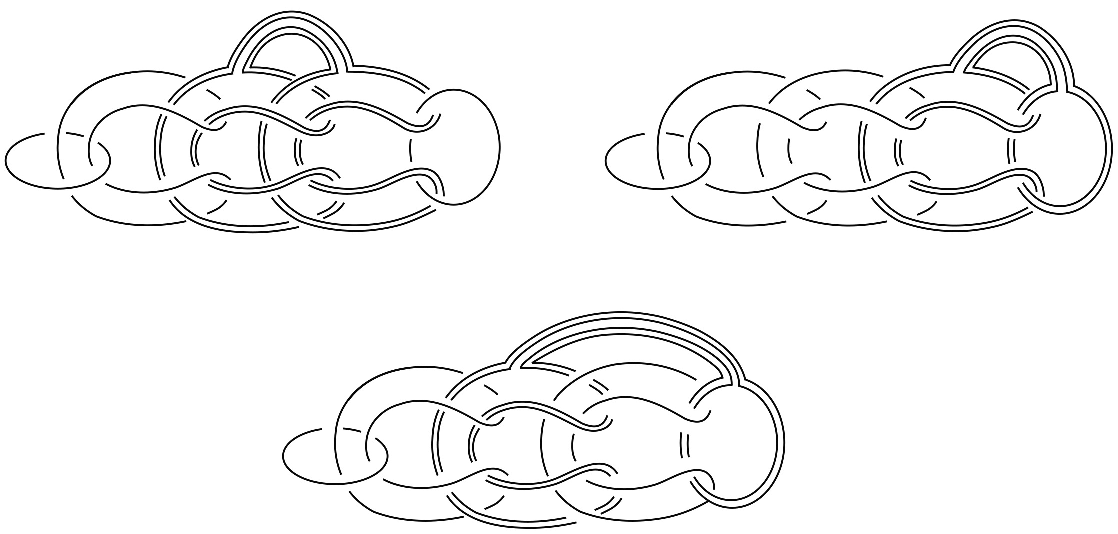}
{\small
\put(-416,125){${\gamma}_1$}
\put(-364,103){$x$}
\put(-318,101){$y$}
\put(-318,115){$z$}
\put(-229,118){$w$}
\put(-197,125){${\gamma}_2$}
\put(-143,104){$x$}
\put(-57,101){$z$}
\put(-100,104){$y$}
\put(-57,116){$w$}
\put(-414,175){(a)}
\put(-194,175){(b)}
\put(-312,15){${\gamma}_3$}
\put(-260,-7){$x$}
\put(-177,-6.5){$z$}
\put(-215,-6){$y$}
\put(-217,7){$w$}
\put(-308,60){(c)}
}
\caption{(a): $\!{\gamma}_1=[x,yz,yz,w]$, (b): $\!{\gamma}_2=[x,y,zw,zw]$,
(c): $\!{\gamma}_3=[x,yw,z,yw]$.}
\label{Engel links1}
\end{figure}

The ``elementary'' $2$-Engel links, geometric analogues of the commutators $[h_1,\ldots, h_4]$ that come up in the proof of corollary \ref{Engel corollary} for $k=4$, are obtained by band-summing two components at a time and then taking a parallel copy, as shown in figure \ref{Engel links1}.
All the links in figures \ref{Milnor link}, \ref{Engel links1} have non-zero non-repeating $\bar\mu$ invariants. It is a very salient feature, as we shall see, that for the links in figure \ref{Engel links1} the $3$-manifold obtained by $0$-framed surgery is also obtained by $0$-framed surgery on a link with vanishing non-repeating $\bar\mu$ invariants: merely slide one of the parallel copies over the other to obtain a homotopically trivial link (see figure \ref{Engel links}).

\section{The A-B slice problem} \label{ABslice section}
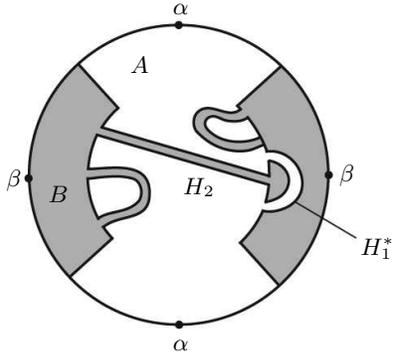
\begin{figure}[ht]
\includegraphics[width=4.5cm]{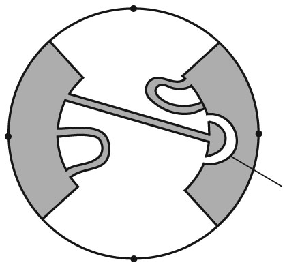}
{\scriptsize
    \put(-71,121){${\alpha}$}
    \put(-71,-6){${\alpha}$}
    \put(-87,99){$A$}
    \put(-134,56){${\beta}$}
    \put(-8,58){${\beta}$}
    \put(-118,50){$B$}
    \put(0,30){${ H}^*_1$}
    \put(-67,53){${ H}_2$}}
 \caption{A decomposition in two dimensions:
$D^2=A\cup B$, $B$ is shaded. $\alpha,\beta$ are linked $0$-spheres in $\partial D^2$. The notation for handles is discussed in section \ref{relative slice subsection}.}
\label{2D figure}
\end{figure}

We start by recalling the definition of an AB slice link from \cite{F3}. Section \ref{relative slice subsection}  summarizes the relative-slice formulation of the AB slice problem, and section \ref{definitions section} defines the notion of a homotopy AB slice link, used in theorem \ref{ABslice theorem}.
\begin{definition} \label{decomposition definition}
A  {\it decomposition} of $D^4$ is a pair of compact codimension zero smooth submanifolds with boundary $A,B\subset D^4$, satisfying conditions $(1)$-$(3)$ below. Denote $$\partial^{+} A=\partial A\cap \partial D^4, \; \; \partial^{+} B=\partial B\cap \partial D^4,\; \; \partial A=\partial^{+} A\cup {\partial}^{-}A, \; \; \partial B=\partial^{+} B\cup {\partial}^{-}B.$$ (1) $A\cup B=D^4$,\\ (2) $A\cap B=\partial^{-}A=\partial^{-}B,$ \\
(3) $S^3=\partial^{+}A\cup \partial^{+}B$ is the standard genus $1$
Heegaard decomposition of $S^3$.
\end{definition}
Each side $A, B$ of a decomposition has an attaching
circle (a distinguished curve in the boundary:
${\alpha}\subset\partial A, {\beta}\subset\partial B$) which is the
core of the solid torus $\partial^{+}A$, respectively  $\partial^{+}B$. The two curves ${\alpha}, {\beta}$ form the Hopf
link in $S^3=\partial D^4$. Figure \ref{2D figure} illustrates the notion of a decomposition in $2$ dimensions. 
The ``trivial'' decomposition of $D^4$ is given by $(A,{\alpha})=$unknotted $2$-handle: $(D^2\times D^2, \partial D^2\times 0)$ and $(B, {\beta})$ = collar: $(S^1\times D^2\times I, S^1\times 0\times 0)$. 
See \cite{FL}, \cite{K}, and section \ref{homotopy solution} below for interesting examples of decompositions.

Given an $n$-component link $L=(l_1,\ldots,l_n)\subset S^3$, consider
its untwisted parallel copy $L'=(l'_1, \ldots, l'_n)$. 

\begin{definition}  \label{AB slice definition}
The link
$L$ is {\it $A$-$B$ slice} if there exist decompositions $(A_i, B_i)$,
of $D^4$ and self-homeomorphisms ${\phi}_i, {\psi}_i$
of $D^4$, $i=1,\ldots,n$ such that all sets in the collection
${\phi}_1 A_1, \ldots, {\phi}_n A_n, {\psi}_1 B_1,\ldots, {\psi}_n
B_n$ are disjoint and satisfy the boundary data:
${\phi}_i({\partial}^{+}A_i)$  is a tubular neighborhood of $l_i$
and ${\psi}_i({\partial}^{+}B_i)$ is a tubular neighborhood of
$l'_i$, for each $i$.
\end{definition}

The restrictions ${\phi}_i|_{A_i}$, ${\psi}_i|_{B_i}$ give disjoint embeddings  of the
entire collection of $2n$ manifolds $\{A_i, B_i\}$ into $D^4$. Moreover, these
re-embeddings are {\it standard}: they are restrictions of
self-homeomorphisms of $D^4$, so the complement
$D^4\smallsetminus {\phi}_i(A_i)$ is homeomorphic to $B_i$, and
$D^4\smallsetminus {\psi}_i(B_i)\cong A_i$. Analogues of this condition in the homotopy context are introduced in definitions \ref{link homotopy ABslice definition}, \ref{homotopy AB slice}.

References \cite{F2, F3} reformulated $4$-dimensional topological surgery conjecture for free groups in terms of the AB slice problem for GBRs (the Generalized Borromean rings). Figure \ref{5 components} shows a representative link from this family.  
The proof of theorem \ref{ABslice theorem} will use the relative-slice formulation of the AB slice problem, discussed next.

\subsection{The relative slice problem.} \label{relative slice subsection}
Our summary of this approach to the  AB slice problem follows \cite{FL}; the reader is referred to this reference for further details.

Given a decomposition $D^4=A\cup B$, without loss of generality it may be assumed \cite{FL} that each side $A, B$  has a handle decomposition 
(rel. collar $S^1\times D^2\times I$) with only $1$- and $2$-handles. Fix the notation: $A=({\partial}^+ A)\times I\,\cup\, {\mathbf H}_1\,\cup\, {\mathbf H}_2$.
As usual in Kirby calculus \cite{GS}, the $1$-handles will be considered as standard $2$-handles ${\mathbf H}_1^*$ removed from the collar, $A=({\partial}^+ A\times I\setminus  {\mathbf H}^*_1)\cup {\mathbf H}_2$. In the illustration in figure \ref{2D figure} the side $A$ has three $2$-handles and a single $1$-handle.

\begin{definition} \label{throgh definition}
If $h_2$ is an embedded $2$-handle in $D^4$ (within any handle decomposition) and $h_1^*$ is a $2$-handle deleted from a collar as above, we say $h_2$ {\it does not go through} $h_1^*$ if $h_2$ is disjoint from co-core$(h_1^*)$.
\end{definition}

Consider a slightly smaller $4$-ball $D'$, equal to the original $D^4$ minus the collar $({\partial}^+ A)\times I$. The removed handles ${\mathbf H}^*_1$ may be considered as $2$-handles attached (with zero framing) to $D'$.
Note that in figure \ref{2D figure} none of the $2$-handles ${\mathbf H_2}$ go through ${\mathbf H}_1^*$. 
This condition does not have to be satisfied for an arbitrary handle decomposition of a given submanifold of $D^4$, but it will hold for all decompositions constructed in this paper, as stated in Condition \ref{handle condition}.

\begin{remark}
The details of the embeddings of $\{A_i , B_i \}$ into $D^4$ are important in the A-B slice problem. It was shown in \cite{K} that any link $L=(l_1,\ldots,l_n)$ with trivial linking numbers is {\it weakly} A-B slice: there exist $n$ decompositions $D^4=A_i\cup B_i$ and disjoint embeddings of the entire collection of $\{ A_i, B_i\}$ into $D^4$ with the boundary data given by $L$ and its parallel copy. To be relevant to the surgery conjecture, these disjoint embeddings have to be {\it standard}, as discussed in the paragraph following definition \ref{AB slice definition}.
We record the relevant information  about embeddings in condition \ref{handle condition}; analogous statements in the relative-slice setting and in the homotopy context are given in conditions 
\ref{standard reembedding condition} and \ref{standard embedding definition} respectively.
\end{remark}

\begin{condition} \label{handle condition}
For each side $C=A$, $B$ of a decomposition $D^4=A\cup B$,
the $2$-handles ${\mathbf H}_2$ of $C$ do not go over the handles 
${\mathbf H}^*_1$ corresponding to the $1$-handles of $C$. (N.B. $2$-handles of $A_i$ may certainly pass through both $1$-handles $H_1^*$ of $B_j$ and the $1$-handles of $A_j$, for all $j\neq i$.)
\end{condition}

Beyond this requirement, in our decompositions each $2$-handle of $C$ is embedded in a standard way (i.e. is unknotted) in $D^4\, \setminus$collar on $\partial^+ C$.
It follows that (except for a single $2$-handle) the $1$-handles of each side are in
one-to-one correspondence with the $2$-handles of the complement. 
We treat this as an additional convenient property, but not part of the definition. (But note that an easy subdivision argument allows any decomposition to be reduced to one of this type at the expense of increasing the number of handles.)
In the decompositions introduced in section \ref{homotopy solution}, the $A$-side has
a zero framed $2$-handle attached to the core of the solid torus ${\partial}^+ A$; this is the ``distinguished'' handle of $A$ which does not have a counterpart on the $B$-side. 

Usually a Kirby diagram is a labeled link in $S^3$, thought of as $S^3=\partial D^4$. To describe $A$ and $B$ in Kirby notation we consider labeled links in the solid torus $\cong S^1\times D^2$, now thought of as the inner boundary, ${\partial}^-$,  of a collar $(S^1\times D^2)\times I$. In our applications  the (isotopy class of the) product structure $S^1\times D^2$ is fixed by the embedding in $S^3$, so integer framing coefficients are interpreted  as usual (self-linking numbers). Similarly, an unlink consisting of dot-bearing components means: delete standard $2$-handles from the collar.
The outer boundaries of these collars are $\partial^+ A$ and $\partial^+ B$, together $\partial^+ A\cup \partial^+ B=S^3$. One may think of surgery on the core circle of $\partial^+ A$ (i.e. the effect of $2$-handle attachment) as transforming $\partial^+ A$ into (something isotopic rel boundary to) $\partial^+ B$.
 
 A Kirby diagram for $B$ may be obtained by taking a Kirby diagram in the solid torus for $A$, performing the surgery as above, and replacing all zeros with dots, and conversely all dots with zeros. (Note that the
$2$-handles in all our decompositions are zero-framed.) To fix the notation, denote the distinguished $2$-handle of $A$ by $H_2$  (as in figure \ref{2D figure}), and the rest of the $2$-handles of $A$ by $\overline{{\mathbf H}}_2$.

Suppose an $n$-component link $L$ is AB slice, with decompositions $D^4=A_i\cup B_i$ and homeomorphisms ${\phi}_i, {\psi}_i$, for $i=1,\ldots, n$. Denote by $D^4_0$ the smaller $4$-ball obtained by removing from $D^4$ the collars on the attaching regions ${\phi}_i(\partial^+ A_i), {\psi}_i(\partial^+ B_i)$ of all $2n$ submanifolds $\{ {\phi}_i(A_i)$, ${\psi}_i(B_i)\}$. Let ${\mathcal H}_2$ denote the $2$-handles of all these submanifolds, and ${\mathcal H}_1^*$ the $2$-handles removed from the collars, corresponding to the $1$-handles. As above, consider ${\mathcal H}_1^*$ as zero-framed $2$-handles attached to $D^4_0$.

A more precise description of ${\mathcal H}_2$, ${\mathcal H}_1^*$ may be given as follows. Add  the superscript $i$ to the handle notation: the $2$-handles of $A_i$ are ${\mathbf H}_2^i= H_2^i\cup \overline{{\mathbf H}_2^i}$, where $H_2^i$ is a ``distinguished'' $2$-handle of $A_i$, and $1$-handles of $A_i$ correspond to ${\mathbf H}^{i*}_1$. Then 

\begin{equation} \label{rel slice}
{\mathcal H}_2\; =\; \underset{i=1,\ldots, n}{\mathsmaller{\bigcup}} {\phi}_i({\mathbf H}_2^i)\, \underset{i=1,\ldots, n}{\mathsmaller{\bigcup}} {\psi}_i({\mathbf H}_1^{i*}), \; \; \;
{\mathcal H}_1^*\; =\; \underset{i=1,\ldots, n}{\mathsmaller{\bigcup}} {\phi}_i({\mathbf H}_1^{i*})\, \underset{i=1,\ldots, n}{\mathsmaller{\bigcup}} {\psi}_i(\overline{{\mathbf H}_1^i}).
\end{equation}

\begin{figure}[ht]
\includegraphics[height=3.6cm]{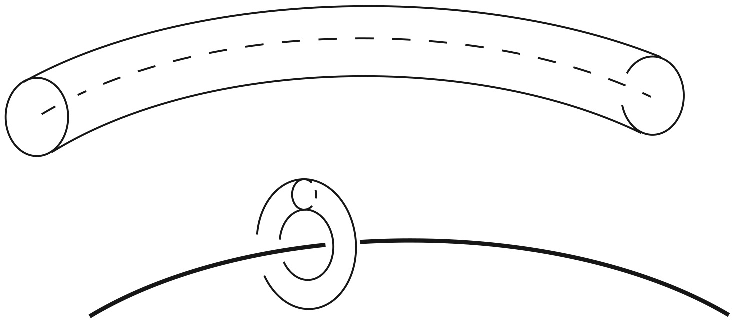}
{\small
    \put(-227,62){$l'_i$}
    \put(-214,-2){$l_i$}
\put(-113,-3){$(J_i, K_i)$}
\put(-110,51){$(\widehat K_i, \widehat J_i)$}}
\put(-96,64){$\cup$}    
\put(-123,5){$\slantedcup$}
 \caption{Stabilizing components added to an A-B slice link $L=\{l_i\}$: link pairs 
$(J_i, K_i)\subset\, $neighborhood of a meridian to $l_i$,
$(\widehat K_i, \widehat J_i)\subset \, $solid torus neighborhood of a parallel copy $l'_i$. (The parallel copy $l'_i$ is not part of the link.)
A 
diffeomorphism between the solid tori exchanging their meridian and longitude takes $K_i$ to $\widehat K_i$ and $J_i$ to $\widehat J_i$. (To relate this to the notation in definition \ref{rel slice definition}: $J=L\cup_i  J_i\cup_i \widehat {K}_i$, $K=\cup_i K_i\cup_i \widehat{J}_i$.)}
\label{RelSlice figure}
\end{figure}
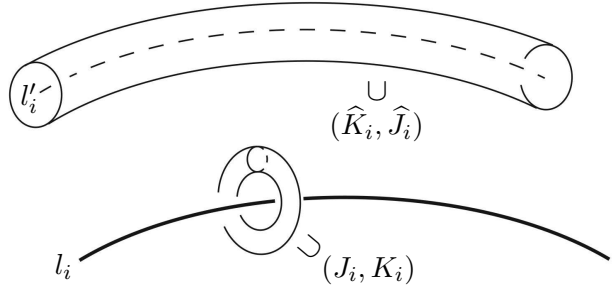

Consider the following two links $J, K$ in $S^3=\partial D^4_0$, which may be read off from the Kirby diagrams of the $\{ A_i, B_i\}$.
Let $J$ denote the attaching curves of the $2$-handles ${\mathcal H}_2$, 
and $K$ the attaching curves of the $2$-handles ${\mathcal H}_1^*$. (Note that ${\mathcal H}_2\subset D^4_0$, and ${\mathcal H}_1^*$ are attached with zero framings to $D^4_0$ along $K$.)  The components belonging to $J$ and $K$ will be labeled accordingly in the figures in this paper; the link $K$ is traditionally drawn in red on blackboards.
We refer to the pair $(J,K)$ as the ``stabilization'' the original link $L$. (Note that $L$ is included in $J$ as the attaching curves of the distinguished $2$-handles $\{ H_2^i\}$.)  The structure of the stabilization links, which is a consequence of the duality between the $1$- and $2$-handles of the two sides of each decomposition, is shown in figure \ref{RelSlice figure}.

\begin{definition} \label{rel slice definition}
A link pair $(J, K)$ in $S^3=\partial D^4_0$ is called {\it relatively slice} if the components of $J$ bound disjoint, smoothly  embedded disks in the handlebody
$$ H_K:=D^4_0\, \cup\,\text {zero-framed 2-handles attached along}\; K.$$
\end{definition}

If a link $L$ is AB slice, by construction the associated link pair $(J, K)$ is then relatively slice. Moreover, since the embeddings 
${\phi}_i(A_i)$, ${\psi}_i(B_i)$ are restrictions of self-homeomorphisms ${\phi}_i(A_i)$, ${\psi}_i(B_i)$ of the $4$-ball, the following analogue of condition \ref{handle condition} holds 
for the relative slicing $(J,K)$.

\begin{condition} \label{standard reembedding condition}
Let $S$ be any submanifold in the collection $\{ {\phi}_i(A_i)$, ${\psi}_i(B_i)\}$. Then 
after an isotopy  (depending on $S$) of the embedded handlebody $H_K=D^4_0\, \cup_K\, $($2$-handles) the slices for the components of $J$ corresponding to $S$ do not go through the $2$-handles $\{ {\mathcal H }\}$
attached to $D^4_0$ along the components of $K$ corresponding to the same submanifold $S$.
\end{condition}
\begin{figure}[ht]
\includegraphics[height=4.2cm]{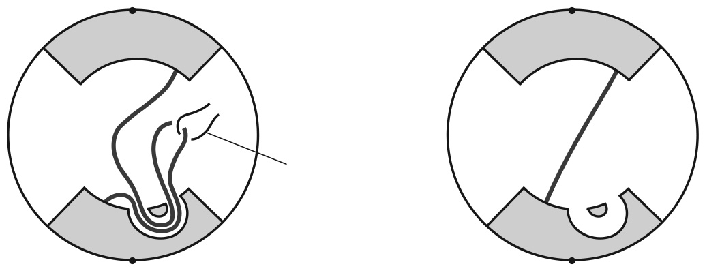}
{\small
    \put(-280,14){$S$}
    \put(-183,39){$S'$}
}
 \caption{An illustration of condition \ref{standard reembedding condition}: the $2$-handles of a submanifold $S$ may ``go over
its $1$-handles'' and link another submanifold $S'$ in $D^4$ (left). However they do not go over its $1$-handles after an isotopy (right), where other submanifolds are disregarded.}
\label{standard embedding figure}
\end{figure}

Note that the statement of condition \ref{standard reembedding condition} in general indeed requires an isotopy: as illustrated in figure \ref{standard embedding figure},
$2$-handles of $S$ may link other submanifolds $S'$. There is an isotopy ``straightening out'' the $2$-handles of $S$ as shown on the right in the figure, but the condition may not be achieved simultaneously for all submanifolds $\{ {\phi}_i(A_i)$, ${\psi}_i(B_i)\}$.

\subsection{Homotopy A-B slice problem} \label{definitions section} We now turn to the definition of  a homotopy  AB slice link, referred to in the statement of theorem \ref{ABslice theorem}. In fact, we state two natural versions of the definition. It will be shown in section \ref{homotopy solution} that theorem \ref{ABslice theorem} holds in both contexts.  The first notion is motivated by link-homotopy theory (section \ref{Milnor group section}):

\begin{definition} \label{link homotopy ABslice definition} ({\bf Link-homotopy A-B slice})
An $n$-component link $L$ is {\it link-homotopy A-B slice} if there exist decompositions $D^4=A_i\cup B_i$, $i=1,\ldots, n$ and handle decompositions of
the submanifolds $A_i, B_i$ so that the corresponding relative-slice problem $(J, K)$ has a link-homotopy solution. That is, in the notation of definition \ref{rel slice definition} the components of $J$ bound disjoint maps of disks $\Delta$ in the handlebody $H_K$. Moreover, the disks $\Delta$ are subject to condition \ref{standard embedding definition} below.
\end{definition}

Recall that the free group action in the context of the AB slice problem is encoded in condition \ref{standard reembedding condition}. A  stronger version of that condition is to omit a reference to an isotopy and require that for no $S$ do its $2$-handles pass over the dual representation of its $1$-handles. We use this stronger version to define an analogue for a link-homotopy A-B slice link:

\begin{condition} \label{standard embedding definition} 
Let $S$ be any submanifold in the collection $\{ {\phi}_i(A_i)$, ${\psi}_i(B_i)\}$. 
Then the maps of disks $\Delta$ for the components of $J$ corresponding to $S$ do not go through the $2$-handles 
attached to $D^4_0$ along the components of $K$ corresponding to the same submanifold $S$. 
\end{condition}

Theorem \ref{ABslice theorem} will be established for a stronger version of a homotopy solution which is defined next. Suppose a possibly disconnected codimension zero submanifold
$(C,\partial^+ C)\subset (D^4, S^3)$ is given, together with a  $1, 2$-handle decomposition of $(C,\partial^+ C)$. 
As in the beginning of section \ref{relative slice subsection}, consider the $4$-ball $D'=D^4\smallsetminus (\partial^+ C)\times I$, and consider the $1$-handles of $C$ as standard slices ${\mathbf H}_1^*$ removed from the collar.
We say that a map
$f\co (C, \partial^+ C)\longrightarrow (D^4, S^3)$  is {\it homotopy standard} if there exists a $1$-parameter family of maps $f_t$ connecting $f$ and ${\rm id}\co C\subset D^4$ such that 
\begin{enumerate}
\item  
the restriction of $f_t$ to $(\partial^+ C)\times I \smallsetminus {\mathbf H}_1^*$ is the identity map for all $t$, and 
\smallskip
\item 
the images of the $2$-handles are contained in $D'$ and are disjoint for all $t$.
\end{enumerate} 

\begin{definition} \label{homotopy AB slice} ({\bf Homotopy A-B slice})
An $n$-component link $L$ is {\it homotopy A-B slice} if there exist  decompositions $D^4=A_i\cup B_i$, a $1,2$-handle structure for each submanifold $A_i, B_i$,  and disjoint
maps of all $2n$ submanifolds $\{ A_i, B_i\}$ into $D^4$ with the boundary data corresponding to $L$ and its parallel copy (as in definition \ref{AB slice definition}), 
such that the restriction to each $A_i, B_i$, $i=1,\ldots, n$ is homotopy standard.
\end{definition}

Here is a brief outline of the way (singular) slices will be found for a homotopy solution to the relative-slice problem in theorem \ref{ABslice theorem}.
One may band sum the
components of $J$ with (an arbitrary number of) parallel copies of the components of $K$. These bands correspond to index $1$ critical points 
of the slices with respect to the radial Morse function on $D^4_0$. Parallel copies of each component $K_i$ bound disjoint copies of the core of the $2$-handle attached to $K_i$. For a suitable choice of band sums, the resulting band-summed link $J'$ will be null-homotopic.
Then the construction of singular slices is completed by capping off the components of $J'$ with disjoint maps of disks in $D_0^4$. Of course the crucial part of the proof is the construction of the decompositions $D^4=A_i\cup B_i$ enabling this strategy to succeed.

\section{Proof of theorem \ref{ABslice theorem}: A homotopy solution to the AB slice problem.}\label{homotopy solution}

As discussed in the introduction, the generalized Borromean rings 
$\{\!${\sl Bing(Hopf)}$\!\}$
form a collection of model surgery problems. We start by noting that highly Bing doubled links in $\{\!${\sl Bing(Hopf)}$\!\}$ are still universal for surgery. In the setting of capped gropes this follows from grope height raising: for any $n\geq 3$ the attaching curve of a capped grope $g^c$ of height $3$ bounds a capped grope of height $n$ in the (untwisted) thickening of $g^c$ \cite[Proposition 2.7]{FQ}. (The proof of \cite[Theorem 5.1A]{FQ} explains how to get a capped grope of height $3$, starting from  a surgery kernel.)

\begin{prop} \label{higher prop} \sl 
Let $n\geq 5$ be fixed. 
Denote by $\{\!$Bing(Hopf)$\!\}_n$ the links $L$ in
$\{\!$Bing(Hopf)$\!\}$ satisfying 
$ML/(ML)^n\cong MF/(MF)^n$.
 Then $\{\!$Bing(Hopf)$\!\}_n$
forms a collection of model surgery problems.  
\end{prop}
\begin{figure}[ht]
\includegraphics[height=3.2cm]{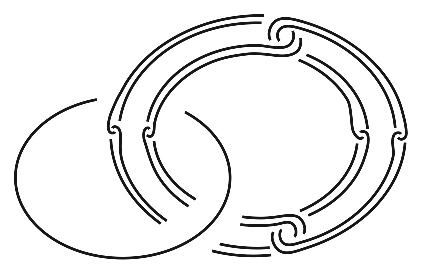}
{\small
\put(-146,27){$l_1$}
\put(-102,73){$l_2$}
\put(-107,20){$l_3$}
\put(-15,72){$l_4$}
\put(-14,12){$l_5$}
}
 \caption{A link $L\in\{\!${\sl Bing(Hopf)}$\!\}_5$}
\label{5 components}
\end{figure}
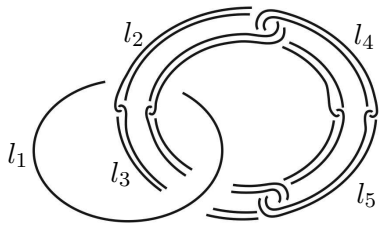
The Milnor group condition picks out the class $\{\!${\sl Bing(Hopf)}$\!\}_n$,
obtained from the Hopf link by (ramified) Bing doubling performed at least $n-2$ times.
The usually referenced \cite{F1} class of universal problems is, in this notation,
$\{\!${\sl Bing(Hopf)}$\!\}_3$, but grope height raising \cite{FQ} allows one to restrict to any coinitial segment, such as $\{\!${\sl Bing(Hopf)}$\!\}_n$, $n\geq 5$.
For our purposes, in light of proposition \ref{weak length 4}, it suffices to consider $n=5$.
In the proof of theorem \ref{ABslice theorem} first consider the case where $L\in  \{\!${\sl Bing(Hopf)}$\!\}_5$ is almost homotopically trivial, meaning absent any component the link becomes homotopically trivial.  (In this collection of links this is equivalent to $L$ being Brunnian). This means that $L$ is obtained from the Hopf link by iterated Bing doubling without ramification. This case captures the idea of the proof; at the end of this section we show what adjustments need to be made in the general case. To be very specific, consider one of the smallest representatives of $\{\!${\sl Bing(Hopf)}$\!\}_5$, the $5$-component link in figure \ref{5 components}. The reason $5$ suffices is that $l_1$ lies in the $4$th term of the lower central series of ${\pi}_1(S^3\smallsetminus (l_2\cup\ldots\cup l_5))$, and Corollary \ref{Engel corollary} gives a useful way of writing such elements in terms of $2$-Engel relations.

\begin{figure}[ht]
\includegraphics[width=12.5cm]{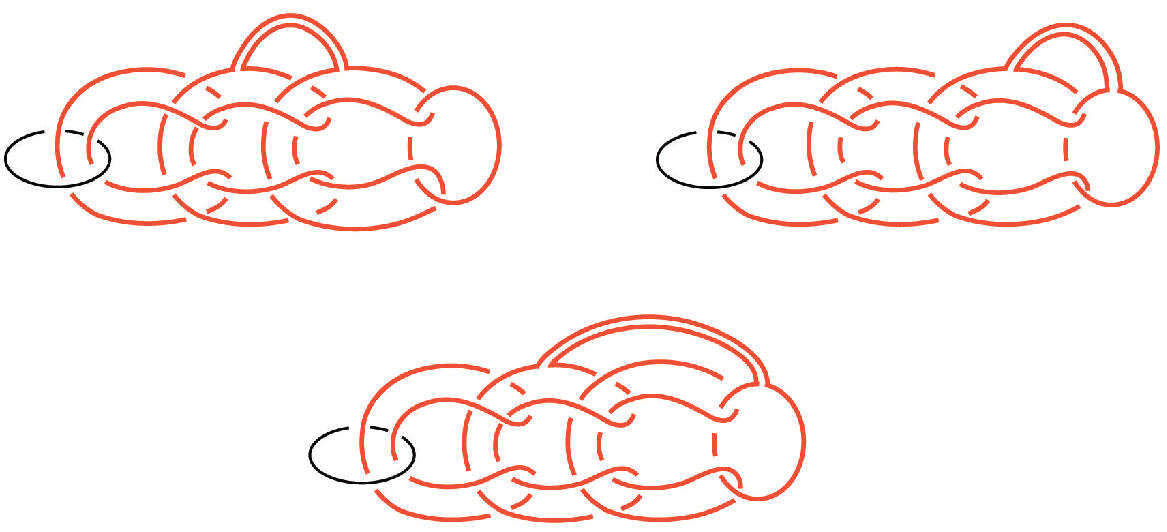}
{\small
\put(-365,138){(a$'$)}
\put(-355,95){$J$}
\put(-280,80){$K$}
\put(-167,138){(b$'$)}
\put(-155,94){$J$}
\put(-90,80){$K$}
\put(-275,48){(c$'$)}
\put(-265,5){$J$}
\put(-190,-10){$K$}
}
\caption{}
\label{Engel links}
\end{figure}
To define the decompositions giving a homotopy A-B slice solution for $L$, consider the links in figure \ref{Engel links}.
These links should be compared with the links in figure \ref{Engel links1}. The crucial difference is that now the band-summed curves are taken {\it without} a parallel copy, and unlike the links in figures \ref{Engel links1} each of the three links in figure \ref{Engel links} is homotopically {\it trivial}. 

\begin{figure}[ht]
\includegraphics[width=13cm]{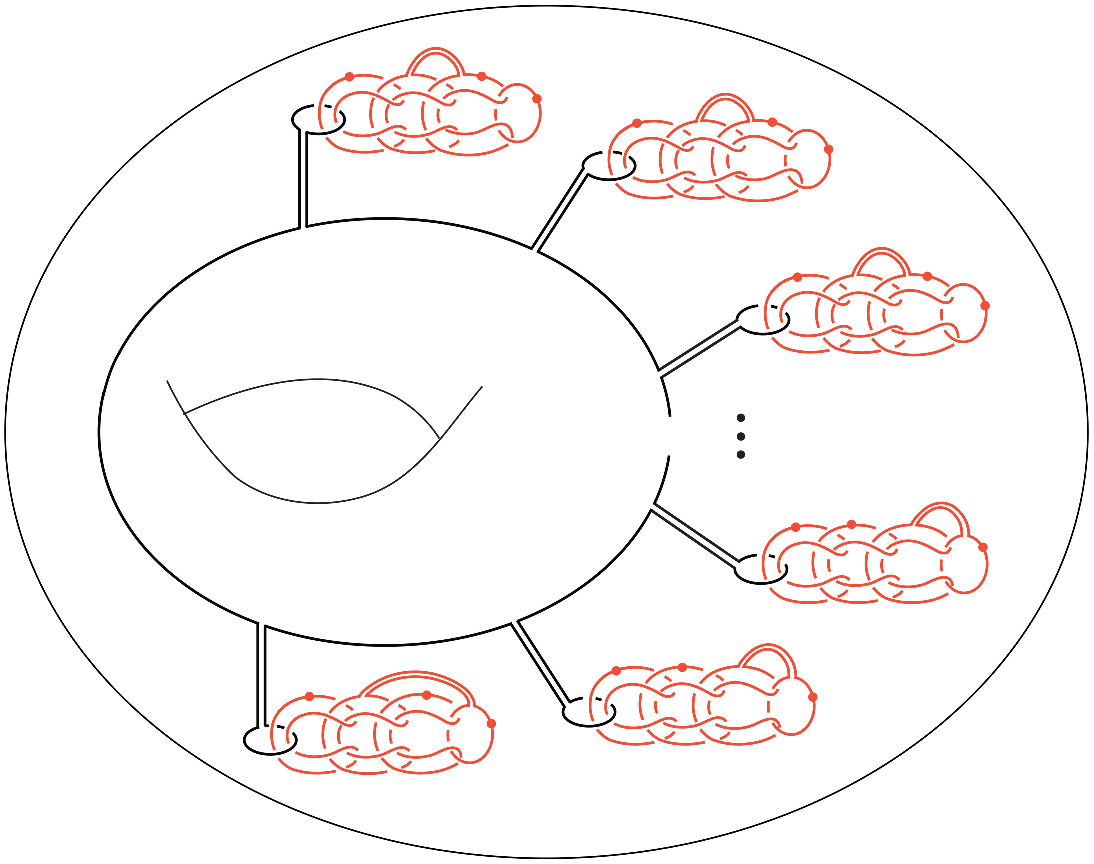}
\put(-333,183){$0$}
{\scriptsize \put(-105,141){$K_4,\ldots, K_9$}}
\put(-325,93){$c$}
\scriptsize{
\put(-252,233.5){$z$}
\put(-230,233.5){$x,y$}
\put(-190,240){$w$}
\put(-152,218){$x$}
\put(-130,218){$y,z$}
\put(-90,226){$w$}
\put(-201,271){$K_1$}
\put(-158,259){$K_2$}
\put(-54,203){$K_3$}
\put(-70,76){$K_{10}$}
\put(-125,26){$K_{11}$}
\put(-220,17){$K_{12}$}
\put(-173,103){$J$}
}
 \caption{Part $A_1$ of the decomposition $D^4=A_1\cup B_1$.}
\label{Torus figure}
\end{figure}
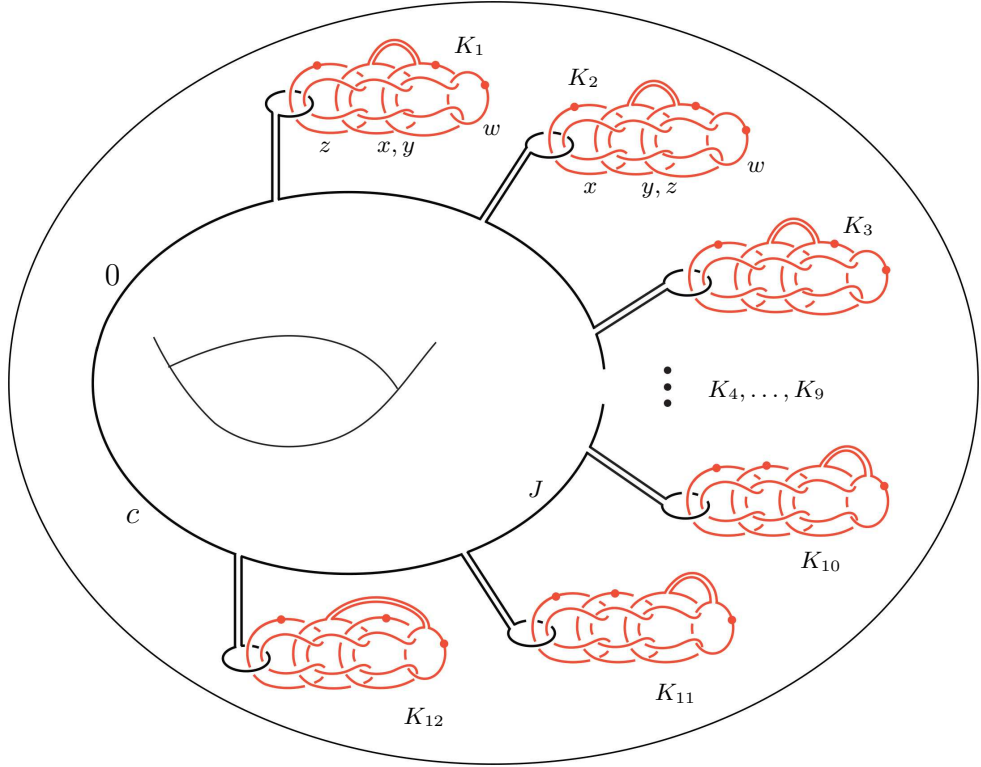

We are now in a position to define the relevant decompositions of the $4$-ball.
For the components $l_i,\, 2\leq i\leq 5$ consider the trivial decomposition $D^4=A_i\cup B_i$ where $A_i\, =\,$unknotted $2$-handle, $B_i\, =\, $collar on the attaching curve.
For the first component consider the decomposition determined by the side $A_1$ shown in figure \ref{Torus figure}. There is one zero-framed $2$-handle and $12\times 3=36$ $1$-handles. (This very specific handle description is given for the $5$-component link in figure \ref{5 components} in part to take a note of the complexity of our homotopy A-B slice solution. The definition for a general $L\in \{\!$Bing(Hopf)$\!\}_5$ is given at the end of the proof.)
The curves representing $1$-handles are dotted (standard notation in Kirby calculus \cite{GS}).  They are traditionally drawn red on the blackboard due to the role they play in the relative-slice setting. 

A precise definition of the construction  in figure \ref{Torus figure} is as follows.
The Kirby diagram is drawn in the solid torus neighborhood of the attaching curve (${\partial}^+A$ in the notation of definition \ref{decomposition definition}).
Start with a single curve, the zero-framed attaching curve for the $2$-handle of $A_1$, given by the core of the solid torus. Consider a total of 12 links embedded in disjoint $3$-balls in (solid torus$\;\smallsetminus\;$its core):
six copies of the link (a$'$) in figure \ref{Engel links}, four copies of (b$'$) and two copies of (c$'$). The proof below shows how the number and types of links are determined by the algebraic structure of the $2$-Engel relation (seen in the proof of lemma \ref{Engel nilpotent}).
In each of these 12 links the left-most component is band-summed with the  previously chosen core curve. It is convenient to use the bands shown in figure \ref{Torus figure}, so that the resulting band-summed curve $c$ is isotopic to the original core curve of the solid torus. To fix the embedding of this handlebody $A_1$ into $D^4$ we specify that the $2$-handle is embedded into the $4$-ball in a standard (unknotted) way. This determines the decomposition $D^4=A_1\cup B_1$.

\begin{prop} \label{homotopy standard prop} \sl
There exist disjoint, homotopy standard maps (in the sense of definition \ref{homotopy AB slice}) of $A_1, \ldots, A_5$ into $D^4$ with the boundary data corresponding to the link $L$ in figure \ref{5 components}.
\end{prop}

{\it Proof.}
The relative-slice problem is an embedding question for disks obtained by embedding the link in figure \ref{Torus figure} in a neighborhood of $l_1$ in figure \ref{5 components}. 
Proposition \ref{homotopy standard prop} solves this problem in the weaker context  of disjoint maps.
We identify the curve $c$ with the component $l_1$. Then the relative-slice problem concerns the pair of links $(J, K)$ where $J=L$ is the link in figure \ref{5 components}, and $K$ consists of the 36 (red) dotted curves in figure \ref{Torus figure}, embedded in a neighborhood of $l_1$. More precisely, $K=\cup _i K_i$ consists of 12 links $K_i$.

Using the commutator notation from section \ref{Engel subsection}, the component $l_1$ of the link $L$ in figure \ref{5 components} represents the element
\begin{equation} 
\label{commutator eq}
l_1\, =\, [[m_2,m_3],[m_4,m_5]]\, =\, 
[m_2,m_3,m_4,m_5]\, \cdot \, [m_3,m_2,m_4,m_5]^{-1}
\end{equation}
in the free Milnor group $M{\pi}_1(S^3\smallsetminus (l_2\cup\ldots\cup l_5))\cong MF_{m_2,\ldots, m_5}$. 
The second equality above follows from the Hall-Witt identity  (\ref{Hall-Witt eq}) where, as usual in the Milnor group setting, conjugation may be omitted. 
The choice of a basepoint for $l_1$ also does not affect the expression (\ref{commutator eq}) since these commutators are of maximal length in the Milnor group, this is discussed further in remark \ref{conjugation}. 
Detailed calculations of this type may be found, for example, in \cite{K1}, where the algebra is related to grope geometry.

Apply corollary \ref{Engel corollary} to $g=l_1\in (MF_{m_2,\ldots,m_5})^4$.
It is not difficult to see the number and the types of such commutators $[h_1,\ldots, h_4]$ that come up, going through the proof of lemma \ref{Engel nilpotent}. Specifically, consider each of the two basic commutators in (\ref{commutator eq}), for example start with $[m_2,m_3, m_4, m_5]$. To match the current notation with that of lemma \ref{Engel nilpotent}, set 
\begin{equation} \label{correspondence}
x=m_2,\; y=m_3, \; z=m_4, \; w=m_5.
\end{equation}
Lemma \ref{Engel nilpotent} establishes that $P:=[x,y,z,w]\equiv 1$ (is trivial mod the $2$-Engel relation) by showing $P^3\equiv 1$ and $P^4\equiv 1$. $P^3\equiv 1$ is proved in (\ref{another equation}) - (\ref{another equation2}) by representing $P$ in the free Milnor group as a product of two commutators $[h_1,\ldots, h_4]$ corresponding to the link (a) in figure \ref{Engel links1}. Then $P^4\equiv 1$ is proved in the two paragraphs following (\ref{commutators}) using commutators corresponding to one copy of (a), two copies of (b) and one copy of (c). Establish a 1-1 correspondence between the commutators $[h_1,\ldots, h_4]$ appearing in the proof and six of the links $K_i$ in the definition of $A_1$.

We implement the algebraic argument above geometrically as follows. In the relative-slice setting the slices for $J$ may go multiple times over the $2$-handles attached to $K$; we exploit this by band-summing $J$ to the components of $K$ and their parallel copies. Each link $K_i$ consists of three (dotted) components; denote by $K'_i$ the $4$-component link obtained by adding to it a parallel copy of the ``long'' band-summed curve in $K_i$. (Note that this ``reconstructs'' the links in figure \ref{Engel links1}.) 
Every time a commutator $[h_1,\ldots, h_4]$ is used in the proof of lemma \ref{Engel nilpotent}, perform a band-sum joining $l_2,\ldots,l_5$ with the corresponding link $K'_i$, paying a careful attention to the order of indices discussed next.

The first such commutator that comes up in the proof (line (\ref{first eq}), understood as a $4$-fold commutator as in (\ref{another equation2})) is $[z,xy,xy,w]$. Keeping in mind the notation (\ref{correspondence}), take a band sum of $l_4$  with the component of $K_1$ labeled $z$ in figure \ref{Torus figure}. Then band sum $l_2$ (resp. $l_3$) with the ``long component'' of $K_1$ labeled $x,y$ (resp. its parallel copy). Finally band sum $l_4$ with the component of $K_1$ labeled $w$. There is a $\pm$ choice for each band sum depending on orientations, this choice is discussed in remark \ref{orientation conjugation}.

The next commutator appearing in the proof is $[x,yz,yz,w]$, and there is a corresponding link $K_2$ reserved for band-summing into, as indicated in figure \ref{Torus figure}. 
Proceeding in this manner, perform band-summing into $K'_1, \ldots, K'_6$ corresponding to the proof of lemma \ref{Engel nilpotent} for $[m_2,\ldots,m_5]$. There is another elementary commutator, $[m_3,m_2,m_4,m_5]$ in the expression (\ref{commutator eq}) for $l_1$. Its triviality modulo $2$-Engel relations similarly gives rise to band sums into $K'_6,\ldots, K'_{12}$.

Denote the result of band-summing $l_2,\ldots, l_5$ with all $\{ K'_i \}$ by $l'_2, \ldots, l'_5$.  The link $(l_2,\ldots, l_5)$, as well as each $K'_i$, is an unlink. For a suitable choice of bands, $l'_2, \ldots, l'_5$ is also the unlink. (But even with an arbitrary choice of bands, its Milnor group is free: $M{\pi}_1(S^3\smallsetminus (l'_2\cup \ldots\cup l'_5))\cong MF_{m_2,\ldots,m_5}$.)

\begin{remark} \label{orientation conjugation} {\it (Orientations)}
One aspect of commutator calculus and of band-summing was implicit in the argument above. The choice of a meridian (generator of the Milnor group) to each link component depends on the orientation of the based loop representing it. 
For example, we used expressions for commutators, such as (\ref{commutator eq}) and the commutators in figure \ref{Engel links1}, which do not involve negative exponents of the meridians. These expressions assumed a particular choice of orientations.
Similarly, different choices of orientations of the link components result in different orientation-preserving band-sums. The commutator identity  
\begin{equation} \label{inverse eq} 
[x^{-1},y]=[y,x]^{x^{-1}}
\end{equation}
(cf. \cite{MKS}) is useful in this context. Conjugation may be omitted in the Milnor group (see remark \ref{conjugation}), so an iterated application of this identity implies that changing the orientation of any one meridian $m_{i_j}$ in a commutator $[m_{i_1},\ldots, m_{i_k}]$ inverts the commutator. Therefore various choices of orientations for each such commutator have two possible outcomes overall:  $[m_{i_1},\ldots, m_{i_k}]^{\pm 1}$. We make a choice of orientations so that the commutators $[h_1,\ldots,h_4]$ in the paragraph following remark \ref{conjugation} have the correct exponent to match the calculation in the proof of lemma \ref{Engel nilpotent}.
\end{remark}

\begin{remark} \label{conjugation} {\it (Conjugation)}
It is useful to note another basic fact that conjugations that come up at various points in the proof do not affect calculations in the Milnor group in our setting. The key point is that the component $l_1$ of the link $L$ (figure \ref{5 components}) is in the $4$th term of the lower central series $(MF_4)^4:=(MF_{m_2,\ldots,m_5})^4$. The same comment applies to each curve ${\gamma}_i$ in the links in figure \ref{Engel links1}.
The Milnor group $MF_4$ is nilpotent of class $4$, that is $(MF_4)^5=\{ 1\}$. All calculations take place in the {\it abelian} group $(MF_4)^4$, so conjugation does not have any effect. In particular, the conjugations (of the commutators $[h_1,\ldots, h_4]$) that appear in the statement of corollary \ref{Engel corollary} may be omitted in our present case. 
\end{remark}

The link $(l_1,l_2',\ldots, l_5')$ defined before remark \ref{orientation conjugation} may be viewed as a band-sum of the $5$-component link $L$ with 12 five-component links contained in disjoint $3$-balls: the component $l_1$ is band-summed with the left-most component of each link as in figure \ref{Torus figure}, and the components $l_2,\ldots,l_5$ are band-summed with the dotted curves and their parallel copies as described above. 
The algebraic outcome is that the original element (\ref{commutator eq}) representing $l_1$ in $MF_{m_2,\ldots,m_5}$ is multiplied by the inverse of the product of the commutators $[h_1,\ldots,h_4]$ that appear in the proof of corollary \ref{Engel corollary} for $g=l_1$. 
This may be seen directly by reading off the element represented by $l_1$ in $M{\pi}_1(S^3\smallsetminus (l'_2\cup \ldots\cup l'_5))$; this is also a special case of the additivity of $\bar\mu$-invariants of links under band-summing \cite{Cochran, K2}.
Therefore $l_1$ is trivial in $M{\pi}_1(S^3\smallsetminus (l'_2\cup \ldots\cup l'_5))$, and the link $(l_1,l'_2,\ldots,l'_5)$ is homotopically trivial.
 Capping it off with disjoint null-homotopies in $D^4_0$ gives the desired disjoint singular slices for the link $L$. 

It is immediate from the construction that condition \ref{standard embedding definition} is satisfied: the slice for $l_1$ does not go over the $2$-handles attached along $K$. (And there is nothing to check  for $A_2,\ldots, A_5$ since they do not have any $1$-handles.) 
Moreover, each constructed map $f_i\co A_i\longrightarrow D^4,$ $i=1,\ldots, 5$, is homotopy standard. Indeed, the singular disks completing the construction of the slices are contained in $D^4_0$ (the original $4$-ball minus collars on the attaching regions), a contractible space. The null-homotopy completing the construction is homotopic within $D^4_0$ to the unknotted disks corresponding to the original embedding $A_i\subset D^4$, and since each $A_i$ has a single $2$-handle there is no disjointess to keep track of.
This concludes the proof of proposition \ref{homotopy standard prop}.
\qed

To complete the proof of theorem \ref{ABslice theorem} for the link $L$ we need to consider all 10 submanifolds $\{ A_i, B_i\}$. Since $B_2, \ldots, B_5$ are collars, they do not affect the embedding problem. The handlebody $B_1$ has $36$ two-handles and no $1$-handles. A Kirby diagram for $B_1$  is obtained from that of $A_1$ (figure \ref{Torus figure})
by performing  zero-framed surgery on the curve $c$ and replacing the dots with zeros. Denote the zero-framed link corresponding to $K_i$ by $\overline K_i$.
Composing null-homotopies in the 12 disjoint solid tori, we see that the resulting link $\overline K:=(\overline K_1\cup \ldots\cup \overline K_{12})$ is null-homotopic (in the sense of Milnor) in the solid torus $\partial^- B_1$. The null-homotopies give rise to a map $(B_1, \partial^- B_1)\longrightarrow (D^4, S^3)$ whose image is contained in a collar on the attaching solid torus, and where the $2$-handles of $B_1$ are mapped in disjointly. The result thus far shows that the link $L$ is {\it link-homotopy A-B slice} (definition \ref{link homotopy ABslice definition}). 

The missing ingredient in establishing the {\it homotopy A-B slice} condition is checking that the constructed map $f\co B_1\longrightarrow D^4$ is homotopy standard, as required in definition \ref{homotopy AB slice}. Since there are no $1$-handles, one has to check only that there is a homotopy from $f$ to the original embedding $B_1\subset D^4$ (corresponding to the decomposition $D^4=A_1\cup B_1$), keeping all $2$-handles disjoint. This original embedding 
consists of disjoint, unknotted $2$-handles attached to $\overline K$. (The link $\overline K\subset( $solid torus $\partial^- B_1)$ is the unlink when considered in the ambient $3$-sphere$=\partial D^4$.) 

\begin{figure}[t]
\includegraphics[width=7cm]{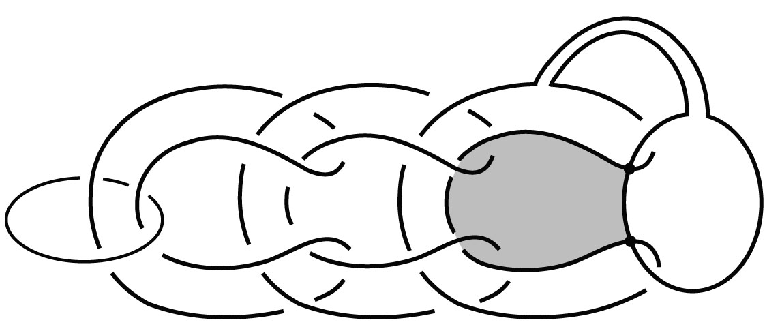}
\put(-205,33){$d$}
{\small 
\put(-157,67){$c_1$}
\put(-110,68){$c_2$}
\put(-72,68){$c_3$}
\put(-65,30){$W$}
\put(0,50){$\overline K_{10}$}
}
\caption{The link $\overline K_{10}$ in the solid torus $S^3\smallsetminus$(neighborhood of $d$). $W$ is a Whitney disk for a null-homotopy.}
\label{WhitneyDisk figure}
\end{figure}
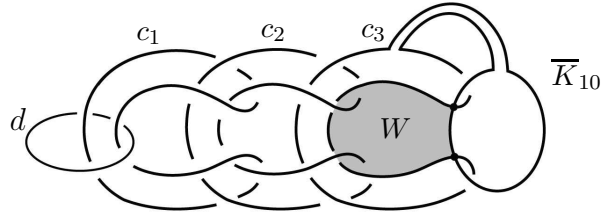

For example, the link $\overline K_{10}$ is shown in figure \ref{WhitneyDisk figure}. This link is considered in the solid torus $=$ complement of  a neighborhood of the curve $d$. 
We remind the reader that to establish the {\it homotopy standard} condition (see the paragraph preceding definition \ref{homotopy AB slice}), one needs to find a homotopy between the  immersed disks bounded by $\overline K_{10}$ in the solid torus $\partial^- B_1$, constructed in the link-homotopy A-B slice solution on one hand, and the standard disks bounded by $\overline K_{10}$ in $D^4$ on the other hand. Moreover, the disks are required to be disjoint during the homotopy. The entire $4$-ball (minus the collar on $\partial B^+_1$) may be used for the homotopy: all other submanifolds $A_i$ and $B_j$, $j\neq 1$ are disregarded here.

The required homotopy will be guided by a suitably chosen Whitney disk.
The obvious Whitney disk $W$ for the null-homotopy of $c_3$, seen in figure \ref{WhitneyDisk figure}, intersects $c_2$. The link $(d, c_1, c_2, c_3)$ is homotopically essential, so there is no Whitney disk disjoint from the rest of the link. However, {\it after $d$ is omitted}, $W$ is isotopic rel boundary to a Whitney disk $W'$ whose interior is disjoint from $c_1,c_2, c_3$. 
Therefore $c_1, c_2$ bound disjoint embedded, unknotted disks $D_1, D_2$ in $D^4$ and $c_3$ bounds a disk $D_3$ with self-intersections which are paired up with a Whitney disk $W'$ whose interior is disjoint from each $D_i$. A Whitney move on $D_3$ along $W'$ gives a homotopy from the constructed map $\coprod D_i\longrightarrow D^4$  to three disjoint unknotted disks bounding $(c_1, c_2, c_3)$. This is a homotopy rel boundary keeping the disks disjoint at all times, as required in the definition of ``homotopy standard''.

Analogous arguments apply to each link $\overline K_i$, and moreover all resulting disks and Whitney disks in $D^4$ are disjoint from each other.
It follows that the map $f\co B_1\longrightarrow D^4$ is homotopy standard, establishing that the link $L$ in figure \ref{5 components} is homotopy A-B slice.

The proof for an arbitrary link $L=(l_1, \ldots, l_n)\in \{\!$Bing(Hopf)$\!\}_5$, defined by Bing doubling without ramification, is directly analogous. 
The decompositions $D^4=A_i\cup B_i$, $i=2, \ldots, n$ 
may be taken to be the trivial decomposition, $2$-handle$\, \cup\, $collar. The decomposition $D^4=A_1\cup B_1$ has one $2$-handle and its $1$-handles are defined by the links $K_i$ whose number and type are determined, as above, by the proof of corollary \ref{Engel corollary}.

Finally, consider an arbitrary $n$-component link $L=(l_1,\ldots, l_n)$ in $\{\!$Bing(Hopf)$\!\}_5$. In this general setting $L$ is not assumed to be almost homotopically trivial, but 
by assumption all $\bar\mu$-invariants of $L$ of length $\leq 4$ vanish. Let $k+1$ be the minimal number of components in $L$ forming a homotopically essential link (by assumption $k\geq 4$). Choose such a $(k+1)$-component sub-link $L'$ of $L$; renumbering the components if necessary, $L'=(l_1,\ldots, l_{k+1})$. Note that $L'$ is almost homotopically trivial.

Based component $l_1$ represents an element of the $k$-th term of the lower central series of $M{\pi}_1(S^3\smallsetminus (l_2\cup\ldots\cup  l_{k+1}))$. Corollary \ref{Engel corollary} expresses it as a product of 
conjugates of the commutators of the form $[h_1,\ldots,h_k]$. 
Recall the  ``elementary Engel links'' in figure \ref{Engel links1}, corresponding to $4$-fold commutators $[h_1,\ldots,h_4]$ in corollary \ref{Engel corollary}.
There are analogous $(k+1)$-component links corresponding to $k$-fold commutators $[h_1\ldots,h_k]$ for any $k\geq 4$. 
Now consider the proof in the $5$-component case above, directly adapted to the $(k+1)$-component link $L'$.
Specifically, for each factor $[h_1\ldots,h_k]$ in the expression for $l_1$, consider a corresponding $k$-component dotted link $K_i$ as in figure \ref{Torus figure}, contained in a solid torus neighborhood of a meridian to $l_1$. Just as in the main body of the proof above, band sums\footnote{The band sums are taken in the complement of the remaining components $l_{k+2}, \ldots, l_n$ of $L$.}
of $l_2, \ldots, l_{k+1}$ with the dotted links $K_i$ yield a homotopically trivial link $(l_1,l'_2, \ldots, l'_{k+1})$. 
(Note that the band sums are taken with a dotted link which by itself (omitting $l_1$) is an unlink. Such band summing may create  {\it higher} non-repeating $\bar\mu$-invariants which involve all of the components $l_1,\ldots, l_{k+1}$, and in addition some of the other components of $L$.)

Now consider the entire resulting link $l_1\cup\overline{L}:= (l_1, l'_2, \ldots, l'_{k+1}, l_{k+2}, \ldots, l_n)$. The argument in the preceding paragraph was used to kill $\bar\mu_{1,2,\ldots,k+1}$. If $l_1$ is trivial in the Milnor group of the rest of the entire link
$M{\pi}_1(S^3\smallsetminus \overline{L})$, perform link-homotopies on $\overline{L}$ and contract $l_1$ in the complement. If $l_1$ is part of another homotopically essentially sub-link, pick such a homotopically essential, almost trivial $l$-component sub-link
($l\geq k+1$) and repeat the argument from the previous paragraph.

With $l_1$ gone, consider the remaining $(n-1)$-component link $\widetilde L$. Since the only operations applied to $l_2,\ldots, l_n$ to form $\widetilde L$ were band-sums with unlinks and link-homotopy, the non-repeating $\bar\mu$-invariants of $\widetilde L$ of length $\leq 4$ are trivial. The same argument as above now applies to $\widetilde L$.
Proceeding by induction, the link is eventually reduced to a $4$-component link. The condition on non-repeating $\bar\mu$-invariants ensures that it is homotopically trivial.

The number and types of elementary Engel links that come up in this process for $L$ are determined by the proof of lemma \ref{Engel nilpotent}.
Each of the elementary links has a {\it null-homotopic} counterpart, illustrated in figure \ref{Engel links}. 
The decompositions $D^4=A_i\cup B_i, i=1,\ldots, n$ are then defined analogously to figure \ref{Torus figure}: each $A_i$ has a single zero-framed $2$-handle, and the $\{ A_i\}$ 
incorporate all dotted links that are needed in the algebraic argument.
The proof of proposition \ref{homotopy standard prop} then goes through to give disjoint, homotopy standard maps of $A_1,\ldots, A_n$ into $D^4$. The proof that the submanifolds $B_i$ admit 
homotopy standard maps into the collar is directly analogous to that in the almost trivial case considered above.
\qed
\begin{remark}
The proof of theorem \ref{ABslice theorem} essentially relied on the asymmetry of the roles played by the $1$- and $2$-handles of the submanifolds $A_i, B_i$. Specifically, many parallel copies of each dual $1$-handle core are used in the relative slicing (in the proof of Proposition \ref{homotopy standard prop} the link $J$ is band-summed into the components of $K$ {\it and} their parallel copies), whereas a single copy of each $2$-handle must be mapped disjointly.
This should be compared with the statement of the Round Handle Problem in section \ref{s-cobordism section} where the analogous links are in fact symmetric, see remark \ref{comparison remark}.
\end{remark}

\subsection{A link-homotopy solution for the Borromean rings} \label{Borromean homotopy solution}
In light of proposition \ref{higher prop}, in the proof of theorem \ref{ABslice theorem} it sufficed to consider links in $\{\!$Bing(Hopf)$\!\}_5$.
The Borromean rings (Bor)  is the simplest and best known example of a link in $\{\!$Bing(Hopf)$\!\}$, however since Bor is not in $\{\!$Bing(Hopf)$\!\}_5$,
the proof did not apply directly to this link.
We sketch a modification needed to give a link-homotopy solution of the A-B slice problem 
for Bor.

\begin{figure}[ht]
\includegraphics[width=13.5cm]{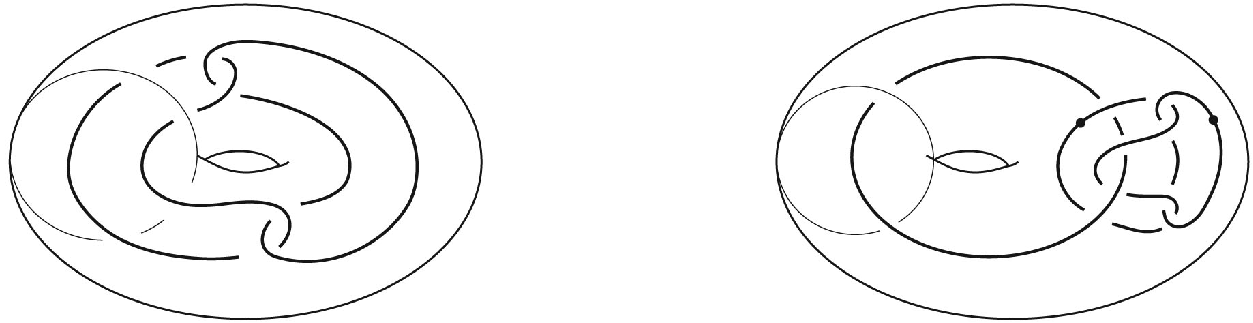}
\put(-386,5){$A_i$}
\put(-155,5){$B_i$}
{\small
    \put(-336,8){$0$}
    \put(-290,7){$0$}
  \put(-82,9){$0$}}
 \caption{$D^4=A_i\cup B_i$, $i=2,3$.}
\label{torus Kirby}
\end{figure}
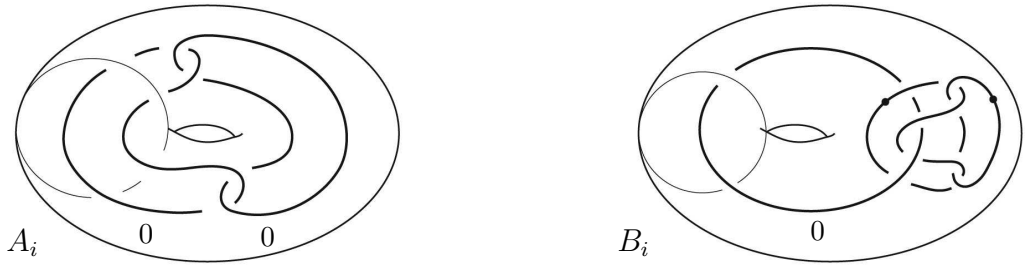

Consider the decomposition $D^4=A_1\cup B_1$ as in the proof of theorem \ref{ABslice theorem}. For $i=2,3$ 
let $D^4=A_i\cup B_i$ where 
$B_i=T\times D^2$ is a thickening of the torus $T$ with a single boundary component, embedded in a standard way in $D^4$. Then  $A_i$ is obtained from the collar on the solid torus ${\partial}^+ A_i$ by attaching two zero-framed $2$-handles to the Bing double of the core, figure \ref{torus Kirby}.

The link in figure \ref{5 components} is obtained from Bor by Bing doubling two of the components.
This Bing doubling is incorporated in the definition of $A_2, A_3$, so proposition \ref{homotopy standard prop} 
applies directly to yield a  link-homotopy solution for $A_1, A_2, A_3$: disjoint maps $f_i\co A_i\longrightarrow D^4$ so that
the attaching curves $\{ {\alpha}_i\}_{i=1,2,3}$ form the Borromean rings, all $2$-handles of $\{A_i\}$ are mapped in
disjointly, and condition \ref{standard embedding definition} is satisfied.

\begin{figure}[ht]
\includegraphics[width=5cm]{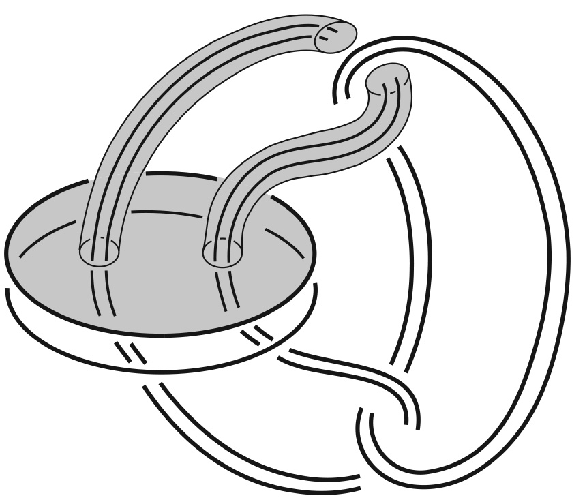}
{\scriptsize
\put(-32,17){${\alpha}_3$}
\put(-83,13){${\beta}_1$}
\put(-103,5){${\alpha}_1$}
\put(-150,70){${\beta}_2$}
\put(-150,39){${\alpha}_2$}
\put(-15,5){${\beta}_3$}
}
{\small
\put(-123,102){$B_2$}
}
\caption{}
\label{WhitenyDisk figure}
\end{figure}

Since each component of Bor bounds a genus one surface in the complement of the other components (figure \ref{WhitenyDisk figure}), 
the submanifolds $B_2, B_3$ admit disjoint embeddings in the complement of $\coprod f_i(A_i)$. Since $B_1$ is a collar, it does not affect the embedding problem.
This completes a link-homotopy solution for Bor.

\newpage

\section{The Round Handle Problem, $5$-dimensional $s$-cobordisms, \\ and general doubles.} \label{s-cobordism section}
\subsection{}
We call attention to an explicit construction of a smooth $4$-manifold $M$ with 
$\partial M\cong {\mathcal S}_0(Wh(L))$, the zero-framed surgery on the Whitehead double\footnote{As seen in figure \ref{Whitehead double} $Wh(L)$ is not defined until the clasp signs $\pm 1$ on each component are specified. We drop this detail from our notation, but point out where in the construction of $M$ these signs are seen. Assertions about $Wh(L)$ apply to {\it all} sign choices.} of an initial $k$-component link $L\subset S^3$. Whitehead doubling replaces each component $l_i$ of $L$ with a satellite:

\begin{figure}[ht]
\includegraphics[height=2cm]{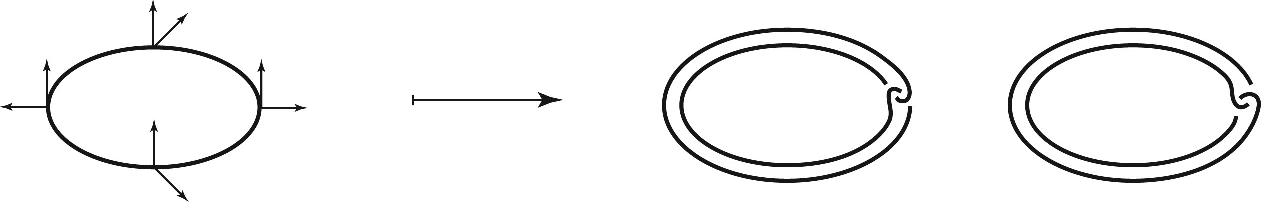}
{\small
    \put(-348,-11){tube around $l_i$}
    \put(-348,-21){with zero framing}
    \put(-140,-20){the $+$, $-$ (resp.) $Wh(l_i)$}}
 \caption{Whitehead doubling}
\label{Whitehead double}
\end{figure}

$M$ is the source of an $F_k$ (free group) - surgery problem, rel boundary, with target $\natural_{i=1}^k S^1\times D^3$. The problem has vanishing (Wall) surgery obstruction iff all the linking numbers $(l_i, l_j)=0$, $1\leq i, j\leq k$. Solving this surgery problem constructs a slice complement for $Wh(L)$. The well-known universal surgery problems \cite{F1} arise when $L$ is some ramified Bing double of the Hopf link: 
$\HopfLink$\hspace{-2.5em}. It is routine to build the maps and cover them with required normal data, so we will only describe the construction of $M$. 

\begin{figure}[ht]
\includegraphics[height=5cm]{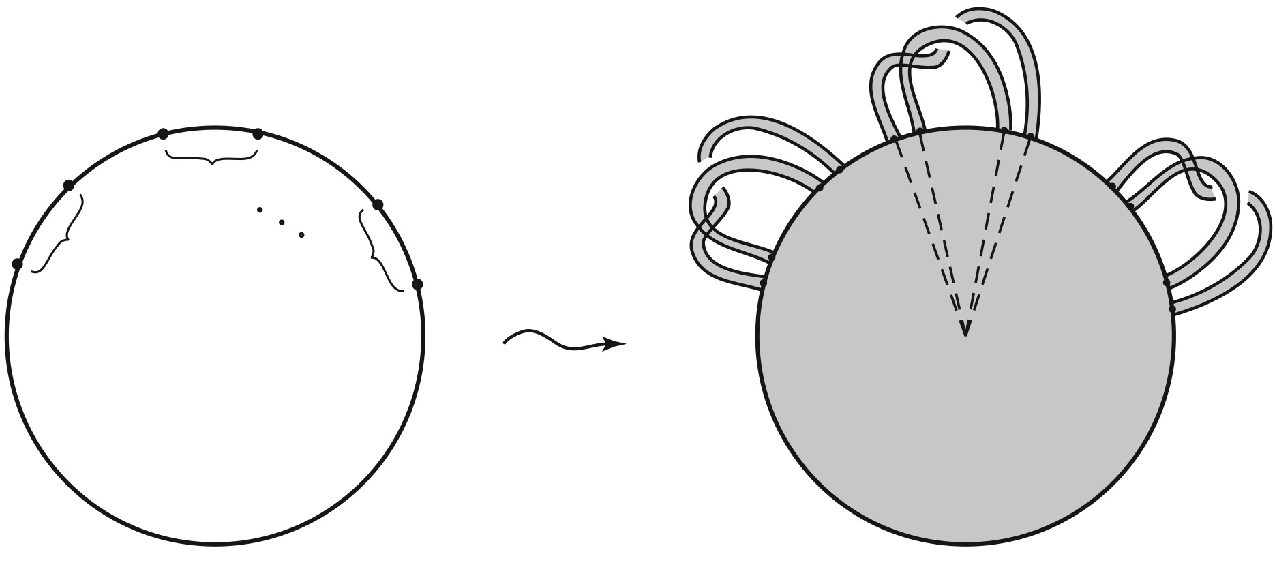}
{\small
    \put(-294,77){$l_1$}
    \put(-264,91){$l_2$}
\put(-233,73){$l_k$}
\put(-266,40){$D^4$}
    \put(-85,35){$M$}}
 \caption{Schematic of $M$}
\label{M figure}
\end{figure}
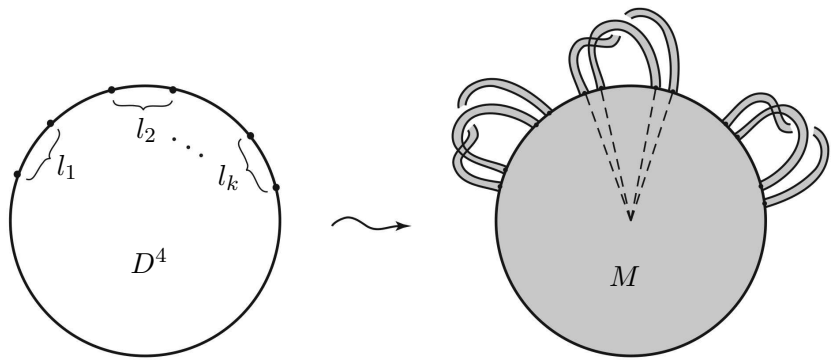

By definition $M=M(L)$ is obtained (see figure \ref{M figure}) by attaching $k$ pairs of plumbed $2$-handles to $D^4$. The attaching circles are $2L$, i.e. the link $L$ and a parallel copy with all framings equal to zero. There is a sign choice, $\pm$, at each plumbing point. This is the place where the sign of the clasp in the Whitehead double is determined and will not be commented on again. 

\begin{lemma} \label{M lemma} 
$\partial M\cong {\mathcal S}_0(Wh(L))$.
\end{lemma}

The proof is given below (following definition \ref{RHP definition}) and repeats \cite{FT2}.
Henceforth assume all linking numbers $(l_i, l_j)=0$. From this we see $k$ hyperbolic pairs (one displayed using dotted lines in figure
\ref{M figure}) over ${\mathbb Z}[F_k]$, i.e. $2k$ spherical classes of the form 
($2$-handle core)$\cup $cone to origin ($\partial$core), representing $\oplus_k \bigl(\begin{smallmatrix}
0&1\\ 1&0 \end{smallmatrix}\bigr)$. $F_k$ is the free fundamental group of $M$ generated by the plumbings. The nonsingularity of this form is
equivalent to the natural map ${\alpha}\co \partial M\longrightarrow \sharp_k S^1\times S^2$ being a ${\mathbb Z}[F_k]$-homology isomorphism.

A consequence of the (still open) topological surgery conjecture is that there exist a topological $4$-manifold $N$, $\partial N = \partial M$, with a homotopy equivalence $\beta$ extending $\alpha$:
$${\beta}\co (N, \partial) \longrightarrow (\natural_k \, S^1\times D^3, \sharp_k \, S^1\times S^2).$$

The entire thrust of the A-B slice discussion was to find a way of contradicting the existence of $N$ using ``low-tech''\footnote{It is an interesting question whether there is a ``high-tech'' approach to refuting the surgery conjecture. We know no useful reformulation in gauge theory.}
 nilpotent invariants of $L$. The philosophy was that  $Wh(L)$ has little to grab onto, certainly no nilpotent invariants so it was preferable  to ``undouble'' the problem and work directly with $L$ where there are $\bar\mu$ invariants to work with. In this section we describe a variant of this approach which we call the ``round handle problem'' (RHP). Like the A-B slice problem RHP can be translated into a question about slicing some ``stabilized'' version of $L$. The advantage of the RHP variant is that the stabilization is better controlled - the ramification of dotted (red) curves featured in the homotopy A-B slicing (section \ref{homotopy solution}) cannot occur in the RHP context.
 Thus even non-repeating $\bar\mu$-invariants, associated with GBR, might possibly be used to formulate an obstruction.
The ``disadvantage'' of the new context is that such an obstruction to ``stable slicing''  contradicts only  the logical union ($4$D surgery conjecture $\!\wedge\!$ $5$D s-cobordism conjecture) and we would not gain any information on which fails\footnote{However see remark \ref{aside}.} - merely that something goes seriously wrong. A statement of these two famous conjectures are given concisely in a footnote.
\footnote{{\it Surgery conjecture}: Any degree one normal map $(M, \partial M)\longrightarrow (X, \partial X)$ from a $4$-dimensional topological manifold to a Poincar\'{e} pair which is a ${\mathbb Z}[{\pi}_1 X]$-homology equivalence over $\partial X$ is topologically normally bordant to a homotopy equivalence iff the Wall obstruction in $L_4[{\pi}_1 X]$ vanishes.\\
{\it $s$-cobordism conjecture}:  Assume $(W; M_1, M_2)$ is a compact, topological $5$-dimensional $s$-cobordism which is a product over the boundary (i.e. $\partial W\smallsetminus {\rm int}(M_1\cup M_2)\cong \partial M_1\times I\cong \partial M_2\times I$, $\cong$ meaning {\it homeomorphic}.) Then $W$ itself admits a compatible product structure $W\cong M_1\times I \cong M_2\times I$.}

\begin{remark} \label{aside} {\it Aside on proper s-cobordism}.
It is an old observation (see the next paragraph) that the ``proper $5$D s-cobordism'' conjecture implies both $4$D surgery and $5$D s-cobordism conjectures. 
So, if one insists, a specific failure could be pointed to (if the RHP has no solution).

The proper or p-s-cobordism theorem was established by L. Taylor (Ph.D. thesis, UC Berkeley, 1972) for p-s-cobordisms of dimension $6$ and higher.  It is an open question whether his thesis result extends in the topological category  to dimension $5$. The algebraic setting for the general obstructions is complicated a bit by properness but the case of greatest interest is when the global fundamental group is free and the fundamental group of the end also (the same) free group. In \cite{F01} a variant of the low dimensional  surgery  sequence is established. It is straightforward that a successful extension of Taylor's thesis  (in the above free case) would convert the published  variant to the full topological surgery sequence in these dimensions. It was considered so unlikely that this extension existed when \cite{F01} was written that this point is not explicitly made in the paper.
\end{remark}

\begin{figure}[ht]
\vspace{.3cm}
\includegraphics[height=3cm]{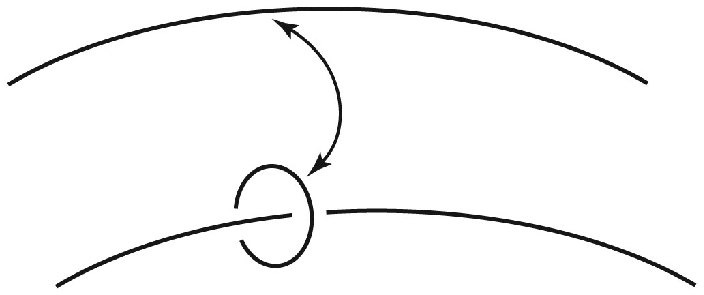}
{\small
    \put(-216,60){$l'_i$}
    \put(-201,-1){$l_i$}
\put(-116,3){$m_i$}
    \put(-103,54){$r_i$}}
\vspace{.3cm}
 \caption{Attaching a round handle}
\label{round handle figure}
\end{figure}
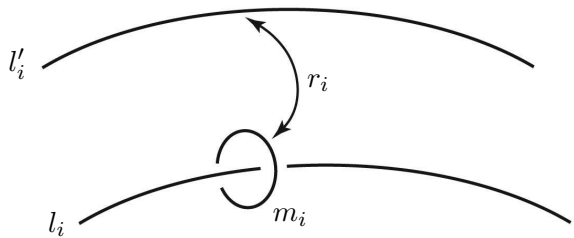

Again given $L$ define $R=D^4\cup k$ round $1$-handles. In this dimension a round $1$-handle is $(D^1\times D^2\times S^1, S^0\times D^2\times S^1)$. $R$ is built by attaching the $i^{\text th}$ round $1$-handle $r_i$ to a meridian $m_i$ to $l_i$ {\it and} a parallel copy $l'_i$ of $l_i$ (lying beyond the meridian), figure \ref{round handle figure}. 

Suppose the link $L$, considered as lying in $\partial R$, is slice in $R$, meaning $L$ bounds $k$ disjoint, topologically flat $2$-disks in $R$, equivalently bounds $k$ disjoint topological $2$-handles in $R$. Let $T$ be the ``slice complement'', i.e. the manifold with boundary obtained by deleting the interior of those $k$ $2$-handles, $T=R\smallsetminus {\text{int}} (\coprod_k 2$-handles$)$.
The proofs of the following two lemmas are postponed until after the definition of the Round Handle problem (definition \ref{RHP definition}).

\begin{lemma} \label{T lemma}
$\partial T\cong \partial M$.
\end{lemma}

Actually $N$ is a candidate for the slice complement $T$. By this we mean, if $N$ exists we can reconstruct a manifold $R'$ very much like $R$ by attaching $k$ $2$-handles to $N$.

\begin{lemma} \label{s-cobordism lemma} \sl
If $N$ exists we can form $R'=N\cup k$ ($2$-handles) so that there exists a $5$D s-cobordism $W$, which is a smooth product over the boundary, joining $R'$ to $R$.
\end{lemma}

Thus if we assume $4$D surgery and $5$D s-cobordism conjectures (we call this package the surgery sequence conjecture (SSC)) then $L$ is slice in $R$. The slice disks $S$ may be taken to be topologically transverse \cite{FQ} to the $k$ cocores $(D^2\times S^1)_i$ of $r_i$, the round $1$-handles. Cutting $R$ open along the cocores recovers the $4$-ball $D^4$ with the promised ``stabilization'' $\widehat L$ of $L$, i.e. two copies of the $1$-manifolds =$(\text{slice disks})\cap\; \coprod_{i=1}^k(D^2\times S^1)_i.$ The components of $\widehat L$ are now seen to co-bound some disconnected planar surface  $P$ made from fragments of the slice disks. To summarize: $L=\partial S\subset R$ yields $\widehat L=\partial P\subset D^4$. One may understand the combinatorial possibilities for $\widehat L$ and $P$ and attempt to see if any are compatible with known properties of $\bar\mu$-invariants. The chief feature of $\widehat L$ is that $\widehat L=L\cup Q\cup \widetilde Q$ where $Q$ is contained in the meridional  solid tori $M_i$ with core circles $m_i$, and $\widetilde Q$ is identical to $Q$ but transported by the zero-frame preserving homeomorphism from $M_i$ to $L'_i$, the parallel solid torus with core $l'_i$ (see figure \ref{round handle figure}). In practice the components of $Q$ ``help'' with the existence of $P$ by canceling the $\bar\mu$-obstructions but they beget harmful $\widetilde Q$ with new $\bar\mu$-obstruction to bounding $P$.

\begin{definition} \label{RHP definition}
The Round handle problem RHP is to determine whether any link $L$ with vanishing linking numbers is slice (bounds disjoint topologically flat disks) in the corresponding manifold $R$. A contradictory possibility is that some
non-trivial $\bar\mu$ of $L$ survive in all possible stabilization processes to prevent any $\widehat L$ bounding $P$ as described above.
\end{definition}

\begin{figure}[ht]
\includegraphics[height=2.9cm]{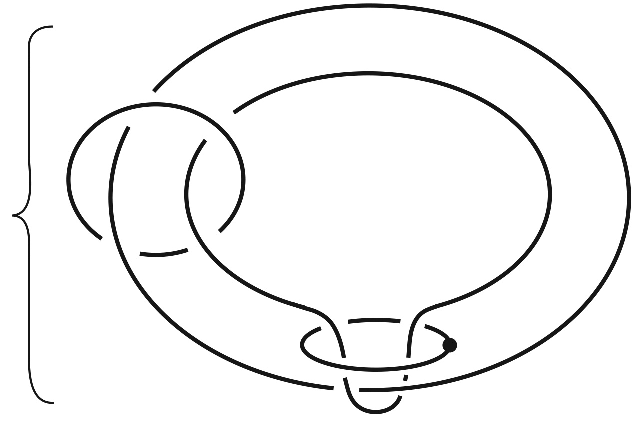} \hspace{2cm} 
\includegraphics[height=2.9cm]{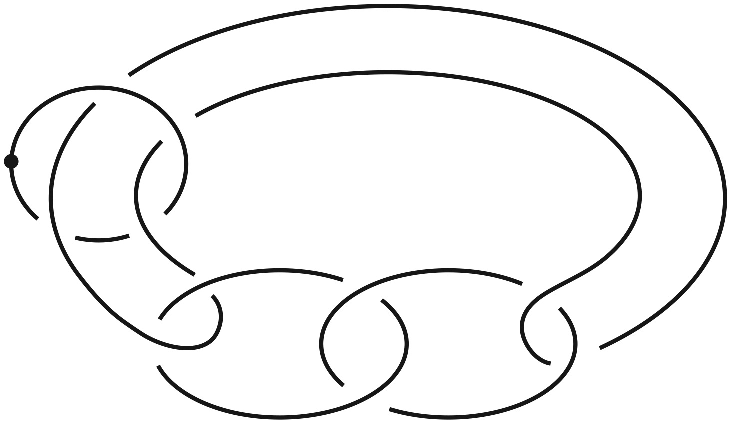}
\put(-362,42){leads}
\put(-362,29){to $M$}
{\large \put(-181,41){$\cong$}}
\put(-321,61){$d$}
{\small \put(-220,11){$0$}
\put(-232,21){$0$}
\put(-90,33){$0$}
\put(-58,33){$0$}}
\put(-7,65){$d$}
 \caption{}
\label{calc1}
\end{figure}

\begin{remark} \label{comparison remark}
We summarize the difference between the link stabilization formulations of the A-B slice problem and of the RHP. 
They have similar set-ups: for each link component $l_i$ of $L$ both problems consider two solid tori $M_i, L'_i$, neighborhoods of a meridian $m_i$ to $l_i$ and of a parallel copy $l'_i$, compare
figures \ref{RelSlice figure}, \ref{round handle figure}. 
To begin with, in both setups there are identical ``stabilization'' links in the solid tori $\{ M_i, L'_i \}$, and the question is whether the link $L$ bounds planar surfaces in $D^4$ whose other boundary components correspond to the stabilization links. The distinction between the two is that in the AB slice problem one is allowed to take an arbitrary number of parallel copies of the ``helping'' red curves (corresponding to the passage of the slices over the $2$-handles attached along these curves), while their counterpart curves in the dual solid tori do not have to be ramified. In the RHP the curves in $M_i,L'_i$ match precisely, corresponding to the passage of the slices over the round handles.
The proof of theorem \ref{ABslice theorem} in section \ref{homotopy solution}, using the $2$-Engel relation, crucially relies on taking parallel copies, so it does not go through in the RHP setting. This is the basis for our comment in the introduction: after stabilization non-repeating $\bar\mu$-invariants may reappear as repeating; there is no sharp dichotomy. Embedded disks are stronger than disjoint maps of disks, and intermediate is disjoint maps on {\it many} link parallels. It seems further progress on the homotopy approach to the AB slice problem would require studying parallels, i.e. repeating $\bar\mu$-invariants.
Definition \ref{RHP definition} is motivated by another line of thought: perhaps if one is satisfied with an obstruction to the logical union (surgery $\vee$ $s$-cobordism) then one might still work in the simpler world of non-repeating $\bar\mu$-invariants.
\end{remark}
{\it Proof of Lemma \ref{M lemma}}. We use the usual conventions (cf. \cite{GS}) for handle diagrams (Kirby calculus). The argument is a local handle computation inside each of the $k$ solid torus neighborhoods  $L_i$ of $l_i$. We draw dual circle marked $d$ to define - via its complement - the solid torus $L_i$.
The plumbed pair is diagrammed as two zero-framed $2$-handles and one $1$-handle (circle with dot) in $L_i$, figure \ref{calc1}. 

\begin{figure}[h]
\includegraphics[height=3.3cm]{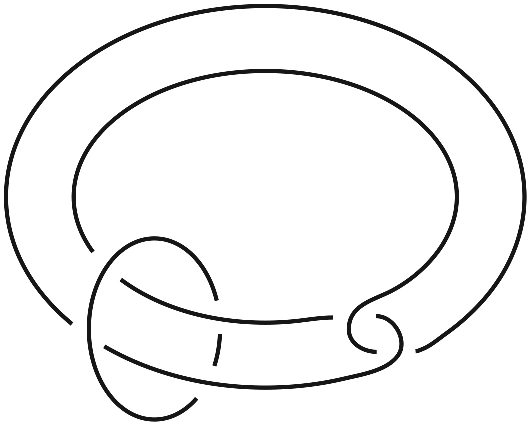} \hspace{2cm} \includegraphics[height=3.3cm]{calc2.eps}
{\large \put(-335,46){$\overset{\partial}{\cong}$}
\put(-155,46){$\cong$}}
\put(-187,70){$d$}
\put(-291,0){$0$}
\put(-4,70){$0$}
\put(-107,0){$d$}
 \caption{}
\label{calc2}
\end{figure}

As far as the boundary is concerned, we may cancel the hyperbolic pair and replace the dot with a zero, obtaining figure \ref{calc2}.

Notice that the calculation did {\it not} assume unknottedness of $l_i$, it takes place in the solid torus $L_i$, re-embedded as unknotted for convenience only. \qed

\begin{figure}[ht]
\includegraphics[height=3.2cm]{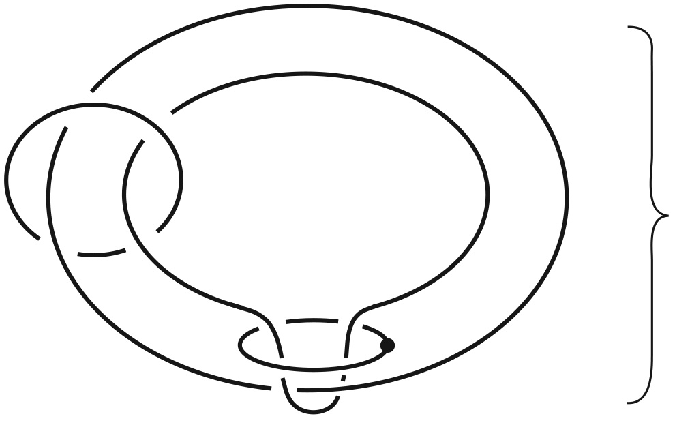} 
\put(3,48){leads}
\put(3,35){to $R$}
\put(-153,62){$d$}
\put(-60,57){$l_i$}
{\small \put(-42,81){$0$}}
 \caption{}
\label{calc3}
\end{figure}

{\it Proof of lemma} \ref{T lemma}. Again the calculation can be localized to the solid tori $L_i$. The round handle attachment indicated in figure \ref{round handle figure} is shown in terms of ordinary handles in figure \ref{calc3}.

Slicing $l_i$ has the effect, in so far as the boundary is concerned, of placing either a dot or a zero on that component, returning us to the first panel of figure \ref{calc1}, i.e. a diagram of $\partial M$. \qed

{\it Proof of lemma \ref{s-cobordism lemma}.}
We define a diffeomorphism $d\co \partial R\longrightarrow \partial R'$ by starting with: $\partial T\cong \partial M \cong \partial N$ and propagating the attachments of the $2$-handles, $h_1,\ldots, h_k$, of $R\smallsetminus T$ from $\partial T$ to $\partial N$.  The first step is to extend this to a simple homotopy equivalence
$h\co (R, \partial R)\longrightarrow (R', \partial R')$. $R'$ is homotopy equivalent to 
$(\vee_k S^1)\vee (\vee_k S^2)$, so the only possible obstruction to extending $d$ as a map is $\mathcal{O}\in H^3(R, \partial R; {\pi}_2 R')$.

The boundary of a basis of these relative $3$-cocycles  are the $2$-sphere $S_i$
factors, $1\leq i\leq k$, in $\partial R\cong {\mathcal S}_0(L)\sharp(\sharp_k S^1\times S^2)$. $N$ has $k$ distinguished tori $T_i$, the natural genus one Seifert surface for
$Wh(l_i)$ capped off by surgery on $l_i$. In $R'=N\cup k$ $2$-handles, the $2$-handles $h_i$ surger (along a copy of $l_i$) $T_i$ into a $2$-sphere $S'_i$ and $d(S_i)=S'_i$. We use the notation: 
$$\partial (T_i\times [0,1]\cup h^-_i)=T_i\times 0\cup S'_i, \;\, {\text{where}} $$
$$h^-_i=(D^2\times D^1, \partial D^2\times D^1)\subset (D^2\times D^1\times D^1,\partial D^2\times D^1\times D^1)=h_i.$$ 
Since each $T_i$ is null-homotopic in $N$, $S'_i$ bounds a singular $3$-ball $b_i$ in $N$ lying in $(T_i\times [0,1]\cup h^-_i)\cup_{T_i\times 0} (\text{Cone}(T_i))$, showing 
${\mathcal{O}}=0$. 

Since $R'$ is also homotopy equivalent to $(\vee_k S^1)\vee (\vee_k S^2)$ and $d\co \partial R\longrightarrow \partial R'$ lines up the generators bijectively, 
$h$ is automatically a homotopy equivalence. Since $Wh_1(F_k)=0$ \cite{Stallings}, 
$h$ is also a simple homotopy equivalence.

Covering each map with (arbitrary) normal data we obtain two structures on $(R', \partial)$:
$$\text{id}\co (R', \partial)\longrightarrow (R', \partial), \; \; \text{and} \; \; 
h\co (R, \partial)\longrightarrow (R', \partial).$$
The possible obstructions to a relative normal cobordism between $id$ and $h$ lies in 
$[(R', \partial), (G/TOP, *)]$. The Postnikov tower 
for $G/TOP$ begins with a $K({\mathbb Z}_2, 2)$ which detects Arf invariants for the possible $2$ dimensional splitting problems, and then a $K({\mathbb Z}, 4)$  which detects signature differences for the possible $4$ dimensional splitting problems \cite{KS}. Let us consider the obstructions to making $h$ normally bordant to $id$. 
The first obstruction is $$o_2\in H^2(R', \partial ; {\pi}_2(K({\mathbb Z}_2, 2)))\cong Hom(H_2(R', \partial ; {\mathbb Z}_2), {\mathbb Z}_2).$$ 
The co-cores $c_1,\ldots, c_k$ of $h_1,\ldots, h_k$ (under the identification $\partial T\cong \partial N$) are a basis for $H_2(R', \partial ; {\mathbb Z})$ and $o_2([h_i])$ is $({\rm Arf}\, h^{-1}(c_i) - {\rm Arf}\, (c_i)) = (0-0)= 0$ since $h^{-1}(c_i)\subset R$ is itself homologous to the co-core $c_i$, again identifying $\partial N\cong \partial T$.
The second obstruction may be identified as signature$(R)-$signature$(R')$=0, so it vanishes as well.
Thus $id$ and $h$ are normally cobordant, rel boundary, via $W^5$. According to Wall \cite{Wall} there is a surgery obstruction ${\sigma}(W)\in 
L^s_5(F_k)\cong L^h_5(\{ e\})\oplus_{i=1}^k L^h_4(\{ e\})$, by the splitting principle \cite{Cappell}. 
$L^h_5(\{ e\})\cong 0$ and $L^h_4(\{ e\})_i\cong 8{\mathbb Z}$ given by the signature of the spin $4$-manifold dual 
to each free generator. It is possible to modify our choice of $W$ to $W'$ to kill these $k$ surgery obstructions. To change the $i$th obstruction by $\pm 8$ replace an embedded $S^1\times D^4$ parallel to the $i$th free group generator with an embedded $S^1\times (E^8$-manifold$\smallsetminus \dot{D}^4)$. Knowing ${\sigma}(W')=0\in L^s_5(F_k)$, $W'$ is normally cobordant, rel its boundary, to an s-cobordism $W''$ from $R'$ to $R$. \qed

\subsection{General doubles} \label{general double section}
In the early days of $4$-manifold topology decomposition theoretic properties of Whitehead doubling played a key role \cite{F0}. But in the current study  of non-simply connected surgery we are completely divorced (and perhaps it is our loss) from point set topology so it becomes a hindrance to adhere to the literal meaning of ``Whitehead double''. We suggest a more algebraic generalization (in fact two) which will be exploited in section \ref{model section}. Surface genus is a natural parameter in the generalization and we will only exploit  the genus one case so the reader may restrict the definition below accordingly. The purpose of Whitehead doubling, from our current perspective is to weaken a link $L$ so that the problem of slicing its replacement $Wh(L)$  can be expressed as an (unobstructed) $4$D surgery problem (which if solved would produce a candidate manifold for the link slice complement.) The construction of $M$ (see figure \ref{M figure}) illustrates this strategy. But if this is all we want then we may define a ``general double'' as follows.

Consider a (usually) disconnected surface $S=\coprod_{i=1}^k S_i\subset S^3$, with $\partial S_i$ a simple closed curve, $\partial S=: K$ is a $k$-component link. The case used in section \ref{model section} is where genus$(S_i)=g_i=1$.
We assume that in some basis of simple closed curves on $S$ the Seifert form is:
$$\bigoplus _{j=1}^J \;\, \;\,
\begin{array}{c|cc|}
\multicolumn{3}{r}{{\scriptstyle{x_j}}\;\;\;{\scriptstyle{ y_j}}}\\
 {\scriptstyle{x_j}} & 0 & \pm1 \\
{\scriptstyle{y_j}} & 0 & 0  \\ 
\multicolumn{3}{c}{}
\end{array}\; ,
\hskip 1.1cm
J=\sum_{i=1}^k g_i.
$$

This makes $K$ a ``good boundary link'' (GBL) \cite{F01}. When $L$ has vanishing linking numbers, $K=Wh(L)$ has an obvious Seifert surface showing it is of this type. Good boundary links are known to admit unobstructed surgery problems for constructing a slice complement. (It is not known, in general, when these surgery problems have topological solutions.)

Assume $g_i=1$, then each surface $S_i$ is a $\pm$ plumbing of two untwisted bands which when pushed off according to the lower left Seifert matrix entry, are disjoint and with linking number $<\! x_i, y_i^- \! >=0$. ($<\! x_i, y_i^+ \! >=\pm 1$, the upper right entry of the Seifert matrix.)

\begin{definition} \label{genus1} ({\it Generalized double, genus one case} ($g_i=1$))
With the notation above, we say a GBL $K$ admitting a Seifert surface $S$ where each component $S_i$ has genus $g(S_i)=1$ is the generalized double of the set of disjoint simple closed curve pairs $\{ (x_1,y_1^-), \ldots, (x_J, y_J^-) \}$. (In the case of unramified Whitehead double $Wh(K), K=\{ K_1, \ldots, K_J\}, x_i=K_i$ and $y_i^-\,=\,0$-linking parallel$(K_i), 1\leq i\leq J$.)
\end{definition}

\begin{lemma} \label{general double lemma} \sl
${\mathcal S}_0(K)\cong \partial \overline{M}$ where $\overline M$ is obtained from $D^4$ by attaching $k$ $\pm$plumbed pairs of $2$-handles. $\overline M$ is the source of an unobstructed surgery problem for building a slice complement for $K$. Each pair is attached to zero framed $x_i$ and $y^-_i$.
\end{lemma}

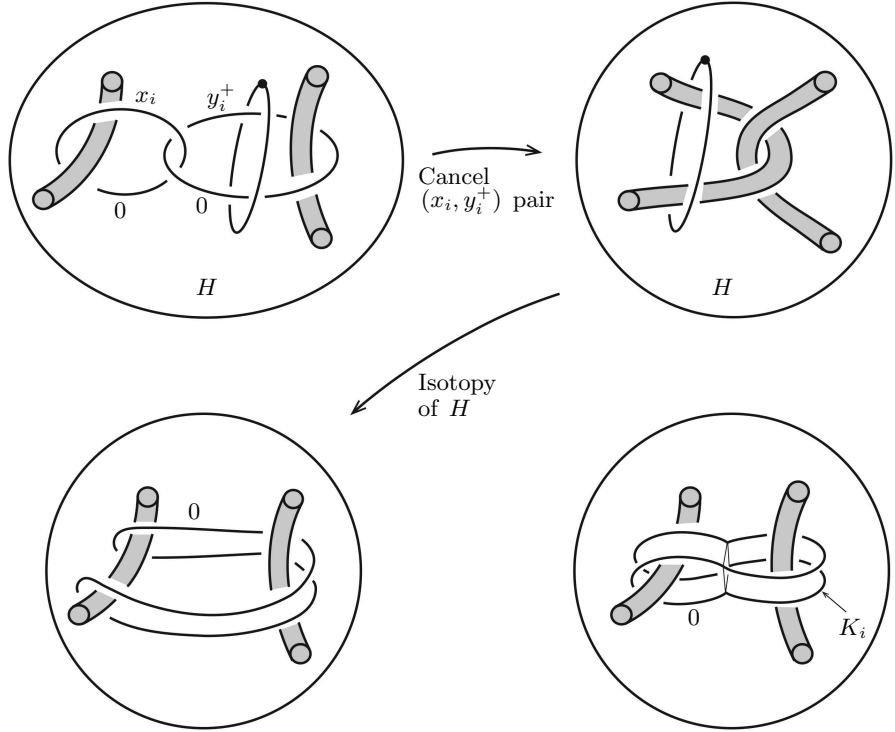
\begin{figure}[ht]
\includegraphics[height=10cm]{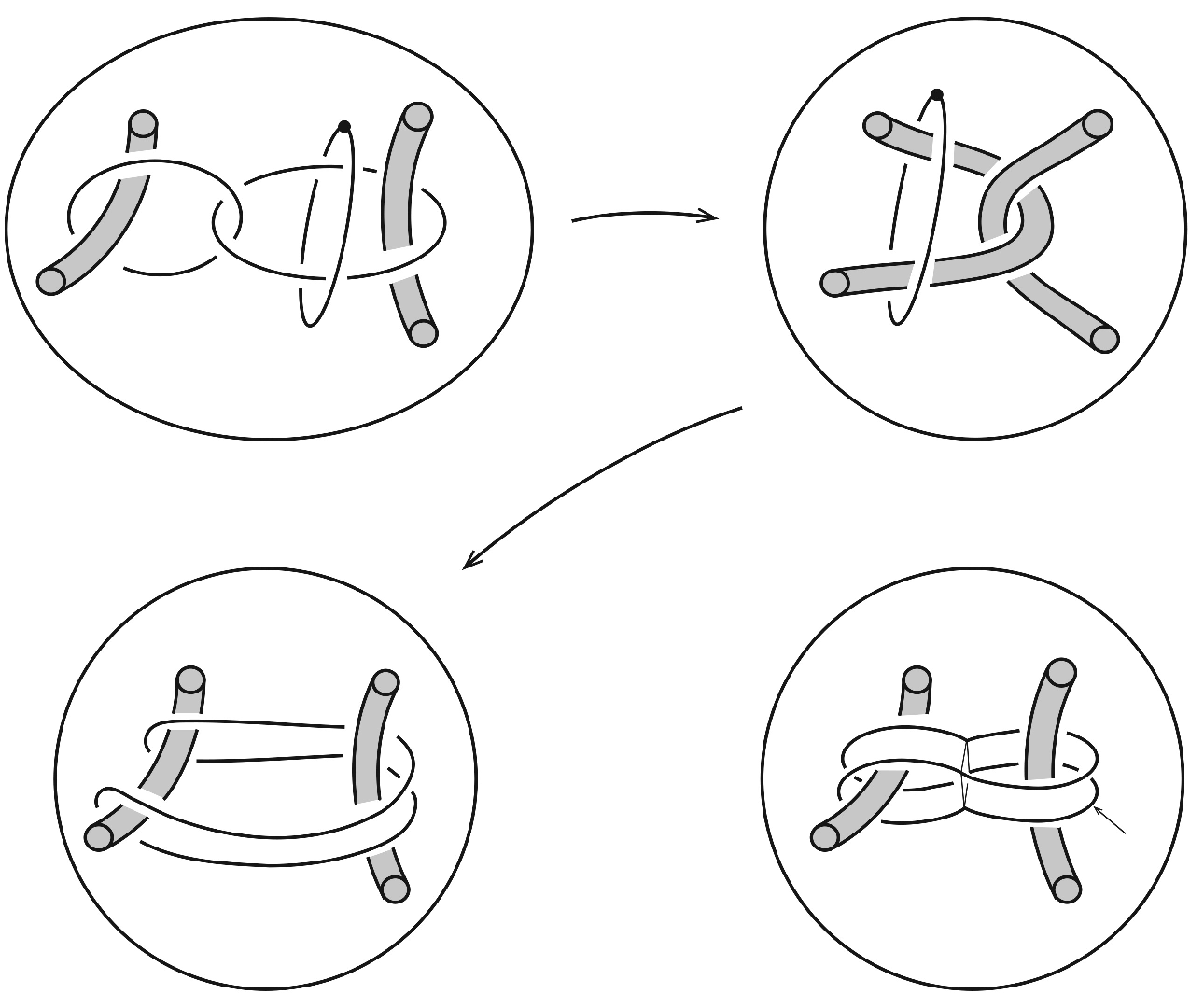} 
{\scriptsize \put(-180,211){Cancel}
                 \put(-180,202){$(x_i, y^+_i)$ pair}
\put(-180,133){Isotopy}
\put(-180,123){of $H$}
\put(-265,169){$H$}
\put(-70,169){$H$}
\put(-79,44){$0$}
\put(-288,242){$x_i$}
\put(-261,241){$y^+_i$}
\put(-296,198){$0$}
\put(-265,200){$0$}
\put(-268,84){$0$}
\put(-22,40){$K_i$}
}
 \caption{Shaded tubes are {\it removed} from their balls, creating genus.}
\label{general Kirby figure}
\end{figure}

{\it Proof of lemma \ref{general double lemma}.} The schematic for $\overline M$ is very similar to figure \ref{M figure}, except now $x_i, y^-_i$ are not necessarily $l_i$ and a parallel copy $l'_i$ as implied in that figure. As in the proof of lemma \ref{M lemma} we can localize the Kirby calculation, this time to a genus two handlebody $H$ (rather than a solid torus) containing a component $K_i$ of $K$. We show initial and final Kirby diagrams in $H$ are boundary equivalent, see figure \ref{general Kirby figure}.

The statement about surgery is immediate: ${\pi}_2$ is a free module over ${\mathbb Z}[F_k]$ and the intersection form is manifestly hyperbolic. \qed

{\it Note.} The proof shows that the diffeomorphism type of ${\mathcal S}_0(K)$ does not depend on the choice of arcs which join $x_i$ to $y^-_i$ and define the plumbing, similarly for $\overline M$ in lemma \ref{general double lemma}.

\begin{figure}[ht]
\includegraphics[height=18.2cm]{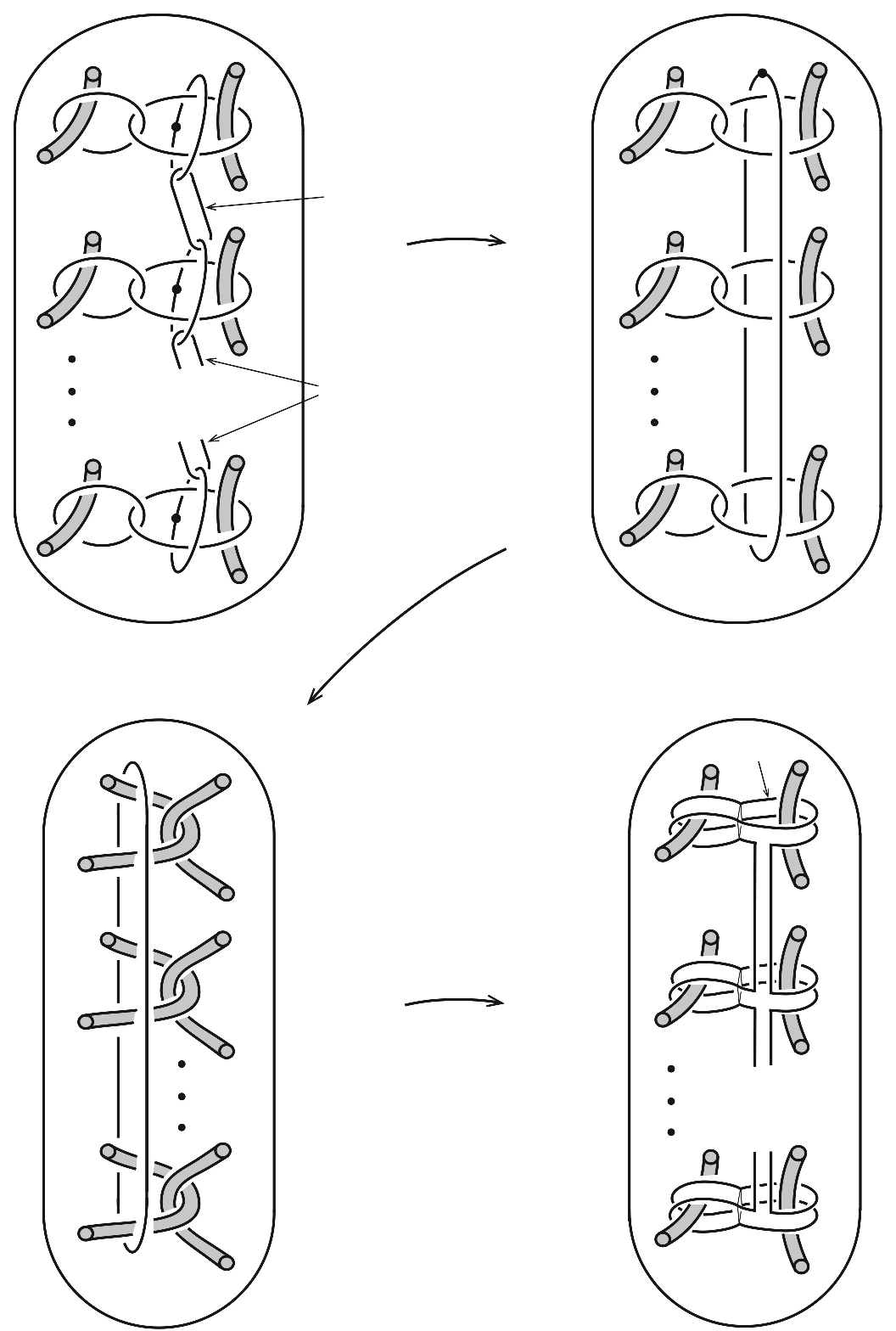} 
{\small 
\put(-200,250){$\overset{\partial}{\cong}\; $, $\cdot \mapsto 0$,}
\put(-200,238){and cancel}
\put(-200,225){$(x_i, y^+_i)$ pairs}
\put(-188,114){Isotopy}
\put(-188,102){of $H$}
\put(-200,410){Morse cancel}
\put(-200,397){$1$ and $2$ pairs}
\put(-295,283){$H$}
\put(-292,10){$H$}
\put(-64,283){$H$}
\put(-60,10){$H$}
\put(-60,225){$K_i$}
\put(-296,223){$0$}
\put(-218,439){new $h$}
\put(-218,363){new $h$}
}
{\scriptsize
\put(-302,481){$0$}
\put(-288,481){$0$}
\put(-309,452){$x_1$}
\put(-249,471){$y_1^+$}
\put(-287,433){$0$}
\put(-302,417){$0$}
\put(-288,417){$0$}
\put(-309,388){$x_2$}
\put(-249,406){$y_2^+$}
\put(-287,371){$0$}
\put(-302,328.5){$0$}
\put(-288,328.5){$0$}
\put(-309,299){$x_{g_i}$}
\put(-249,316){$y_{g_i}^+$}
\put(-287,342){$0$}
}
 \caption{}
\label{GeneralCalc3 figure}
\end{figure}

\begin{definition} \label{highergenus} ({\it Generalized double, no genus restriction})
A general GBL $K$ with Seifert surface $S$ as above is said to be a generalized double of the set 
 $$\{ [ (x^{}_1,y_1^-), \ldots, (x^{}_{g_1}, y_{g_1}^- )], \ldots, [(x^{}_{J-g_k}, y_{J-g_k}^-),\ldots,
 (x^{}_J, y_J^-)  ] \}.$$ (Note the set does not fully specify $K$. For each component $S_i$ of $S$ in addition to the $i^{\rm th}$ bracket $[\ldots ]$ of plumbed pairs, $g_i-1$ bands must be chosen to build $S_i$ as a band sum from these $g_i$ plumbed pairs.)
\end{definition}

Now dropping the genus restriction, the analogous  surgery problem source  $\overline M$ has some further $2$-handles. Beyond the $\pm$ plumbed pairs for each $(x_i, y^-_i)$ when more than one, say $g_i$, pairs lie on $S_i$ we must add an additional $g_i-1$ $2$-handles to add relations collapsing the $g_i$ free generators $e_1,\ldots, e_{g_i}$ dual to these plumbings to one. The relations can be taken to be a chain of loops in the simplest possible form (shown in figure \ref{GeneralCalc3 figure}) which read 
$e_1 e_2^{-1}, \ldots, e_{g_i-1}e_{g_i}^{-1}$. 

\begin{lemma} \label{general double lemma2} \sl
Now considering general doubling  (without genus restriction), we have defined a $4$-manifold $\overline M$ (the details of the attaching circles for the last $2$-handles $\{ h\}$ are as specified in figure \ref{GeneralCalc3 figure}). Again $\overline M$ is the source for an unobstructed surgery problem for building a slice complement (with ${\pi}_1$ freely generated by meridians). ${\mathcal S}_0(K)\cong \partial \overline{M}$.
\end{lemma}

{\it Proof of lemma \ref{general double lemma2}.}
For $g_i>1$ the required calculation  takes place in a genus $2g_i$ handlebody $H$ and is given in figure \ref{GeneralCalc3 figure}. \qed

{\it Note on terminology:} General doubles generalize Whitehead doubles but have little relation to Bing doubles. Bing doubles do not really weaken links but merely push nilpotent invariants further down the lower central series.

\section{Weak link homotopy} \label{weak homotopy section}

\subsection{}
This section defines an equivalence relation, weak link homotopy (WLH$_0$ - we will explain the zero), for $k$-component links $L$ in an oriented $3$-manifold $M$. It is slightly stronger than cobordism  by disjoint embedded surfaces in $M\times I$, but considerably weaker than Milnor's link homotopy. To give the definition we need some preliminary notions. 

We call the following space $\Delta$ a {\it $k$-branched disk}\hspace{.07cm}:
$${\Delta} \, =\, ( \amalg_{i=1}^k D^2_i) \, /\, ({\rho}, {\theta}, i) \equiv  ({\rho}', \, {\theta}', i') \; {\rm iff} \; {\rho}={\rho}'\leq 1/2, {\theta}={\theta}', \, 1\leq i, i'\leq k,$$

where $D^2_i$ is one of $k$ copies of the unit $2$-disk, figure \ref{weak homotopy figure}.
\begin{figure}[ht]
\includegraphics[height=2.5cm]{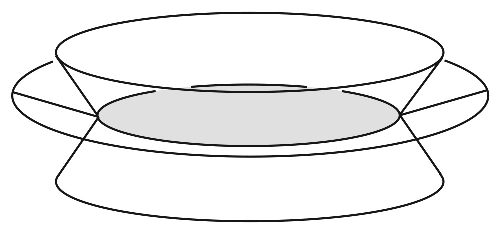}
{\small
\put(-190,30){${\Delta}\; =$}
\put(10,30){for $k=3$}
}
 \caption{A $3$-branched disk. The central disk $D^2_{1/2}$ is shaded}
\label{weak homotopy figure}
\end{figure}

We say ${\Delta}=$ ({\it central disk}) $\, \cup\,$({\it fringe}) \, = $D^2_{1/2}\, \cup\,  \overline{{\Delta}\smallsetminus D^2_{1/2}}.$ Consider a map $f\co {\Delta}\longrightarrow M\times I$ which is an ``embedding'' except for transversal double points with preimages in the common disk of radius $1/4$, $D^2_{1/4}\subset {\Delta}$. We put ``embedding'' in quotes since we design in a certain lack of smoothness along $\partial D^2_{1/2}$. Restricted to each $D^2_i$, $1\leq i\leq k$, $f$ is piecewise smooth with a crease along $\partial D^2_{1/2}$,  so that the $k$ inward normals to the fringe are always distinct and always outside the tangent space to $D^2_{1/2}$.

The normal bundle to $f(D^2_{1/2})$ has a unique trivialization up to homotopy (since $\Delta$ is contractible). Since the $k$ non-vanishing sections determined by the fringe are all in the same homotopy class, they  define a homotopy class of trivialization over $\partial D^2_{1/2}$. Comparing the two yields an integer which we call the relative Euler class ${\chi}$ of $f$. Maps $f$ as above with ${\chi}=0$ are called admissible.

\begin{definition} \label{WLH0 definition} (WLH$_0$) Two links $L_0=(l^0_1,\ldots, l^0_n)$ and $L_1=(l^1_1, \ldots, l^1_n)$ in $M$ are {\it WLH$_0$ equivalent} if there exists an embedding $$g\co ((\amalg_{i=1}^n S^1_i)\times[0,1])_- \hookrightarrow M\times [0,1],$$
with $g ((\amalg_{i=1}^n S^1_i)\times j) = L_j,\,  j=0,1$, identifying $M\cong M\times 0 \cong M\times 1$. The ${}_-$ subscript indicates the $n$ cylinders are multiply punctured. The punctures are grouped in an arbitrary manner into sets of cardinality $\{ k_1,\ldots, k_s\}$ and each of these $s$ groups is matched up to disjoint {\it admissible} maps $f_1,\ldots, f_s$, disjoint from $g$, except where they are matched along the $s$ fringes.  So the final result is a map which is an embedding except for double points in $s$ disks $\{ D^2_{1/4,\, l}\,, 1\leq l \leq s\}$ of a complex that is made from $(\amalg_{i=1}^n S^1_i)\times[0,1]$ by pressing in $s$ places $k_j$ sheets together, $1\leq j\leq s$. There is no restriction on the orientation of sheets or how many $k_j$ come together on each of the $s$ occasions. We call the union of maps $f\cup g = h$, and may regard it as $h\co ((\amalg_{i=1}^n S^1_i)\times[0,1]) \longrightarrow M\times [0,1]$.
\end{definition}

\begin{remark} There is an equivalent definition
 of WLH$_0$ where one further requires that $h$ is level preserving.  The reduction to level preserving form is an application of general position, following the argument of \cite{Giffen, Goldsmith} for ``concordance implies homotopy''. A weak null-homotopy may be presented as a sequence of elementary weak homotopies, illustrated in figure \ref{weak homotopy figure}.
 \begin{figure}[ht]
\includegraphics[height=5cm]{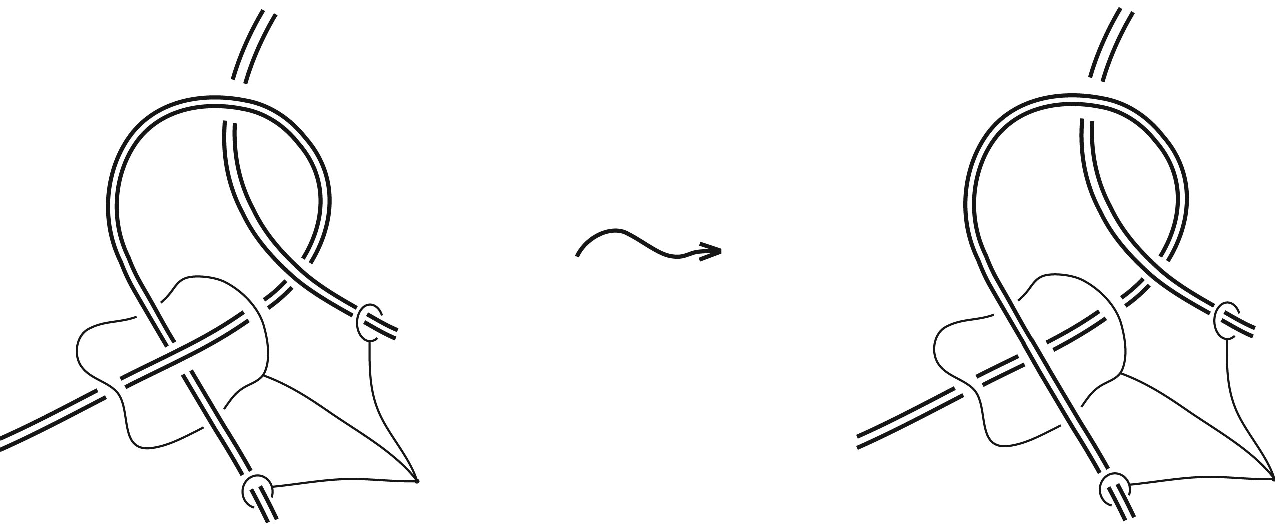}
{\small
\put(-288,7){$x$}
\put(-273,30){$\gamma$}
\put(-253,46){$y$}
}
 \caption{An elementary weak homotopy. The based curve ${\gamma}$ in the link complement, corresponding to the $2$-Engel relation $[x,x,y]=[x, x^y]$, becomes trivial after the indicated move. (Compare with figure \ref{link homotopy figure} illustrating link homotopy.)}
\label{weak homotopy figure}
\end{figure}
\end{remark}

\begin{remark} Instead of {\it zero}, any subset  of allowed integer relative Euler classes, closed under $x\mapsto -x$, could be used to define a version of weak link homotopy. In such generalizations linking is not conserved and we have not investigated these.
\end{remark}

\begin{remark} The WLH$_0$ equivalence relation was our pathway into studying the $2$-Engel relation, but in the end it is not essential to the logical structure of the paper. In section \ref{model section} the $2$-complex called the  ``Engel mess'' $X$ (see figure \ref{spaghetti}) is the image of a null WLH$_0$.
\end{remark}

\begin{remark}  In an earlier version of the paper we did not state the framing assumption ${\chi}=0$. We thank the referee for pointing out this omission.
\end{remark}

We have not fully calculated WLH$_0$ even for $M=S^3$, but consider this a nice open problem. A few comments on WLH$_0$ in $S^3$ come next.

\subsection{Comments on WLH$_\mathbf{0}$ in $\mathbf{S^3}$} \label{comments subsection}
WLH$_0$ is a coarser equivalence relation than Milnor's link homotopy \cite{M}, so it is of interest which $\bar\mu$-invariants (of link homotopy) survive to WLH$_0$, and which do not.

The zero-framed condition implies linking numbers $\bar\mu_{ij}, i\neq j$ survive, i.e. distinguish WLH$_0$ classes of links. Less trivially, $\bar\mu_{123}$ also survives to WLH$_0$.

\begin{proposition} \label{weak mu bar prop} \sl Milnor's invariant $\bar\mu_{123}$, with its usual indeterminacy, is an invariant of  WLH$_0$.
\end{proposition}

{\it Proof.} It is sufficient to consider three component links. The relations in ${\pi}_1(S^3\times [0,1]\smallsetminus {\rm image}( h) )$, associated to the transverse double points of $f_1,\ldots, f_s$, are of the $2$-Engel form $[x,x,y]$.

If we denote by $2E({\pi})$ the normal closure of $2$-Engel relations within a group ${\pi}$, and recall ${\pi}^2=[{\pi}, {\pi}], \ldots, {\pi}^j=[{\pi}^{j-1}, {\pi}]$ as our notation for the lower central series, then
\begin{equation}\label{inclusions}
{\pi}^3 \supseteq 2E({\pi}) \supseteq {\pi}^4.
\end{equation}

The first inclusion is immediate since the $2$-Engel relations are $3$-fold commutators. The second was proved in section \ref{Engel section}. (Corollary \ref{Engel corollary} contains more refined information.)

Turning now to $\bar\mu_{123}$, let $m={\rm g.c.d.}({\bar\mu}_{12},{\bar\mu}_{23},{\bar\mu}_{13})$. Up to sign, for a $3$-component link $L_0=(l^0_1,l^0_2, l^0_3)$
\begin{equation} \label{weak eq}
\bar\mu_{123}=r\in {\mathbb Z}/m{\mathbb Z}, \; {\rm where} \; [l^0_1]=rg\in C,
\end{equation}
where $g$ is a generator of the center $C$ of the class two nilpotent group $N_0$,
\begin{equation} \label{center}\begin{CD}
1 @>>> {\mathbb Z}_m @>>> N_0 @>>> {\mathbb Z}_m\times {\mathbb Z}_m @>>> 1\\
@. @| @. @.\\
@. C @. @. @.
\end{CD}
\end{equation}
corresponding to a generator of $H^2({\mathbb Z}_m\times {\mathbb Z}_m; {\mathbb Z}_m)$, 
\begin{equation} \label{def N0}
N_0=\left( {\pi}_1(S^3\smallsetminus (l^0_2\cup l^0_3))/ {\pi}_1(S^3\smallsetminus (l^0_2\cup l^0_3))^3 \right) \otimes {\mathbb Z}_m\footnote{This tensor operation may be defined on nilpotent groups via the composition series, or more simply by imposing the relation $g^m=1$, for all $g$.}.
\end{equation} 

Let $h\co (S^1_1\cup S^1_2\cup S^1_3)\times I \longrightarrow M\times I$ be a WLH$_0$ from $L_0$ to $L_1$ and $h_{23}$ be the restriction of $h$ to the last two components. (And let $h_1, h_2$ and $h_3$ be the restrictions to the first, second and third sheets respectively.) Note that $N_0\cong N_1$, where
\begin{equation} \label{def N0}
N_1=\left( {\pi}_1(S^3\smallsetminus (l^1_2\cup l^1_3))/ {\pi}_1(S^3\smallsetminus (l^1_2\cup l^1_3))^3 \right) \otimes {\mathbb Z}_m.
\end{equation} 

Of course, 
\begin{equation} \label{N}
N:= N_0/2E(N_0)\cong N_1/2E(N_1)\cong N_1\cong N_0.
\end{equation}

The only non-spherical homology of $M\times I\smallsetminus {\rm image}(h_{23})$, not coming from the top or bottom levels, is generated by the linking tori around the double points of image$(h_{23})$. These realize $2$-Engel relations which by (\ref{N}) are trivial in $N$, and by \cite[Lemma 13]{K2}
do not affect the class two nilpotent quotient. Thus
$N$ accepts a map ${\pi}_1(M\times I\smallsetminus {\rm image}(h_{23}))\longrightarrow N$.
Now $l^0_1$ and $l^1_1$ bound a multiply punctured cylinder $\gamma$ mapping into $M\times I\smallsetminus {\rm image}(h_{23})$. The top and bottom are $\partial_0 {\gamma}=l^0_1$ and ${\partial}_1{\gamma}=l^1_1$. The punctures $\{ {\partial}_i{\gamma} \}, i>1$, correspond to sheets of $h_1(S^1\times I)$ which are involved in weak homotopies (i.e. in the source of a map $g$ as above.)

Since $\partial_i {\gamma}$, $i>1$, has zero linking number with image$(h_2)$ and image$(h_3)$, each $[\partial_i{\gamma}]\in C$. In fact, examining the local geometric model, one sees that each $[\partial_i {\gamma}],$ $i>1$, is a conjugate of a loop of the form $[w,w^{\alpha}]$, a $2$-Engel relation. The required local model for $k=2$ sheets is shown in part (d) of figure \ref{singularities}. In this model $[\partial_i {\gamma}]$ becomes  a $(k+1)^{\rm st}$ ($3^{\rm rd}$) parallel round circle in the rightmost panel of the figure. It bounds a punctured torus $T_{-}$ in the complement of other components. The curves $w$ and $w^{\alpha}$ are visible on $T_{-}$ as a dual pair of small linking circles to the dotted component. This verifies that $[\partial_i {\gamma}]=[w,w^{\alpha}]$.
It follows that  $[{\partial}_i {\gamma}],$ $i>1,$ is trivial in $C$. Therefore $[{\partial}_0 {\gamma}]=[{\partial}_1 {\gamma}]\in C$, so ${\bar\mu}_{123}$ is invariant. \qed

The calculations in section \ref{Engel section} suggest that certain fourth order $\bar\mu$ invariants might also persist as mod $3$ invariants of WLH$_0$, but we have not checked this. On the other hand, we have:

\begin{proposition} \label{weak length 4} \sl Any link $L\subset S^3$, whose $\bar\mu$-invariants of length four and less vanish, is trivial under WLH$_0$.
\end{proposition}

{\it Proof.} By reordering components, it is sufficient to consider a first non-zero invariant, $\bar{\mu}_{12\ldots n}$ for $n\geq 5$. The vanishing of lower $\bar\mu$-invariants implies that $[l_1]\in {\pi}^{n-1}$, ${\pi}:=M{\pi}_1(S^3\smallsetminus (l_2\cup\ldots\cup l_n))$, and $\bar\mu_{12\ldots n}$ depends only on this element.
Since $n\geq 5$, line (\ref{inclusions}) implies ${\pi}^{n-1}\subseteq 2E({\pi})$.

Choose a sequence of elementary weak homotopies (figure \ref{weak homotopy figure}) which add the necessary relations which trivialize $[l_1]$. Then contract $l_1$ in the complement of this weak homotopy of $l_2\cup\ldots \cup l_n$. Finally reverse the sequence to obtain a link with $n-1$ components.  Proceed in this way until the number of components reaches four. At this point the link is actually homotopically trivial since the $\bar\mu$ invariants of length $\leq 4$ vanish. Composing these homotopies and weak homotopies provides the claimed weak null homotopy of $L$. 
\qed

\section{New universal surgery problems}   \label{model section}
For us surgery problems (``problems'') are in $4$ dimensions and, as studied in \cite{Wall}, non-singular over the integral group ring of the target.
It is irrelevant whether they arise in a closed or bounded context, the set up being a normal map
$f\co M\longrightarrow P$ from a $4$-manifold to a Poincar\'{e} space (or pair, but we drop the boundary from the notation). A ``solution'' is a normal cobordism to a simple homotopy equivalence.
In attempting to solve a given problem  one typically struggles  to embed $2$-complexes $X$ in the source $M$ of the surgery problem, where $X$ in some sense approximates $S^2\vee S^2$. There is considerable universality governing which $2$-complexes turn up, e.g. $S^2\vee S^2$-like  capped gropes\footnote{the right hand side of figure \ref{SurgeryGrope figure} is a height$=\! 1$ $S^2\vee S^2$-like capped grope.} \cite{FQ}, and the neighborhoods ${\mathcal N}(X)=M'$ of these standard examples are then themselves sources of bounded surgery problems to a $P'\simeq \vee S^1$. This set $\{ M'\longrightarrow P'\}$ or in abbreviation $\{ M'\}$ constitutes a countable, but easily parametrized list  of surgery problems, called ``universal''\footnote{also called ``model'' or ``complete''} which if solvable imply {\it all} unobstructed problems are solvable.
Solutions are transitive (see lemma \ref{transitivity}): if $M\subset M'$ (and captures the surgery kernel) then a solution to $M'$ solves $M$ as well.

\begin{figure}[ht]
\includegraphics[height=3cm]{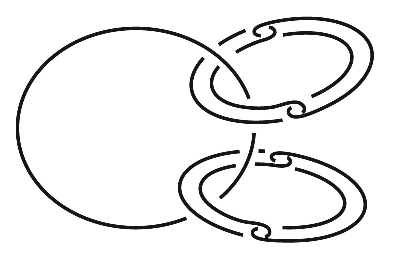}
 \caption{}
\label{RamifiedBor figure}
\end{figure}

\smallskip

The usual list of universal surgery problems is $\{\!${\sl Wh(Bing(Hopf))}$\!\}$. One starts with the Hopf link:$\HopfLink$\hspace{-2.5em}, then does any (non-zero) amount of iterated ramified Bing doubling. For example, Bing doubling the right component with ramification two (and no iteration) yields the link in figure \ref{RamifiedBor figure}.

Finally, $Wh$ means apply (without iteration) ramified $\pm$ Whitehead doubling to {\it each} component. See figure \ref{SurgeryGrope figure} for an example.
\begin{figure}[ht]
\includegraphics[height=3.7cm]{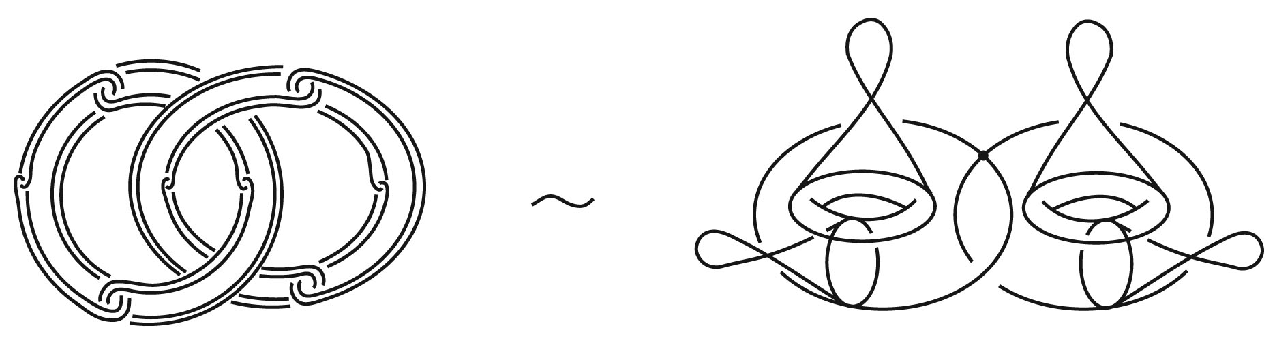} 
 \caption{The link on the left is obtained by Bing doubling both components of the Hopf link, and then Whitehead doubling each resulting component. The figure on the right illustrates the corresponding surgery kernel.}
\label{SurgeryGrope figure}
\end{figure}

The surgery problems $M$ are exactly the $M$ as in figure \ref{M figure} and lemma 
\ref{M lemma} where $L=\;${\sl Bing(Hopf)} (possibly ramifying the components with parallel copies), and by lemma \ref{M lemma} solving $M$ produces a ``free'' slice complement, i.e. one with free ${\pi}_1$, freely generated by meridians, for {\sl Wh(Bing(Hopf))}, i.e. constructs a manifold $\simeq\vee S^1$ (generators = meridians) with boundary ${\mathcal S}_0${\sl (Wh(Bing(Hopf)))}.

There is a natural partial order on capped gropes under inclusion: higher genus for surface stages and more double points for caps means ``inside''. Contraction \cite{F1, FQ} of capped gropes at grope tips, together with grope height raising \cite{F1, FQ}, show that both lower and higher stage capped gropes may lie inside a given one. (In terms of links, more stages corresponds to more ramified Bing doubling.) Any co-initial segment $\{ M'\}$ will be universal for surgery , so to solve surgery it is sufficient to handle cases with arbitrarily high (ramified) Bing doubling before ramified Whitehead doubling.

\begin{figure}[ht]
\includegraphics[height=2.4cm]{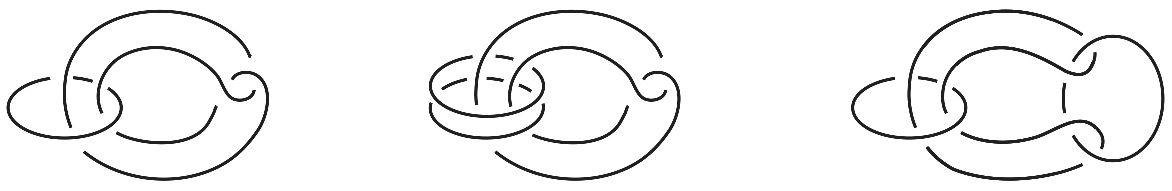}
{\small
\put(-420,15){$2$}
\put(-322,15){$1$}
\put(-270,43){$2$}
\put(-270,10){$2'$}
\put(-178,12){$1$}
\put(-110,10){$1$}
\put(-60,21){$2$}
\put(-7,8){$3$}
\put(-420,-15){$\bar\mu_{1122}=1$}
\put(-420,-28){Wh link is triv$^+$}
\put(-420,-41){M(Wh) solved.}
\put(-290,-15){$\bar\mu_{1122'}=1$}
\put(-290,-28){Wh$_{\text{\rm ram}}$, slightly ramified,}
\put(-290,-41){not triv$^+$.}
\put(-290,-55){M(Wh$_{\text{ram}}$): open question.}
\put(-135,-15){$\bar\mu_{123}=1$}
\put(-135,-28){Bor. rings$\in\{$Bing(Hopf)$\}$.}
\put(-135,-41){M(Bor) a presumably typical}
\put(-135,-55){universal problem.}
}
\caption{}
\label{three cases}
\end{figure}

There seems to be a hierarchy of links (with all linking numbers vanishing). The strongest - hardest to slice their Whitehead doubles - are the ramified iterated Bing doubles of Hopf, Bing(Hopf). These are universally hard. At the other extreme are the (homotopically trivial)$^+$ - links $L$ (triv$^+$). By definition these are the $k$-component links $L$ which if turned into $L^+$, a $(k+1)$-component link obtained by adding a parallel copy to any single component, $L^+$ is homotopically trivial in the sense of Milnor (see section \ref{Milnor group section}). It is known \cite{FT2} that for such $L$ the associated problem $M$ can be solved. To be very explicit consider three cases, figure \ref{three cases}.

The basic Whitehead link is ``held together'' by a nontrivial $\bar\mu_{1122}$ but it is triv$^+$ since both indices must be repeated to obtain a non-zero value.
Wh$_{\text{ram}}$, a temporary notation for the $3$-component link obtained by ramifying one component of Wh, shown as the middle figure in \ref{three cases}, still has no nontrivial $\bar\mu$ invariants with non-repeating indices, but it is no longer triv$^+$. Despite 20 years of work its status is still open: we do not know if the surgery problem M(Wh$_{\text{ram}}$)  is solvable or even if Wh(Wh$_{\text{ram}}$) is slice (with any fundamental group). 
Finally, the Borromean rings seems typical (though at the top of the partial order) of links whose associated surgery problems $M$ form the universal family. Bor (and other members of $\{$Bing(Hopf)$\}$) are held together by a non-repeating $\bar\mu$ invariant.

In this section we exhibit new universal links which are general (see section \ref{s-cobordism section}) doubles of links $L$ which are more like Wh$_{\text{ram}}$ than Bor. For convenience of the reader we restate the theorem from the introduction:

{\bf Theorem \ref{new models theorem}}. {\sl
There is a family of links $\{ K\}$  for which the problem of constructing free slices constitutes a universal problem, where each $K\in \{ K\}$ is of the form: $${\text{D(Ram(h-triv)),}}$$ a generalized (genus one) double of a ramified homotopically trivial link.}

\begin{remark}  \label{universal theorem remark}
Any 
link $L$ with a non-trivial $\bar\mu$ may be ramified to produce $L_{\rm ram}$ with a nontrivial non-repeating $\bar\mu$, 
so theorem \ref{new models theorem} does not yield a universal set  of links based on general doubling of a link with all 
non-repeating $\bar\mu$ invariants vanishing. However 
it was a surprise to see any collection of universal problems tied so closely to homotopically trivial links. The new universal problems certainly focus attention on the impact of ramification. This is the second surprise and is perhaps the ``other side of the coin'' to section \ref{homotopy solution} where ramification is exploited to construct an unexpected homotopy solution to the A-B slice problem.
\end{remark}
\begin{remark}
Part of what we have learned is that what 20 years ago was regarded as a minor technical distinction - homotopically trivial versus (homotopically trivial)$^+$ - may instead lie at the heart of the matter. (See \cite{CP} for recent results related to this problem.) At present, Wh(Wh$_{\text{ram}}$) is certainly the most interesting link slice problem.
\end{remark}
The next lemma expresses the transitivity of surgery solutions; we will use it to pass from the original to the new universal problems.
\begin{lemma} \label{transitivity} \sl (Transitivity of surgery)
Let $f\co M\longrightarrow P$ be an unobstructed problem with surgery below the middle dimension completed so that $f_\sharp$ is an isomorpshim on ${\pi}_1$ and so that
the kernel module  $$K(f)=\text{\rm ker}(H_2(M;{\mathbb Z}[{\pi}_1P])\longrightarrow H_2(P;{\mathbb Z}[{\pi}_1P]))$$
is a free module with hyperbolic intersection form $\lambda$ and standard $\mu$. Suppose $K(f)$ is represented by an embedded (possibly disconnected) $4$-manifold
$(W, \partial)\subset \text{\sl interior}(M)$. Here ``represented'' means that 
$$H_2(W;{\mathbb Z}[{\pi}_1P])\overset{\text{inc}_*}{\longrightarrow} H_2(M;{\mathbb Z}[{\pi}_1P])$$
is an isomorphism. We assume that ${\pi}_1(W)\cong F_k$ is free and $W$ itself is the source of an unobstructed problem, $g\co W\longrightarrow Q$ where $Q=(\natural _k S^1\times D^3, \sharp_k S^1\times S^2)$. Suppose $g'\co V\longrightarrow Q$ is a solution to $g$, then there is a solution to $f$ of the form $f'\co M':= (M\smallsetminus W)\cup V\longrightarrow P$.
\end{lemma}
{\it Note.} We have used notations in the statement and proof appropriate to the case of orientable $M$. There is no difficulty extending to the non-orientable case.

{\it Proof.} Consider the two braided Mayer-Vietoris sequences below, where for ease of reading we have suppressed the functor 
$H_*({}_{-};{\mathbb Z}[{\pi}_1P])$ applied to each space and simply written a subscript for the value of $*$.
\begin{equation} \label{braided}
\xymatrix@C=1pc@R=.5pc
{
& &  (M\smallsetminus W)_2\oplus W_2    \ar[r] & M_2  \ar[dr]!UL^(.6){\partial}  &     &    (M\smallsetminus W)_1\oplus W_1    \ar[r]  & M_1          \ar[dr]!UL^(.6){\partial} & \\
\ldots  \ar[r] & \partial W_2  \ar[ur]!DL \ar[dr]!UL & & & \partial W_1   \ar[ur]!DL \ar[dr]!UL & & & \partial W_0 \\
& & (M\smallsetminus W)_2\oplus V_2    \ar[r] & M'_2   \ar[ur]!DL^(.8){\partial}      &  &    (M\smallsetminus W)_1\oplus V_1    \ar[r] & M'_1 \ar[ur]!DL^(.8){\partial   }  & \\
}
\end{equation}
A theorem of Wall \cite{Wall0} constructs $1$-manifold $1$-skeleton for Poincar\'{e} spaces $P$ so we may embed $Q\subset P$ and thereby obtain maps from line (\ref{braided}) to line (\ref{sequence2}), the corresponding Mayer-Vietoris sequence of Poincar\'{e} spaces.
\begin{equation} \label{sequence2}
\begin{split}
\ldots\longrightarrow \sharp_k(S^1\times S^2)_2\longrightarrow (P\smallsetminus Q)_2\oplus Q_2\longrightarrow P_2 \overset{\partial}{\longrightarrow} 
\sharp_k(S^1\times S^2)_1 \\
\longrightarrow (P\smallsetminus Q)_1\oplus Q_1\longrightarrow P_1\longrightarrow  \sharp_k(S^1\times S^2)_0
\end{split}
\end{equation}
Thus we now reinterpret line (\ref{braided}) not as homology groups but the kernel groups, $K_*(\; \;\; ; {\mathbb Z}[{\pi}_1 P])$. 

Since $V\simeq \vee S^1 \simeq Q$, regardless of the cover ${\pi}_1 P$ induces on $V$ and $Q$ the homotopy equivalence lifts to a homotopy equivalence $\widetilde V \simeq \widetilde Q$ implying $K_*(V;{\mathbb Z}[{\pi}_1 P])\cong 0$. In fact, the {\it only} non-trivial kernel group in line (\ref{braided})  is $K_2(W;{\mathbb Z}[{\pi}_1 P])$. Thus $K_*(M';{\mathbb Z}[{\pi}_1 P])\cong 0$ for all  $ *$. Finally, a calculation using van Kampen's theorem implies that $f_{\sharp}\co {\pi}_1 M'\longrightarrow {\pi}_1 P$ is an isomorphism, so $f'\co M'\longrightarrow P$ is a homotopy equivalence. \qed

Armed with this lemma, and knowledge of the strength of the $2$-Engel relation, our task is now to draw a detailed schematic ``spaghetti picture" of a new class of problems $W$ inside $M$ and carrying its kernel. $M$ is a typical member of a (coinitial) segment of the original universal problems based on {\sl Wh(Bing(Hopf))} and $W$ will be a new problem based on {\sl D(Ram(h-triv))}. Our spaghetti picture is intended to capture all the important features. Once drawn we will convert it into a schematic Kirby link diagram about which  we have enough precise knowledge to prove theorem \ref{new models theorem}.

We remember there is ramification but do not draw it in schematics. We start with any element of {\sl Bing(Hopf)} in which one component, e.g. $l_1$ in figure \ref{5 components}, is trivial in the Milnor group of the remaining components modulo the $2$-Engel relation, as in corollary \ref{Engel corollary}.
The simplest of these has $5$ components, and is drawn in figure \ref{5 components}.
We use this $L$ for illustration but formally choose any $L\in \{${\sl Bing(Hopf)}$\}$ with all $\bar\mu$-invariants of length $\leq 4$ vanishing. 

\begin{figure}[ht]
\includegraphics[height=8cm]{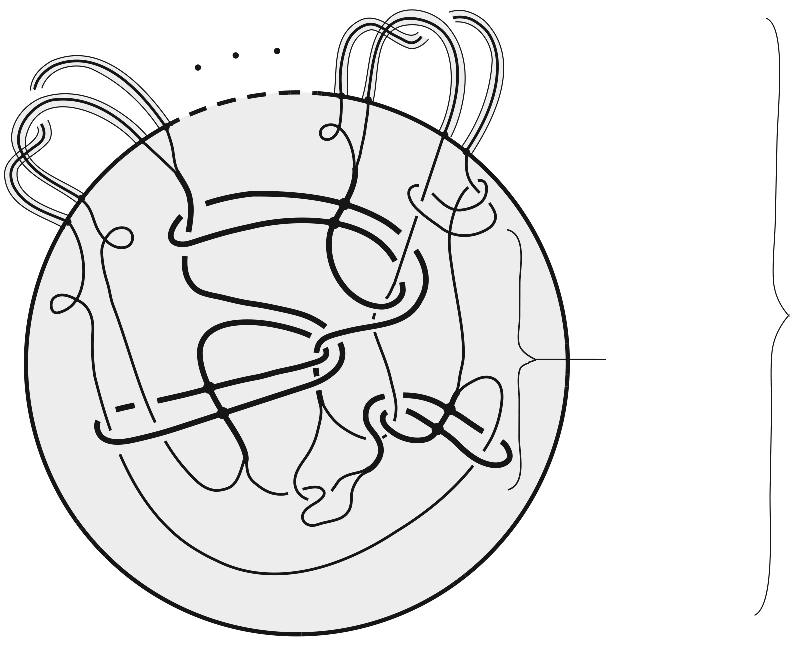}
{\scriptsize
\put(-184,9){$M(L)$}
\put(-66,97){``Engel mess'',}
\put(-50,85){$X_-$}
\put(3,112){$X$}
\put(-232,240){$5$ plumbed}
\put(-232,229){handle pairs}
\put(-232,218){for $L$ in fig. \ref{5 components}}
}
 \caption{}
\label{spaghetti}
\end{figure}

This property that all $\bar\mu$-invariants of length $\leq 4$ vanish  for a link $L$ is inherited by all links $PL$, by definition the class of links where the initial $k$ components of $L$ are arbitrarily ramified with parallel copies.
$P_2L$ means two parallels are taken.  Proposition \ref{weak length 4} implies that not only is any element of $P_2 L$ trivial under WLH$_0$ but trivial via geometric moves, $0$-framed paired link homotopies, which are the geometric counterpart of the algebraic factors enumerated in corollary \ref{Engel corollary}. The totality of these moves is weak link homotopy (of a special form) and its general position image in $D^4$ will be denoted by $X_-$ and called the ``Engel mess''.  When $X_-$ is joined to the cores of $k$ plumbed  $2$-handles the result, $X$, is a singular image of $\sqcup_{i=1}^k (S^2\vee S^2)$ in $M(L)$, as pictured in figure \ref{spaghetti} (compare figure \ref{M figure}).

The weak null-homotopy in figure \ref{spaghetti} has (arbitrarily) 5 ordinary self-intersections pictured and 3 sections where two strands run together for a bit (they are called a ``packet'' in section \ref{weak homotopy section}, and drawn heavily in figure \ref{spaghetti}) and while fused have 6 (again arbitrary) self-intersections. Note two is not arbitrary but the number of strands in a packet required by corollary \ref{Engel corollary}. 
\begin{remark} The commutator identity  
(\ref{inverse eq}) implies that the commutators $[h_1,\ldots, h_k]$ in the statement of corollary \ref{Engel corollary} may be assumed to involve only the given normal generators $\{g_i\}$, and not their inverses. It follows that the two-strand ``packets''
may be assumed to be of a single type, where (for some choice of orientations) the orientations of the two strands match. 
\end{remark}
\begin{figure}[ht]
\hspace*{-6.8cm}\includegraphics[height=10cm]{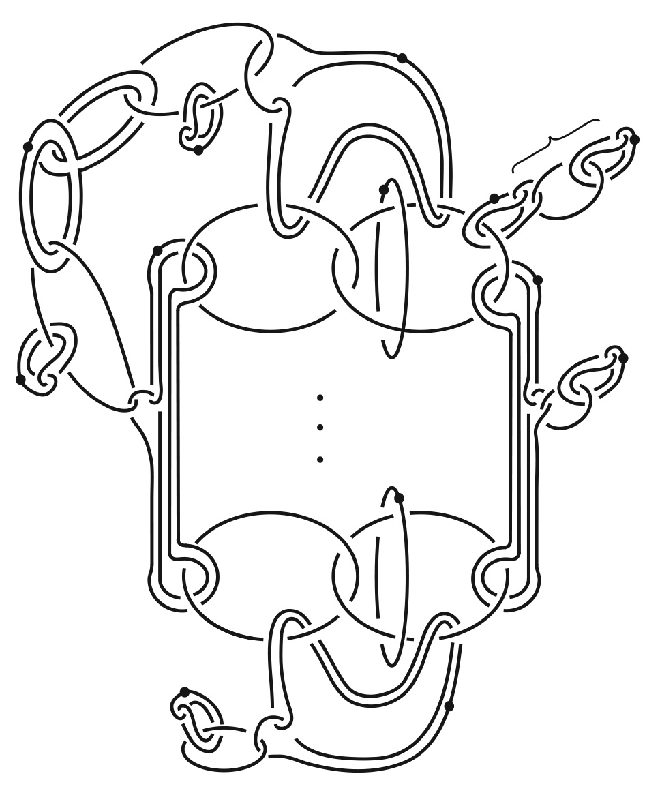}
{\small
\put(-214,253){a}
\put(-183,270){c}
\put(-166,264){$0$}
\put(-175,221.5){b}
\put(-232,176){c}
\put(-202,186){$0$}
\put(-230,133){b}
\put(-168,140){d}
\put(-116,128){k}
\put(-63,140){d}
\put(-148,-1){$0$}
\put(-183,22){b}
\put(-150,154){f}
\put(-81,154){f}
\put(-150,170){$0$}
\put(-81,170){$0$}
\put(-150,88){$0$}
\put(-150,101){f}
\put(-81,102){f}
\put(-99,114){g}
\put(-99,147){g}
\put(-81,88){$0$}
\put(-43,237){e}
\put(-25,129){c}
\put(-13,140){b}
\put(-38,121){$0$}
\put(-34,197.5){$0$}
\put(-171,2){c}
\put(-74,245){d}
\put(-74,17){d}
\put(21,200){a: an intersection point of}
\put(35,187){caps (only one pictured)}
\put(21,174){b: a self-intersection point of caps}
\put(21,161){c: a cap}
\put(21,148){d: a double pt. of paired sheets}
\put(21,135){e: an ordinary capped double point}
\put(35,122){of $X$ (only one pictured)}
\put(21,109){f: the $2k$ essential $2$-spheres}
\put(21,96){g: the $k$ double points}
\put(35,83){(plumbings) outside $D^4$}
}
\caption{Note: type d $1$-handles occur (numerous) times between all $\binom{2k}{2}$ pairs of of types f curves.}
\label{EngelKirby figure}
\end{figure}
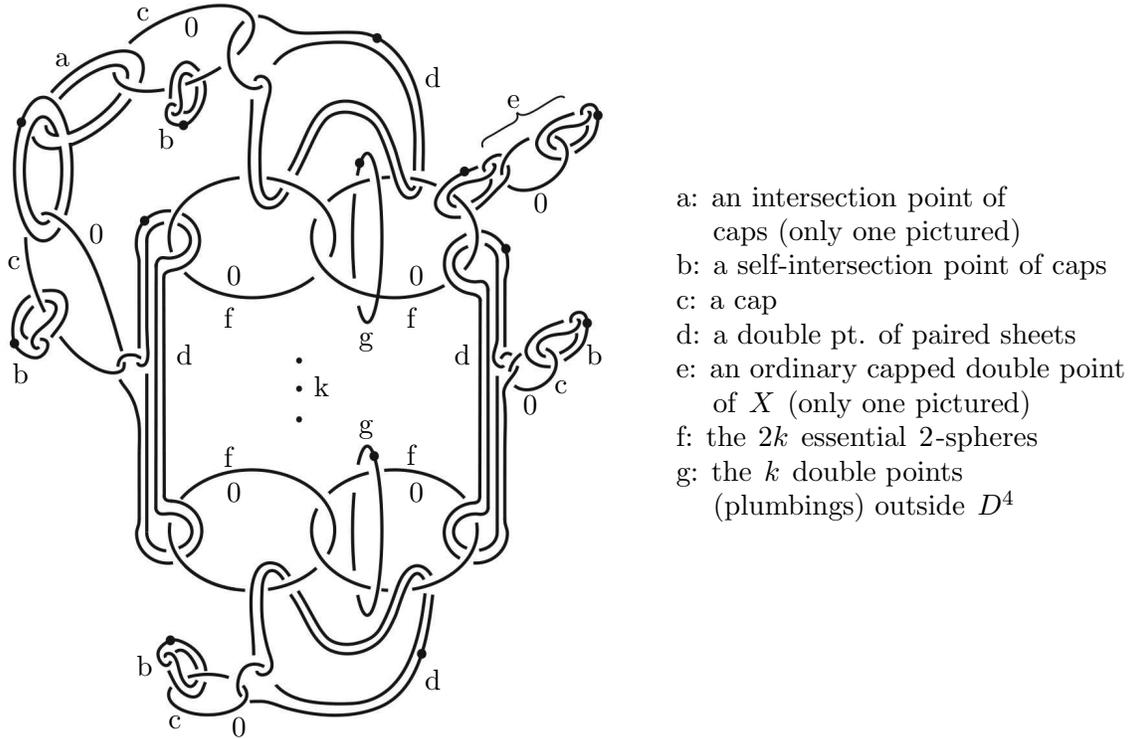
We call the $2$-complex $X$ pictured in figure \ref{spaghetti} the closed ``Engel mess''. We need to make it still messier so that it contains dual spheres. Actually corollary \ref{Engel corollary} allows us to build the weak null homotopy so that each of the $2k$ longitudes ${\gamma}_1,\ldots, {\gamma}_{2k}$ to $P_2 L$ lies in the kernel $K$:
\begin{equation} \label{triangle eq}
1\longrightarrow K\longrightarrow 
{\pi}_1(S^3\smallsetminus P_2L) \longrightarrow {\pi}_1(D^4\smallsetminus X) \longrightarrow 1
\end{equation}
This is because each ${\gamma}_i\in ({\pi}_1(S^3\smallsetminus P_2L))^4$, which follows from the formula \cite{M} for the behavior of $\bar\mu$-invariants under ramification.

Thus each of the $2k$ basic spherical classes $S_i$ in $X$ contains a geometrically dual sphere $S_i^{\perp}$ in $M$ (meeting $X$ only in one point of the basis - sphere.) The null-homotopy of ${\gamma}_i$ in $(D^4\smallsetminus X)$ glued to the $2$-handle core parallel to ${\gamma}_i$ is one the geometric dual spheres. Of course these $2k$ dual spheres $\{ S^{\perp} \}$ intersect and self-intersect.

It is now a standard construction \cite{F0, FQ} that the $\{ S^{\perp} \}$ can be exploited  to add an additional layer of zero-framed caps (intersecting and self-intersecting) to ``kill'' all double points of $X$ in $D^4$ (the $k$ double points in the plumbed handles cannot be capped.) Let $X^+=X\cup$(caps) and $W=$Neighborhood$(X^+)$; figure \ref{EngelKirby figure} exhibits the essential features of a handle diagram for any such $W$.

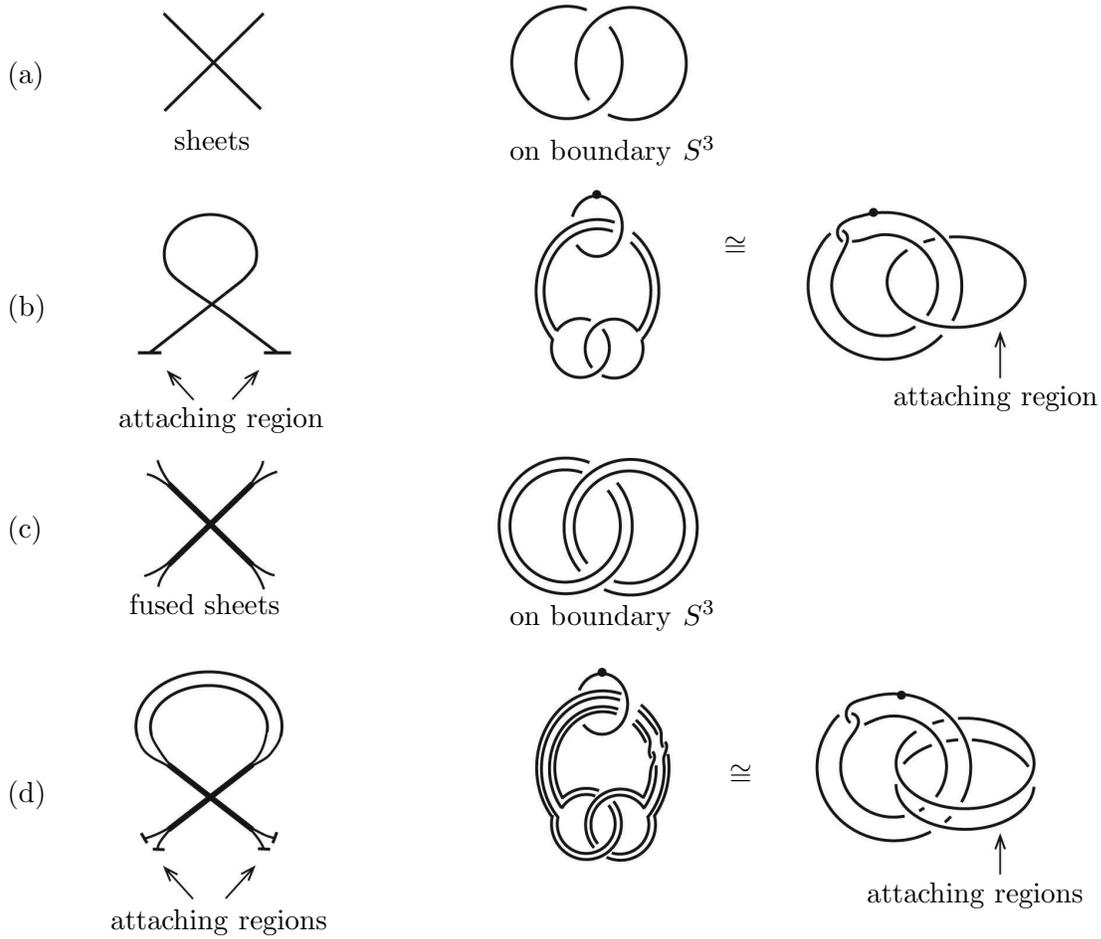
\begin{figure}[ht]
\includegraphics[height=12cm]{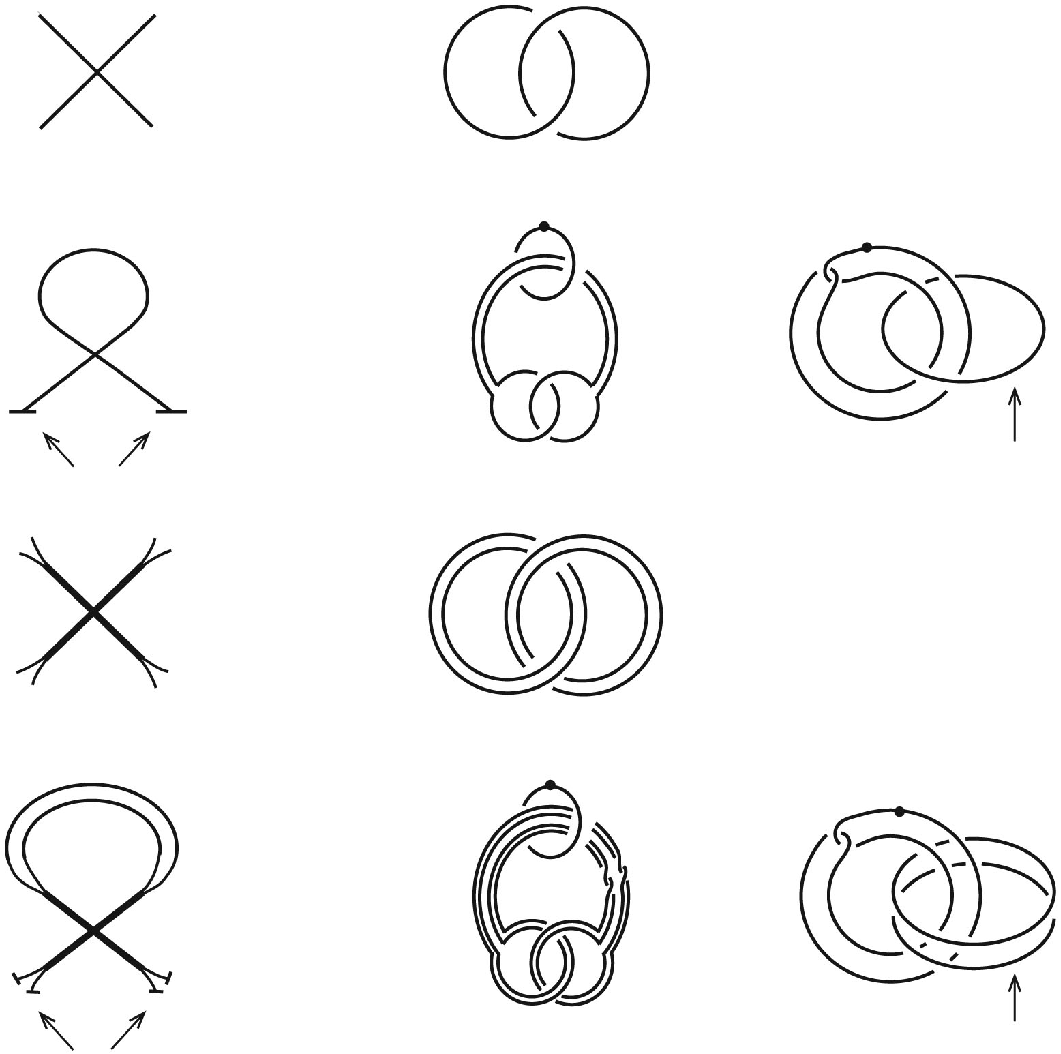}
{\small
\put(-390, 310){(a)}
\put(-390, 222){(b)}
\put(-390,138){(c)}
\put(-390,38){(d)}
\put(-327,285){sheets}
\put(-348,180){attaching region}
\put(-200,281){on boundary $S^3$}
\put(-55,188){attaching region}
\put(-344,109){fused sheets}
\put(-200,105){on boundary $S^3$}
\put(-351,-10){attaching regions}
\put(-65,0){attaching regions}
\put(-120,245){$\cong$}
\put(-118,45){$\cong$}
}
\caption{Four types of singularities and their associated Kirby diagrams}
\label{singularities}
\end{figure}

A slightly novel identity needed to draw figure \ref{EngelKirby figure} is the effect of a double point on a paired sheet, calculated in figure \ref{singularities}.

\begin{figure}[ht]
\includegraphics[height=8.7cm]{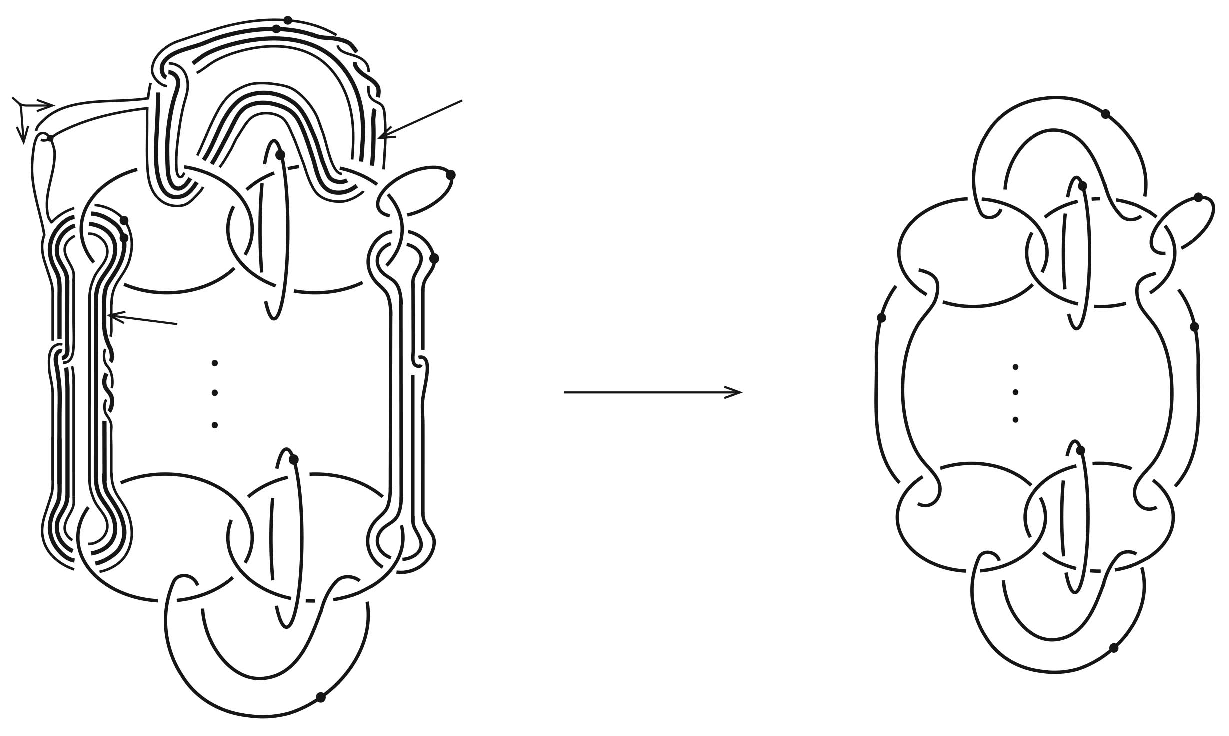}
{\small
\put(-38,209){$D^2$}
\put(-4,175){$D^2$}
\put(-136,95){$D^2$}
\put(-9,95){$D^2$}
\put(-32,25){$D^2$}
\put(-90,94){$0$}
\put(-90,134){$0$}
\put(-42,94){$0$}
\put(-42,134){$0$}
\put(-280,15){(a)}
\put(-120,15){(b)}
\put(-420,216){$D$}
\put(-257,213){$D$}
\put(-354,135){$D$}
\put(-372,18){$D^2$}
\put(-269,95){$D$}
\put(-260,184){$D^2$}
\put(-240,125){or more compactly}
\put(-320,135){$p_1$}
\put(-320,100){$p_k$}
}
\caption{}
\label{EngelKirbyCompact figure}
\end{figure}

Next we should Morse cancel $1$-handles and $2$-handles whenever possible in figure \ref{EngelKirby figure}. The result is a bit complicated to draw but we can use a short hand writing $D$ whenever a component is ``generalized double'' (see section \ref{general double section}).  In some cases (look at e in figure \ref{EngelKirby figure}) this is precisely a Whitehead double - the curve drawn in figure 6.7 should be doubled with itself (i.e. Whitehead doubled) to produce the original curve (e). In other case, a, two separate curves are generalized doubled, this indicated by the forked arrow in figure \ref{EngelKirbyCompact figure} which points to two curves representing bands to be plumbed. The boundary of this plumbing is the dotted curve (a) in figure \ref{EngelKirby figure}. 
The b-curves in figure \ref{EngelKirby figure} account for the dotted curves labeled by $D$ with a straight arrow. 
$D^2$ means ``double twice''. With this notation figure \ref{EngelKirby figure} becomes figure \ref{EngelKirbyCompact figure} with ramification permitted before each doubling step.

\begin{figure}[ht]
\includegraphics[height=2.7cm]{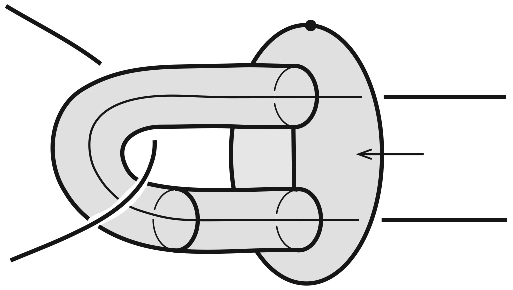}
{\small
\put(-21,32){Seifert surface}
\put(-43,-2){$p$}
\put(-89,35){$x$}
\put(-95,2){$y$}
}
\caption{}
\label{Seifert figure}
\end{figure}

The only $1$-handle components $p_1,\ldots, p_k$ ($p$ for plumbing) not (twice) doubled  in figure \ref{EngelKirbyCompact figure} (b) can be thought of as doubles of a pair $x, y^-$ using the obvious Seifert surfaces these components bound disjoint from the rest of the link diagram, figure \ref{Seifert figure}.

\begin{figure}[ht]
\includegraphics[height=6.5cm]{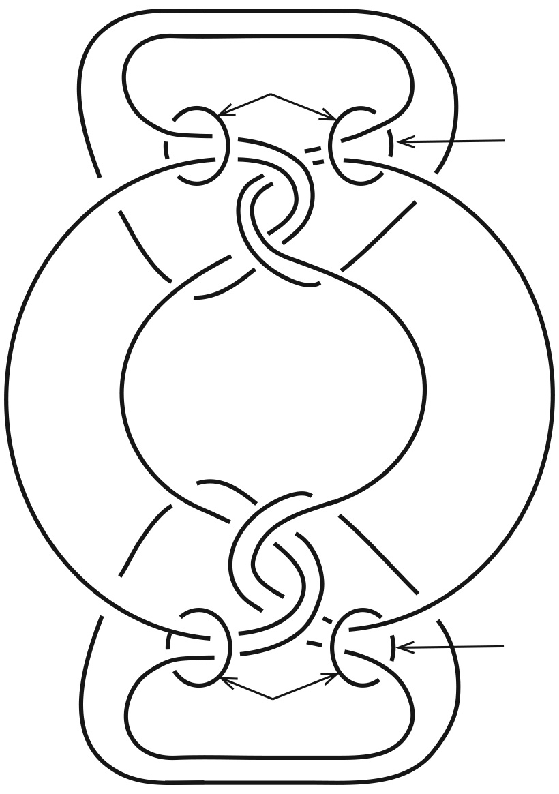}
{\scriptsize
\put(-11,32){$y_k^-$}
\put(-71.5,13){$D$}
\put(-91,46){$x_k$}
\put(-126,9){$D^2$}
\put(-144,85){$D^2$}
\put(1,85){$D^2$}
\put(-126,167){$D^2$}
\put(-71.5,166){$D$}
\put(-11,151){$y_1^-$}
\put(-91,164){$x_1$}
}
\caption{}
\label{FinalKirby figure}
\end{figure}

Canceling all hyperbolic pairs and replacing $p_i$ as Double ($x_i, y_i^-$) we get figure \ref{FinalKirby figure}.
 
Figure \ref{FinalKirby figure} actually displays a link of the type claimed in theorem \ref{new models theorem}. Leave $x_i$ and $y_i^-$ as they are and double each component labeled by a $D^2$ {\it once}. The result is a ramification of a homotopically trivial link - but not better (in our heuristic hierarchy of Whitehead doubles introduced in this section; ``better'' meaning {\it closer} to being sliced by our current techniques.) Now {\it all} these components are singly doubled (which includes the possibility of ramifying before taking a generalized double) to obtain the link in figure \ref{FinalKirby figure}, 
showing that $\partial W\cong{\mathcal S}_0$(D(Ram(h-triv))), as desired. \qed

\medskip

{\bf Acknowledgments}.  
We would like to thank Ian Agol for sharing his expertise on Engel groups. We would like to thank the referee for helpful suggestions which
improved the paper.

VK was partially supported by the NSF grant DMS-1309178 and by the Simons Foundation grant 304272.


\bigskip

\appendix

\section{$n$-Engel relations} \label{Appendix}
Since the $2$-Engel relation is rather interesting in $4$D topology, it is natural to wonder what use higher $n$-Engel relations might have. First we summarize a bit of what is known algebraically and raise some new questions natural when one works (as so often in $4$D topology) in the free Milnor group, rather than in the free group.

The fact 
about $2$-Engel groups that is crucially used in 
applications in this paper 
is that the 
free group $F_k$ modulo the $2$-Engel relation is nilpotent of a fixed class, independent 
of the number of generators $k$. It is known that $3$-Engel groups \cite{Hei} and also $4$-Engel 
groups  \cite{HVL} are locally nilpotent (i.e. every finitely generated subgroup is nilpotent). The question for $n$-Engel groups, $n>4$ presently appears to be open.
We do not know how the (local) nilpotency class of the free group $F_k$ mod the $3$- or $4$-Engel relation depends on $k$. 
Recall (section \ref{Milnor group section}) that the free Milnor group on $k$ generators, $MF_k$, is nilpotent of class $k$.
\begin{question} \sl
Fix $n\geq 3$. Is the nilpotency class of $MF_k$ modulo the $n$-Engel relation less than $k$? More specifically, is it independent of $k$?
\end{question}

An affirmative answer to this question could lead to an improvement of the results of geometric applications of the $2$-Engel relation in this paper.  
Another possible way to refine the algebraic structure is to see if there is a way to restrict (for example to a certain term of the lower central series) the group elements $y$ that come up in applications of the $2$-Engel relation $[y,x,x]$.

\begin{figure}[ht]
\includegraphics[height=3.9cm]{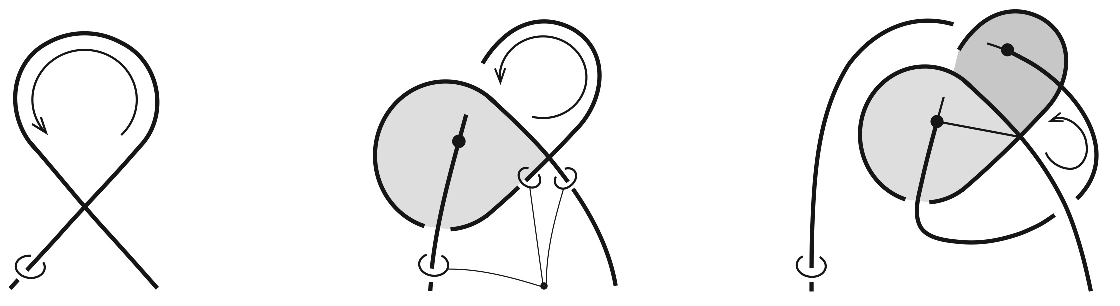}
{\scriptsize
\put(-378,73){$y$}
\put(-410,12){$x$}
\put(-212,79){$y$}
\put(-284,12){$m_1\!=\!x$}
\put(-24,56){$y$}
\put(-126,11){$x$}
\put(-232,45){$m_2$}
\put(-196,45){$m_3$}
}
\put(-420,-22){(a) kinky handle}
\put(-300,-22){(b) ``secondary kinky handle''}
\put(-125,-22){(c) ``tertiary kinky handle''}
\caption[caption]{\\\hspace{\textwidth} (a) 
The double point implies $[x,x^y]=1$ or equivalently $[y,x,x]=1$. \\\hspace{\textwidth} 
(b) 
 The second double point corresponds to $[y,x,x,x,x]=1$. \\\hspace{\textwidth} 
(c) 
The third double point corresponds to
$[y,x,x,x,x,x,x]=1$, etc.}
\label{CassonHandles figure}
\end{figure}

Returning to the first sentence of the introduction, we should follow Casson's philosophy  and determine  the local singular disk structure which enforces $[y,x,\ldots,x]$, the $n$-Engel relation. This turns out to be rather easy beginning with the ``kinky handle'' and then elaborating. The correspondence is laid out in figure \ref{CassonHandles figure}.
(The $n$-ary kinky handle gives raise to the $2n$-Engel relation.)

We sketch the calculation in (b). Denoting the meridian $m_1$ to the secondary kinky handle $H$ in $D^4$ by $x$, the meridian $m_2$ equals $x^y$, and $m_3=(x^y)^x$.
Then a meridian $m_{\text{\rm w}}$ to the (shaded) Whitney disk reads off $$m_{\text{\rm w}}=[m_2,m_3]=[x^y,(x^y)^x]=[x^y,[x,x^y]]=[x^y,[x,[y,x]]].$$
Finally, the $2$-cell of the Clifford torus of the second double point gives raise to a relation
$[m_{\text{\rm w}}, m_1]=[[x^y,[x,[y,x]]], x]$ in ${\pi}_1(D^4\smallsetminus H)$.
Using the commutator identities (\ref{product id}), modulo higher order commutators this gives the $4$-Engel relation
$[y,x,x,x,x]$.

\begin{figure}[ht]
\includegraphics[height=2.45cm]{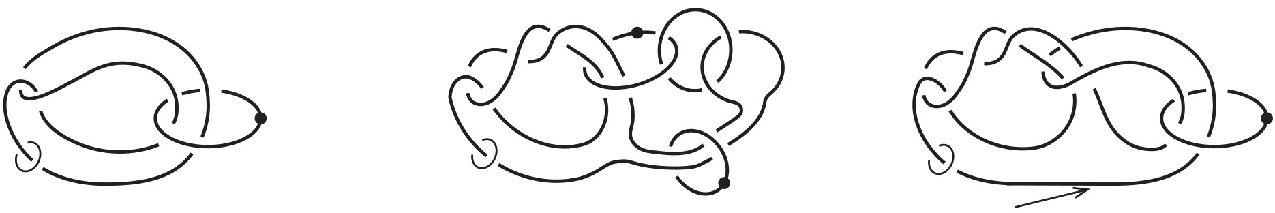}
{\small
\put(-333,20){$y$}
\put(-423,12){$x$}
\put(-177,8){$y$}
\put(-273,12){$x$}
\put(-183,65){$0$}
\put(-2,23){$y$}
\put(-122,10){$x$}
}
\put(-408,-10){kinky handle}
\put(-146,30){$\cong$}
\put(-275,-12){``secondary kinky handle''}
\put(-112,-8){attaching region}
\caption{}
\label{KinkyHandles figure}
\end{figure}

A more rigorous Kirby diagram description  (or if you like ``definition'') of higher order kinky handles is given in figure \ref{KinkyHandles figure}.
As usual we suppress discussion of the $\pm$ sign at clasps, but we do pay careful attention to framings so that the two components of the final, canceled down, link diagram are each individually unknotted (as with the Whitehead link). 
This is a useful feature when dualizing link diagrams. 

\begin{figure}[h]
\includegraphics[height=2.5cm]{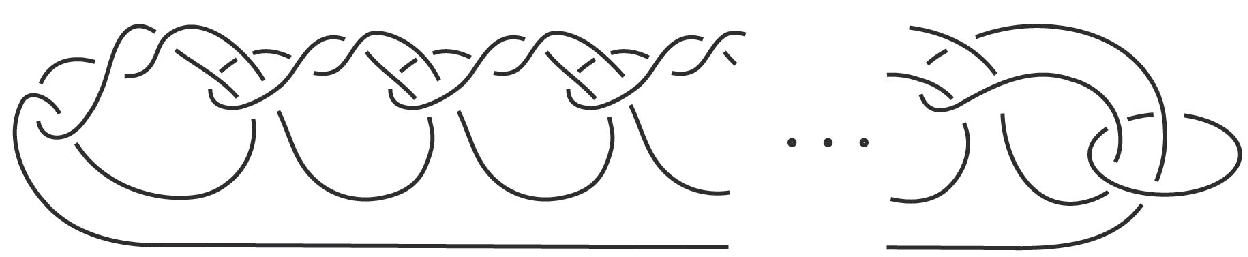}
\caption{}
\label{Engel link figure}
\end{figure}

Continuing in this way the model $n$-ary kinky handle is shown in figure \ref{Engel link figure}.
\smallskip

\begin{remark}
It looks very likely that the least area unknotting disk for the longer component as in \cite{HT} is exponential in $n$. Or put another way, if that component is made round the link diagram will necessarily have $\geq \text{\rm const}^n$ crossings for some const$\geq 1$. The method of \cite{HT} looks relevant to this case as well, but we did not succeed in adapting the argument. The new feature here is that the bridge number  of the diagram is not constant (as in \cite{HT}) but linear in $n$.
\end{remark}

These links are similar to Milnor's family \cite[Figure 1]{M2}, which arise from thickenings of figure \ref{CassonHandles figure} with less favorable choices of framings.

\begin{figure}[h]
\includegraphics[width=13cm]{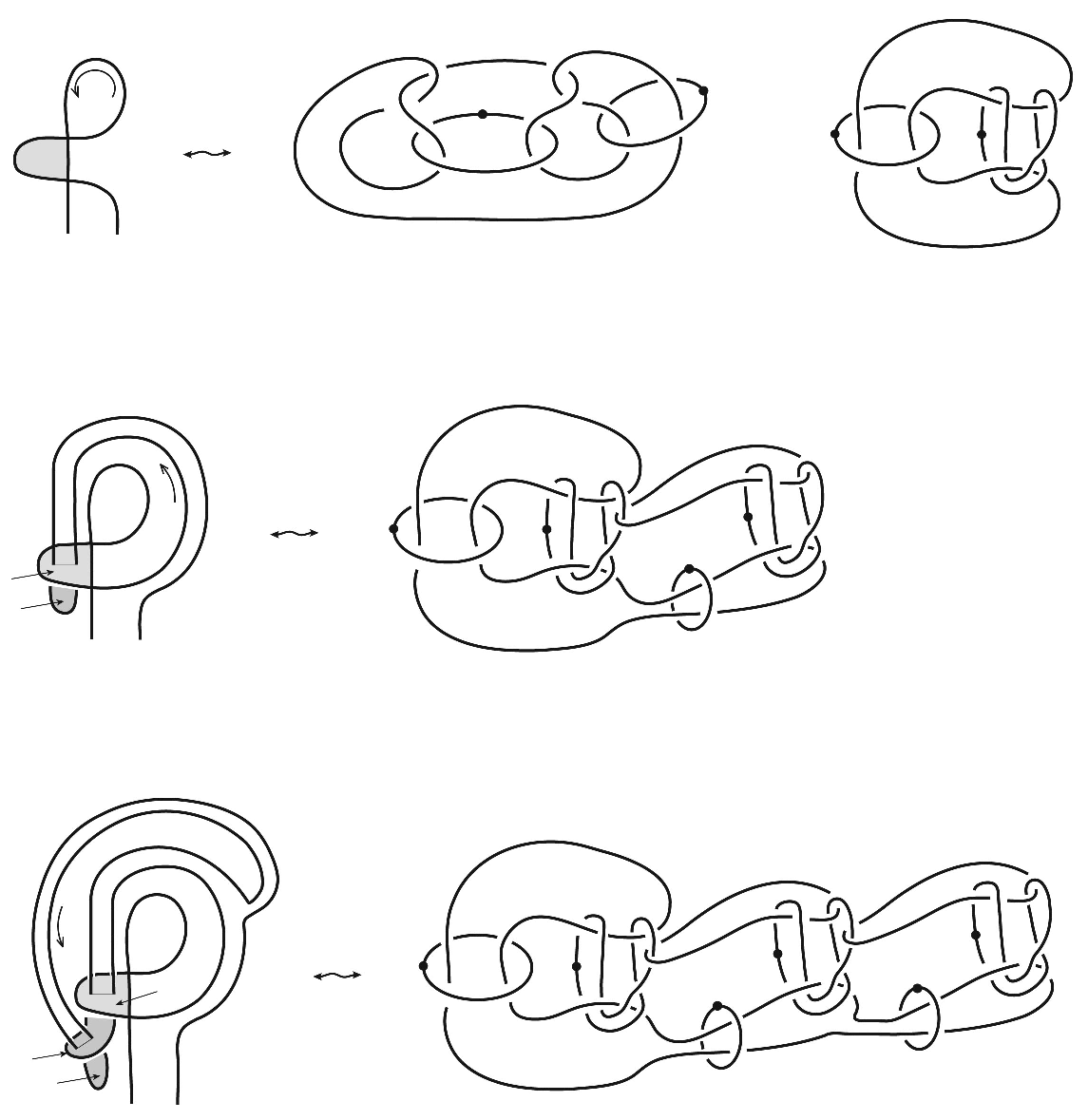}
{\scriptsize
\put(-347,63){$y$}
\put(-373,17){$W_2$}
\put(-126,78){$0$}
\put(-360,3){$W_3$}
\put(-317,42){$W_1$}
\put(-80,46){$W_2$}
\put(-125,10){$A_2$}
\put(-236,37){$A_1$}
\put(-149,39){$W_1$}
\put(-70,78){$0$}
\put(-11,60){$W_3$}
\put(-56,17){$A_3$}
}
{\scriptsize
\put(-318,213){$y$}
\put(-380,180){$W_1$}
\put(-130,230){$0$}
\put(-376,167){$W_2$}
\put(-87,202){$W_2$}
\put(-132,159){$A_2$}
\put(-246,186){$A_1$}
\put(-158,188){$W_1$}
}
{\scriptsize
\put(-341,349){$y$}
\put(-360,324){$W$}
\put(-345,337){$A$}
\put(-10,330){$W$}
\put(-94,323){$A$}
\put(-131,336){$W$}
\put(-211,346){$A$}
}
\put(-116,323){$\cong$}
\put(-394,270){(a) Result of a finger move on a disk, with Whitney disk $W$ and accessory disk $A$.}
\put(-394,133){(b) Secondary finger move}
\put(-394,-19){(c) 3rd order finger move}
\caption{Iterated finger moves giving raise to higher Engel relations: a spine version on the left and a precise description using Kirby diagrams on the right.}
\label{Constructive Engel}
\end{figure}

It is a standard technique in link homotopy theory to perform finger moves on surfaces in $4$-space in order to introduce Milnor relations (\ref{eq:Milnor group}) in the fundamental group of the complement. We conclude this section by pointing out that there is an analogous way of introducing higher Engel relations using iterated finger moves, see figure \ref{Constructive Engel}.

\begin{thebibliography}{10}


\bibitem[Bur02]{Burnside} W. Burnside, {\it On groups in which every two conjugate operations are permutable},
Proc. Lond. Math. Soc. 35 (1902), 28-37.


\bibitem[Cap71]{Cappell} S. Cappell, 
{\it A splitting theorem for manifolds and surgery groups}, 
Bull. Amer. Math. Soc. 77 (1971), 281-286. 


\bibitem[Cas73]{Casson} A. Casson, {\it Three lectures on new infinite constructions in $4$-dimensional manifolds}, 
With an appendix by L. Siebenmann. Progr. Math. 62, \`{A} la recherche de la topologie perdue, 
201-244, Birkh\"{a}user Boston, Boston, MA, 1986.


\bibitem[CF84]{CF} A. Casson and M.  Freedman, {\it Atomic surgery problems}, Four-manifold theory (Durham, N.H., 1982), 181-199, Contemp. Math., 35, Amer. Math. Soc., Providence, RI, 1984.


\bibitem[CP14]{CP} J.C. Cha and M. Powell, {\it Casson towers and slice links}, 
arXiv:1411.1621.

\bibitem[Co90]{Cochran}  T.D. Cochran, {\it Derivatives of links: Milnor's concordance invariants and Massey's products},  Mem. Amer. Math. Soc. 84 (1990), no. 427.



\bibitem[Fr82a]{F0} M. Freedman, {\it The topology of four-dimensional
manifolds}, J. Differential Geom. 17(1982), 357-453.


\bibitem[Fr82b]{F01} M. Freedman,  {\it A surgery sequence in dimension four; the relations with knot concordance}, Invent. Math. 68 (1982),
195-226.


\bibitem[Fr83]{F1} M. Freedman, {\it The disk theorem for
four-dimensional manifolds}, Proc. ICM Warsaw (1983), 647-663.


\bibitem[Fr86a]{F2} M. Freedman, {\it A geometric reformulation
of four dimensional surgery}, Topology Appl., 24 (1986), 133-141.


\bibitem[Fr86b]{F3} M. Freedman, {\it Are the Borromean rings
$(A,B)$-slice?}, Topology Appl., 24 (1986), 143-145.


\bibitem[FK12]{FK} M. Freedman and V. Krushkal, {\it Topological arbiters}, J. Topol. 5 (2012), 226-247.


\bibitem[FL89]{FL} M. Freedman and X.S. Lin, {\it On the $(A,B)$-slice
problem}, Topology Vol. 28 (1989), 91-110.


\bibitem[FQ90]{FQ} M. Freedman and F. Quinn, {\it The topology of
4-manifolds}, Princeton Math. Series 39, Princeton, NJ, 1990.


\bibitem[FT95a]{FT} M. Freedman and P. Teichner,
{\it $4$-Manifold Topology I: Subexponential groups}, Invent. Math. 122 (1995), 509-529.


\bibitem[FT95b]{FT2} M. Freedman and P. Teichner,
{\it $4$-Manifold Topology II: Dwyer's filtration and surgery kernels}, Invent. Math. 122 (1995), 531-557.

\bibitem[Gif79]{Giffen} C. Giffen, {\it Link concordance implies link homotopy}, Math. Scand., 45 (1979), 243-254.


\bibitem[Gol79]{Goldsmith} D. Goldsmith, {\it Concordance implies homotopy for classical links in $S^3$}, Comment.
Math. Helvitici, 54 (1979), 347-355.


\bibitem[GS99]{GS} R.E. Gompf and A.I. Stipsicz,
$4$-manifolds and Kirby calculus.
Graduate Studies in Mathematics, 20. American Mathematical Society, Providence, RI, 1999.


\bibitem[HST03]{HT} J. Hass, J. Snoeyink and W.P. Thurston,
{\it The size of spanning disks for polygonal curves}, 
Discrete Comput. Geom. 29 (2003), no. 1, 1-17. 


\bibitem[HVL05]{HVL} G. Havas and M.R. Vaughan-Lee,
{\it $4$-Engel groups are locally nilpotent}, 
Internat. J. Algebra Comput. 15 (2005), 649-682. 


\bibitem[Hei61]{Hei} H. Heineken, 
{\it Engelsche Elemente der L\'{a}nge drei}, 
Illinois J. Math. 5 (1961), 681-707. 


\bibitem[Ho29]{Hopkins} C. Hopkins, {\it Finite groups in which conjuate operations are commutative}, Am. J.
Math. 51, (1929), 35-41.


\bibitem[KS77]{KS} R. Kirby and L. Siebenmann, Foundational essays on topological manifolds, smoothings, and triangulations, Ann. Math. Stud. 88, Princeton U. Press 1977.


\bibitem[Kr98]{K2} V. Krushkal, {\it Additivity properties of Milnor's $\bar\mu$-invariants}, J. Knot Theory Ramifications 7 (1998), 625-637.


\bibitem[Kr08]{K} V. Krushkal, {\it A counterexample to the strong version of Freedman's conjecture}, Ann. of Math. 168 (2008), 675-693.


\bibitem[Kr13]{K1} V. Krushkal, {\it ``Slicing'' the Hopf link},  Geom. Topol. 19 (2015), 1657-1683.


\bibitem[KQ00]{KQ} V. Krushkal and F. Quinn, {\it Subexponential groups in
$4$-manifold topology}, Geom. Topol. 4 (2000), 407-430.


\bibitem[Lev42]{Levi} F.W. Levi,
{\it Groups in which the commutator operation satisfies certain algebraic conditions}, 
J. Indian Math. Soc. 6 (1942), 87-97.


\bibitem[MKS66]{MKS} W. Magnus, A. Karrass and D. Solitar, Combinatorial group theory: Presentations of groups in terms of generators and relations, Interscience Publishers, New York-London-Sydney 1966.


\bibitem[Mil54]{M} J. Milnor, {\it Link Groups}, Ann. Math 59 (1954), 177-195.


\bibitem[Mil57]{M2} J. Milnor, {\it Isotopy of links}, Algebraic geometry and topology, Princeton Univ. Press, 1957, 280-306.


\bibitem[Sta65]{Stallings} J. Stallings, {\it Whitehead torsion of free products}, Ann. of Math. 82 (1965), 354-363.


\bibitem[Tra11]{Traustason} G. Traustason,  {\it Engel groups},  Groups St Andrews 2009 in Bath. Volume 2, 520-550, London Math. Soc. Lecture Note Ser., 388, Cambridge Univ. Press, Cambridge, 2011.


\bibitem[Wal67]{Wall0} C.T.C. Wall, {\it Poincar\'{e} complexes. I.}, Ann. of Math. 86 (1967), 213-245. 


\bibitem[Wal70]{Wall} C.T.C. Wall, Surgery on compact manifolds. London Mathematical Society Monographs, No. 1. Academic Press, London-New York, 1970.


\end{thebibliography}
\end{document}